\documentclass[12pt]{amsart}
\usepackage{amsmath,amssymb,amscd}
\usepackage{epsfig, euler}

\numberwithin{equation}{section}
\newtheorem{thm}{Theorem}[section]
\newtheorem{prop}[thm]{Proposition}
\newtheorem{coro}[thm]{Corollary}
\newtheorem{lemma}[thm]{Lemma}

\theoremstyle{definition}

\newtheorem{rem}[thm]{Remark}
\newtheorem{eg}[thm]{Example}
\newtheorem{ex}[thm]{Exercise}
\newtheorem{defn}[thm]{Definition}

\def\til{\widetilde}
\def\mbb{\mathbb}
\def\mcl{\mathcal}
\def\mfk{\mathfrak}
\def\ten{\otimes}
\def\ex{\times}

\def\tu{\textup}

\def\ndt{\noindent}

\def\a{\alpha}
\def\b{\beta}
\def\d{\delta}

\def\e{\epsilon}

\def\Gm{\Gamma}

\def\sm{\sigma}

\def\o{\omega}

\def\bA{\mbb A}
\def\A{\mbb A}

\def\bG{\mbb G}
\def\bR{\mbb R}

\def\bP{\mbb P}
\def\bZ{\mbb Z}

\def\bQ{\mbb Q}
\def\cC{\mcl C}
\def\cD{\mcl D}
\def\cO{\mcl O}
\def\cU{\mcl U}
\def\cX{\mcl X}

\def\inj{\hookrightarrow}

\def\SS{\mcl{O}}
\def\deg{\textup{deg} \, }

\def\lra{\longrightarrow}

\def\inv{^{-1}}

\def\spec{\tu{Spec\,}}

\def\ord{\tu{ord\,}}
\def\Proj{\tu{\textbf{Proj}\,}}

\def\chow{Chow}

\def\Hgn{Hilb_{g,n}}
\def\Cgn{Chow_{g,n}}
\def\Hg2{Hilb_{g,2}}
\def\Cg2{Chow_{g,2}}

\def\SL{\textup{SL}}
\def\om2{\omega^{\ten 2}}
\def\Gr{\tu{Gr}}
\def\Mg{\bar{M}_g}

\def\FMps{\overline{{\bf \mcl M}}^{ps}_g}
\def\Mps{ \overline{M}^{ps}_g}
\def\FMg{\overline{{\bf \mcl M}}_g}
\def\FCg{\overline{\mcl C_g}}

\def\ds{\oplus}

\def\inj{\hookrightarrow}
\def\GL{\tu{GL}}

\def\cycle{\varpi}

\def\Sym{\tu{Sym}}
\def\Hom{\tu{Hom}}

\def\dps{\d^{ps}}
\def\Mcs{\Mg^{cs}}
\def\Mhs{\Mg^{hs}}

\def\dhs{\d^{hs}}
\def\dcs{\d^{cs}}

\def\k{k}

\def\ra{\rightarrow}

\def\bar{\overline}

\input xy
\xyoption{all}

\input epsf
\begin{document}

\normalsize


\title[Log minimal model program for $\Mg$: first flip]{Log minimal model program for the moduli space
of stable curves: the first flip}
\author{Brendan Hassett and Donghoon Hyeon}

\date{\today}

\begin{abstract} We give a geometric invariant theory (GIT) construction of the log
canonical model $\Mg(\a)$ of the pairs $(\Mg, \a \d)$ for $\a \in
(7/10 - \e, 7/10]$ for small $\e \in \bQ_+$.  We show that
$\Mg(7/10)$ is isomorphic to the GIT quotient of the Chow variety
bicanonical curves;  $\Mg(7/10-\e)$ is isomorphic to the GIT
quotient of the asymptotically-linearized Hilbert scheme of
bicanonical curves.  In each case, we completely classify the
(semi)stable curves and their orbit closures.  Chow
semistable curves have ordinary cusps and tacnodes as singularities but
do not admit elliptic tails.  Hilbert semistable curves satisfy further conditions, e.g.,
they do not contain elliptic bridges.
We show that there is a small contraction $\Psi:
\Mg(7/10+\e) \to \Mg(7/10)$ that contracts the locus of elliptic
bridges.  Moreover, by using the GIT interpretation of the log
canonical models, we construct a small
contraction $\Psi^+ : \Mg(7/10-\e) \to \Mg(7/10)$ that is the Mori
flip of $\Psi$.
\end{abstract}

\maketitle

\tableofcontents



\section{Introduction}\label{S:intro}
Our inspiration is to understand the canonical model
of the moduli space $\Mg$ of stable curves of genus $g$.
This is known to be of general type for
$g=22$ and $g\ge 24$ \cite{Far,HM,EH}.  In these cases, we can
consider the canonical ring
\[
\ds_{n\geq 0} \Gamma(n (K_{\FMg})).
\]
which is finitely generated by a fundamental conjecture of birational geometry,
recently proven in \cite{BCHM}.
Then the corresponding projective variety
\[
Proj \ds_{n\geq 0} \Gamma(n (K_{\FMg}))
\]
is birational to $\Mg$ and is called its
{\em canonical model}.

There has been significant recent progress in understanding
canonical models of moduli spaces.  For moduli spaces ${\mathcal
A}_g$ of principally polarized abelian varieties of dimension
$g\ge 12$, the canonical model exists and is equal to the first
Voronoi compactification \cite{SB}. Unfortunately, no analogous
results are known for $\Mg$, even for $g\gg 0$.

Our approach is to approximate the canonical models
with {\em log canonical models}.  Consider
$\a \in [0,1]\cap {\mathbb Q}$ so that $K_{\FMg}+\a\d$
is an effective ${\mathbb Q}$-divisor.  We have the graded ring
\[
\ds_{n\geq 0} \Gamma(n (K_{\FMg} +
\a\d))
\]
and the resulting projective variety
$$
\Mg(\a):= Proj \left(\ds_{n\geq 0} \Gamma(n (K_{\FMg} +
\a\d))\right).
$$

Our previous paper \cite{HH1} describes $\Mg(\a)$ explicitly for
large values of $\a$.  For simplicity we assume that $g\ge 4$: 
Small genera cases have been considered in  \cite{Has,HL1,HL2}.
For $9/11<\a \le 1$, $\Mg(\a)$ is equal to $\Mg$. The first
critical value is $\a=9/11$: $\Mg(9/11)$ is the coarse moduli
space of the moduli stack $\FMps$ of pseudostable curves
\cite{Sch}. A pseudostable curve may have cusps but they are not
allowed to have {\em elliptic tails}, i.e., genus one subcurves
meeting the rest of the curve in one point. Furthermore, there is
a divisorial contraction
$$T:\Mg \ra \Mg(9/11)$$
induced by the morphism
$\mcl T : \FMg \to \FMps$
of moduli stacks which replaces an elliptic tail with a
cusp.  Furthermore,
$\Mg(\a) \simeq \Mg(9/11)$ provided $7/10<\a\le 9/11$.

This paper addresses what happens when $\a=7/10$.  Given
a sufficiently small positive $\e \in \mathbb Q$,
we construct
a small contraction and its flip:
\[
\xymatrix{\Mg(\frac7{10}+\e) \ar[dr]^-{\Psi} & &
\Mg(\frac7{10}-\e)
\ar[dl]_-{\Psi^+}\\
& \Mg(\frac7{10}) &\\
}
\]
The resulting spaces arise naturally as geometric invariant theory
(GIT) quotients and admit partial modular descriptions. We
construct $\Mg(7/10)$ as a GIT quotient of the Chow variety of
bicanonical curves;  it parametrizes equivalence classes of {\it
c-semistable curves}.  We defer the formal definition, but these have
nodes, cusps, and tacnodes as singularities. The flip
$\Mg(7/10-\e)$ is a GIT quotient of the Hilbert scheme of
bicanonical curves;  it parametrizes equivalence classes of {\it
h-semistable curves},  which are c-semistable curves not
admitting certain subcurves composed of elliptic curves (see Definition~\ref{D:h-stable}).

\begin{figure}[!htb]
  \begin{center}
    \includegraphics[width=5.5in]{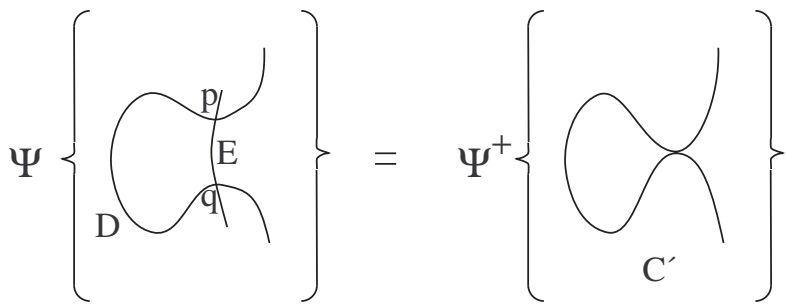}
  \end{center}
\caption{ Geometry of the flip} \label{fig:GF}
\end{figure}

We may express the flip in geometric terms (Figure~\ref{fig:GF}):  Let
$C=D\cup_{p,q}E$ denote an {\em elliptic bridge}, where $D$
is smooth of genus $g-2$, $E$ is smooth of genus one,
and $D$ meets $E$ at two nodes $p$ and $q$.  Let $C'$
be a tacnodal curve of genus $g$, with normalization $D$
and conductor $\{p,q\}$.
In passing from $\Mg(7/10+\e)$ to
$\Mg(7/10-\e)$, we replace $C$ with $C'$.
Note that the descent data for
$C'$ includes the choice of an isomorphism of tangent spaces
$$\iota:T_p D \stackrel{\sim}{\lra} T_q D;$$
the collection of such identifications is a principal homogeneous
space for ${\mathbb G}_m$.  When $C$ is a generic elliptic
bridge, the fiber $(\Psi^+)^{-1}(\Psi(C))\simeq \bP^1$;
see Proposition~\ref{P:const-tac} for an explicit interpretation
of the endpoints.

\smallskip

Here we offer a brief summary of the contents of this paper;  a more
detailed roadmap can be found in Section \ref{SS:detailed}.
Section \ref{S:GITCH} is devoted to a general discussion of the GIT of
Chow points and Hilbert points.  The main applications are the
analysis of tautological classes and polarizations on the Hilbert
scheme, the resulting formulas for Hilbert-Mumford indices
(Proposition~\ref{prop:curveresult}), and cycle maps
(Corollary~\ref{coro:HilbtoChow}).  We also recall various
formulations of the Hilbert-Mumford one-parameter-subgroup
criterion.

Section \ref{S:basin} is a brief review of the basin-of-attraction
techniques used in this paper.  These are important for analyzing when
distinct curves are identified in the GIT quotients.

Section \ref{S:CMS} discusses, in general terms, how to obtain contractions
of the moduli space of stable curves from GIT quotients of Hilbert
schemes.  The resulting models of the moduli space depend
on the choice of linearization;  we express the polarizations in
terms of tautological classes.

Section \ref{S:propchstable} summarizes basic properties of c-semistable
curves: embedding theorems and descent results for tacnodal curves.
Section \ref{S:UBC} offers a preliminary analysis of the GIT of
the Hilbert scheme and the Chow variety of bicanonically embedded curves
of genus $g \ge 4$.   Then in Section~\ref{S:automorphisms} we
enumerate the curves with positive-dimensional automorphism groups.
Section \ref{S:FGITQ} applies this to give a GIT construction of the flip
$\Psi^+ : \Mg(7/10-\e) \to \Mg(7/10)$.

Section \ref{S:hardwork} offers a detailed orbit closure analysis,
using basins of attractions and a careful analysis of the action of
the automorphism group on tangent spaces.  The main
application is a precise
description of the semistable and stable bicanonical curves,
proven in Section
\ref{S:PSA}.

Throughout, we work over an algebraically closed field $k$,
generally of characteristic zero.  However, Sections~\ref{S:GITCH}
and \ref{S:propchstable} are valid in positive characteristic.

\smallskip

\noindent {\bf Acknowledgments:} The first author was partially
supported by National Science Foundation grants 0196187, 0134259,
and 0554491,
the Sloan Foundation,
and the Institute of Mathematical Sciences of the
Chinese University of Hong Kong.
The second author was partially supported by the Korea Institute for
Advanced Study and
the Korea Science and Engineering Foundation (KOSEF) grant funded by the Korea government (MOST) (No. R01-2007-000-10948-0).
 We owe a great deal to S. Keel, who helped
shape our understanding of the birational geometry of $\Mg$
through detailed correspondence.
We are also grateful to D. Abramovich, Y. Kawamata, I. Morrison,
B.P. Purnaprajna,
M. Simpson, D. Smyth, and D. Swinarski for useful conversations.

\section{Statement of results and strategy of proof}
\label{S:SRSP}
\subsection{Stability notions for algebraic curves}
In this paper, we will use four stability conditions:
Deligne-Mumford stability \cite{DM}, Schubert pseudostability \cite{Sch},
c-(semi)stability, and h-(semi)stability.
We recall the definition of {\it pseudostability}, which
is obtained from Deligne-Mumford stability
by allowing ordinary cusps and prohibiting elliptic tails:

\begin{defn}\label{D:p-stable}\cite{Sch}
A complete curve is
{\em pseudostable} if
\begin{enumerate}
\item it is connected, reduced, and
has only nodes and ordinary cusps as singularities;
 \item admits no {\em elliptic tails}, i.e.,
connected subcurves of arithmetic genus one meeting the rest of
the curve in one node;
\item the canonical sheaf of the curve is ample.
\end{enumerate}
The last condition means that each subcurve of genus zero meets
the rest of the curve in at least three points.
\end{defn}

Before formulating the notions of c- and h-(semi)stability,
we need the following definition:

\begin{defn} \label{D:eb}  An {\em elliptic bridge }
is a connected subcurve of arithmetic genus one
meeting the rest of the curve in two nodes.
\end{defn}
\begin{figure}[htb]
\begin{center}
\includegraphics[width=5in]{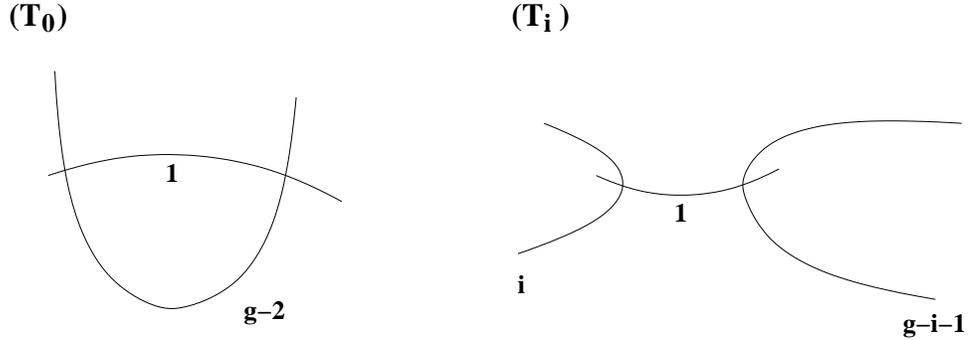}
\end{center}
\caption{Generic elliptic bridges}\label{F:eb}
\end{figure}

In our stability analysis, we will require additional technical definitions:

\begin{defn}\label{D:ech} An {\em open elliptic chain} of length $r$
is a two-pointed
projective curve $(C', p, q)$ such that
\begin{itemize}
\item $C' = E_1 \cup_{a_1}\dots\cup_{a_{r-1}} E_r$ where each $E_i$ is connected of genus one, with nodes, cusps or tacnodes as singularities;
\item $E_i$ intersects $E_{i+1}$ at a single
tacnode $a_i$, for $i = 1, \dots, r-1$;
\item $E_i\cap E_j = \emptyset$ if $|i-j| > 1$;
\item $p, q \in C'$ are smooth points with $p \in E_1$ and $q\in E_r$;
\item $\omega_{C'}(p+q)$ is ample.
\end{itemize}
\end{defn}
An open elliptic chain of length $r$ has arithmetic genus $2r-1$.

\begin{figure}[!htb]
  \begin{center}
    \includegraphics[width=4.5in]{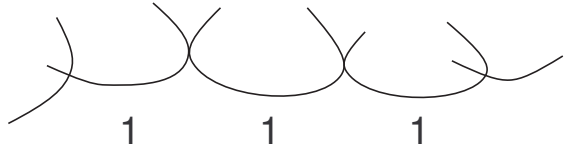}
  \end{center}

  \caption{Generic elliptic chain of length three}
\end{figure}

\begin{defn} \label{D:echain} Let $C$ be a projective connected curve of arithmetic
genus $g \ge 3$, with nodes, cusps, and tacnodes as singularities.
We say $C$ {\em admits an open elliptic chain} if there is an open elliptic chain $(C',p,q)$ and a morphism $\iota : C' \to C$ such that
\begin{itemize}
\item $\iota$ is an isomorphism over $C'\setminus \{p,q\}$ onto its image;
\item $\iota(p), \iota(q)$ are nodes of $C$;
we allow the case $\iota(p) = \iota(q)$, in which case $C$ is said to be a {\em closed elliptic chain}.
\end{itemize}
$C$ admits a {\em weak elliptic chain} if there exists $\iota : C' \to C$ as above with the second condition replaced by
\begin{itemize}
\item $\iota(p)$ is a tacnode of $C$ and $\iota(q)$ is a node of $C$; or
\item $\iota(p) = \iota(q)$ is a tacnode of $C$, in which case
$C$ is said to be a {\em closed weak elliptic chain}.
\end{itemize}
\end{defn}

\begin{figure}[htb]
\centerline{\scalebox{0.6}{\psfig{figure=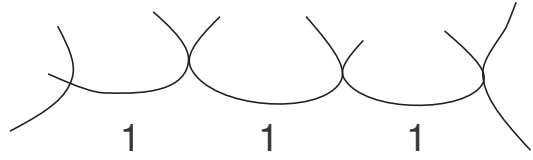}}}
\caption{Generic weak elliptic chain}\label{F:w-echain}
\end{figure}

Now we are in a position to formulate our main stability notions:
\begin{defn}\label{D:c-stable}
 A complete curve $C$ is said to be {\it
c-semistable} if
\begin{enumerate}
\item $C$ has nodes, cusps and tacnodes as singularities;
\item $\o_C$ is ample;
\item a connected genus one subcurve meets the rest of the
curve in at least two points (not counting multiplicity).
\end{enumerate}
It is said to be {\it c-stable} if it is c-semistable and has no tacnodes or elliptic bridges.
\end{defn}
\begin{defn}\label{D:h-stable}
A complete curve $C$ of genus $g$ is said to be {\it
h-semistable} if it is c-semistable and admits no elliptic
chains. It is said to be {\it h-stable} if it is h-semistable and
admits no weak elliptic chains.

\begin{rem} \label{R:whatiscstable}
A curve is c-stable if and only if it is pseudostable and has no elliptic bridges.
\end{rem}

\end{defn}
Table~\ref{Tab:stability} summarizes the defining
characteristics of the stability notions.

\

\scriptsize

\begin{table}[h]
\caption{Stability notions}\label{Tab:stability}
\begin{center}
\begin{tabular}{|c||c|c|c|c|c|}

\hline
& singularity & $\begin{array}{c} \mbox{genus zero} \\ \mbox{subcurve}\\
\mbox{meets}\\   \mbox{the rest in ...}\end{array}$ &  $\begin{array}{c}
    \mbox{genus one} \\ \mbox{subcurve}\\ \mbox{meets} \\
   \mbox{the rest in ...}\end{array}$ &
    $\begin{array}{c} \mbox{elliptic} \\
   \mbox{chain} \\  \end{array}$ & $\begin{array}{c} \mbox{weak elliptic} \\
   \mbox{chain}  \end{array}$\\

\hline

stable & nodes & $\ge 3$ points & -- & -- & -- \\
\hline

   pseudostable & $\begin{array}{c} \mbox{nodes,} \\
   \mbox{cusps} \\  \end{array}$ & $\ge 3$ points & $\ge 2$ points & -- & -- \\

   \hline

   c-semistable & $\begin{array}{c} \mbox{nodes,} \\ \mbox{cusps,} \\
   \mbox{tacnodes}\end{array}$ &  $\begin{array}{c}\mbox{$\ge 3$ points}\\
   \mbox{counting}\\ \mbox{ multipilicity}\end{array}$ &  $\ge 2$ points & -- & --\\

   \hline

   c-stable & $\begin{array}{c} \mbox{nodes,} \\ \mbox{ cusps}\\
   \end{array}$ &  $\begin{array}{c}\mbox{$\ge 3$ points}\\
  \end{array}$ &  $\begin{array}{c}\mbox{$\ge 3$ points}\\
  \end{array}$ & -- & --\\

   \hline

   h-semistable &  $\begin{array}{c} \mbox{nodes,}\\ \mbox{cusps,}\\
   \mbox{tacnodes}\end{array}$  &  $\begin{array}{c}\mbox{$\ge 3$ points}\\
   \mbox{counting}\\ \mbox{ multipilicity}\end{array}$ &
    $\begin{array}{c}\mbox{$\ge 3$ points}\\
   \mbox{counting}\\ \mbox{ multipilicity}\end{array}$  & $\begin{array}{c} \mbox{not} \\
   \mbox{admitted} \\  \end{array}$ & --\\

   \hline

    h-stable &  $\begin{array}{c} \mbox{nodes} \\ \mbox{cusps,}\\
   \mbox{tacnodes}\end{array}$  &  $\begin{array}{c}\mbox{$\ge 3$ points}\\
   \mbox{counting}\\ \mbox{ multipilicity}\end{array}$ &
    $\begin{array}{c}\mbox{$\ge 3$ points}\\
   \mbox{counting}\\ \mbox{ multipilicity}\end{array}$  & $\begin{array}{c} \mbox{not} \\
   \mbox{admitted} \\  \end{array}$ & $\begin{array}{c} \mbox{not} \\
   \mbox{admitted} \\  \end{array}$ \\

   \hline

\end{tabular}
\end{center}
\label{default}
\end{table}%

\normalsize

\subsection{Construction of the small contraction $\Psi$}
\label{SS:CSCP}

We start with some preliminary results.
Recall from \cite{HH1} that $\mcl T:\FMg \ra \FMps$ denotes the functorial contraction and
$T:\Mg \ra \Mps=\Mg(9/11)$ the induced morphism on coarse moduli
spaces, which contracts the divisor $\Delta_1$.

\begin{lemma} \label{L:logterm}
For $\alpha  < 9/11$,
$(\FMps,\alpha\delta^{ps})$
and
$(\Mps,\alpha
\Delta^{ps})$
are log terminal and
$$
\Mg(\a) \simeq  Proj \left(\ds_{n\geq 0} \Gamma(n (K_{\Mps} +
\a\dps))\right).
$$
\end{lemma}
\begin{proof}
Since $g>3$, the locus in $\Mps$ parametrizing curves with nontrivial
automorphisms has codimension $\ge 2$ \cite[\S 2]{HM}.  (Of course,
we have already collapsed $\delta_1$.)  Thus the
coarse moduli map $q:\FMps \ra \Mps$ is unramified in codimension one
and
\begin{equation} \label{E:ramify}
q^*(K_{\Mps}+\alpha \Delta^{ps})=K_{\FMps}+\alpha \delta^{ps}
\end{equation}
for each $\alpha$.
We have the log discrepancy equation \cite[\S 4]{HH1}
\begin{equation} \label{E:logdisc}
K_{\FMg}+\alpha \delta ={\mcl T}^*(K_{\FMps}+\alpha \delta^{ps})+
(9-11\alpha)\delta_1
\end{equation}
and  the pull back
$${\mcl T}^*(K_{\Mps}+ 7/10 \delta^{ps})=
K_{\FMg}+7/10\delta-13/10\delta_1\sim 10\lambda - \delta-\delta_1,$$
where $\sim$ designates proportionality.

Since $\FMg$ is smooth and $\delta$ is normal crossings, the pair
$(\FMg,\alpha\delta+(11\alpha -9)\delta_1)$ is log terminal. The
discrepancy equation implies that $(\FMps,\alpha\delta^{ps})$ is
log terminal for $\alpha \in [7/10,9/11)$.  Applying the
ramification formula \cite[5.20]{KM} to (\ref{E:ramify}) (or
simply applying \cite[A.13]{HH1}), we find that $(\Mps,\alpha
\Delta^{ps})$ is also log terminal.

Since $\Delta_1$
is $T$-exceptional, for each Cartier divisor $L$ on $\Mps$ and $m \ge 0$
we have $ \Gamma(\Mg, T^*L+m\Delta_1) \simeq \Gamma(\Mg^{ps},
L)$. This implies that
\[
\begin{array}{ccl}
\Mg(\a) & = & Proj \oplus_{n\ge 0} \Gamma(\Mg, n( K_{\FMg} + \a \d) ) \\
& = & Proj \oplus_{n\ge 0} \Gamma\left(\Mg, n( T^*(K_{\Mps} + \a \d^{ps})+ (9-11\a) \d_1 )\right) \\
& \simeq &  Proj \oplus_{n\ge 0} \Gamma\left(\Mps, n(K_{\Mps} + \a \d^{ps}) \right).
\end{array}
\]
\end{proof}

We shall construct the contractions 
by using the powerful results of \cite{GKM}:
\begin{prop} \label{P:makePsi}
For $\alpha \in (7/10,9/11]\cap \bQ$,
there there exists a
birational contraction
$$\Psi:\Mg(\alpha) \ra \Mg(7/10).$$
It contracts the codimension-two strata
$T_i, i=0,2,\ldots, \lfloor (g-1)/2 \rfloor,$ where
\begin{enumerate}
\item
$ T_0 = \{ E \cup_{p,q} D \, | \, g(E) = 1, g(D) = g-2\}$;

\smallskip

\item
$T_i = \{ C_1 \cup_p E \cup_q C_2 \, | \, g(C_1) = i, g(E) = 1,
g(C_2) = g-1-i \}$, \, \, $2\le i \le \lfloor (g-1)/2 \rfloor$,
\end{enumerate}
by collapsing the loci $\overline{M}_{1,2} \subset T_i $ corresponding to
varying $(E,p,q)$.
\end{prop}
\begin{rem}
We shall see in Corollary~\ref{C:excep-Psi} that $\Psi$ is an
isomorphism away from  $T_{\bullet}:=\cup T_i$.
\end{rem}
\begin{proof}
Recall that
$K_{\Mps}+\alpha \Delta^{ps}$
is ample provided $7/10<\alpha \le 9/11$;  this
is part of the assertion that
$\Mg(\alpha)=\Mps$ for $7/10<\alpha \le 9/11$
\cite[Theorem 1.2]{HH1}.
However,
$K_{\FMps}+ 7/10 \Delta^{ps}$ is nef but not ample \cite[\S 4]{HH1}.
Indeed, the pull-back to $\Mg$
$$10\lambda - \delta-\delta_1$$
can be analyzed
using the classification of one-dimensional
boundary strata by Faber \cite{F4} and Gibney-Keel-Morrison \cite{GKM}.
It is `F-nef', in the sense that it intersects all these
strata nonnegatively,
and is therefore nef by \cite[6.1]{GKM}.
Later on, we will list the strata meeting it with degree zero.

We apply Kawamata basepoint freeness \cite[3.3]{KM}:
\begin{quote}
Let $(X, D)$ be a proper Kawamata log terminal pair with
$D$ effective. Let $M$ be a nef Cartier divisor such that $a M - K_X -
D$ is nef and big for some $a > 0$. Then $|b M|$ has no basepoint
for all $b \gg 0$.
\end{quote}
For our application, $M$ is a Cartier multiple of
$K_{\Mps}+7/10\Delta^{ps}$ and $D=(7/10-\epsilon) \Delta^{ps}$ for
small positive $\epsilon\in \bQ$. The resulting morphism is
denoted $\Psi$.

We claim that $\Psi$ is birational.
To establish the birationality, we show that each curve
$B\subset \Mg$ meeting the interior satisfies
$$B.(10\lambda - \delta - \delta_1)>0.$$
The Moriwaki divisor
$$A:=(8g+4)\lambda - g \delta_0 - \sum_{i=1}^{\lfloor g/2\rfloor}4i(g-i)\delta_i$$
meets each such curve nonnegatively \cite[Theorem B]{Mor}. We can
write
$$10\lambda -\delta  -\delta_1=(1/g) \, A +(2-4/g)\lambda+(2-4/g)\delta_1+
    \sum_{i=2}^{\lfloor g/2\rfloor} (-1+ 4i(g-i)/g)\delta_i.$$
Each of these coefficients is positive:  Clearly
$1/g,2-4/g>0$ and since $2i/g \le 1$, 
$$-1 + 4i(g-i)/g = -1 + 4i - (2i/g)2i \ge -1 + 4i - 2i > 0.$$
Thus we have
$$B.(10\lambda -\delta -\delta_1) \ge (2-4/g)\lambda.B >0,$$
where the last inequality reflects the fact that the Torelli morphism
is nonconstant along $B$.

We verify the image of $\Psi$ equals $\Mg(7/10)$.
The log discrepancy formula (\ref{E:logdisc}) implies
$$\mathrm{Image}(\Psi)=\Proj \oplus_{n\ge 0}
\Gamma(n(K_{\FMg}+7/10\delta-13/10\delta_1)).$$
However, since $\Delta_1$ is $(\Psi\circ T)$-exceptional
adding it does not change the space of global sections, whence
$$\mathrm{Image}(\Psi)=\Proj \oplus_{n\ge 0}
\Gamma(n(K_{\FMg}+7/10\delta))=\Mg(7/10).$$

Finally, we offer a preliminary analysis of the locus contracted by
$\Psi$.  The main ingredient is the enumeration of one-dimensional
boundary strata in \cite{GKM}  (see also \cite[\S 4]{HH1}). We list
the ones orthogonal to $10\lambda - \delta -\delta_1$; any stratum
swept out by these classes is necessarily contracted by $\Psi$. In
the second and third cases $X_0$ denotes a varying $4$-pointed curve
of genus zero parametrizing the stratum.
\begin{enumerate}
\item{Families of elliptic tails, which sweep out $\delta_1$
and correspond to the extremal ray contracted
by $\mcl T$.}
\item{Attach a $2$-pointed curve of genus $0$ and a
$2$-pointed curve $(D,p,q)$
of genus $g-2$ to $X_0$ and stabilize.  Contracting this
and the elliptic tail stratum collapses $T_0$ along the
$\overline{M}_{1,2}$'s corresponding to fixing $(D,p,q)$
and varying the other components.}
\item{Attach a $1$-pointed curve $(C_1,p)$ of genus $i>1$, a $1$-pointed
curve $(C_2,q)$ of genus $g-1-i>1$, and a $2$-pointed curve of genus
$0$ to $X_0$ and stabilize.  Contracting this and the elliptic tail
stratum collapses $T_i$ along the $\overline{M}_{1,2}$'s
corresponding to fixing $(C_1,p),(C_2,q)$ and varying the other
components.}
\end{enumerate}
Thus the codimension-two strata $T_0,T_2,\ldots,T_{\lfloor (g-1)/2\rfloor}$
are all contracted by $\Psi$.
\end{proof}

\subsection{Construction of the flip $\Psi^+$}
Consider the Chow variety of degree $4g-4$ curves of genus $g$ in
$\bP^{3g-4}$.  Let $\Cg2$ denote the closure of the bicanonically
embedded smooth curves of genus $g$.  Similarly, let $\Hg2$ denote
the closure of these curves in the Hilbert scheme.

\begin{prop} \label{P:cyclebicanonical}
The cycle class map
\begin{equation} \label{E:cycle}
\cycle:\Hg2 \ra \Cg2
\end{equation}
induces a morphism of GIT quotients
$$\Hg2^{ss}/\!\!/\SL_{3g-3} \ra \Cg2^{ss}/\!\!/\SL_{3g-3},$$
where the Hilbert scheme has the asymptotic linearization
introduced in \S\ref{subsect:CPHS}.
\end{prop}
This is a special case of Corollary~\ref{coro:HilbtoChow},
which applies quite generally to cycle-class maps from
Hilbert schemes for Chow varieties.  (See \S\ref{subsect:CPHS}
for background information on the cycle class map.)
Let $\Mhs$ and $\Mcs$ denote the resulting GIT quotients
$\Hg2^{ss}/\!\!/\SL_{3g-3}$ and
$\Cg2^{ss}/\!\!/\SL_{3g-3}$, and
\begin{equation}\label{E:psiplus}
\Psi^+ : \Mhs \to \Mcs
\end{equation}
the morphism of Proposition~\ref{P:cyclebicanonical}.

\begin{thm}\label{T:main12}  Let $\e \in \bQ$ be a small
positive number.   There exist isomorphisms
\begin{equation} \label{eqn:Tmain1}
\Mg(7/10) \simeq \Mcs
\end{equation}
and
\begin{equation} \label{eqn:Tmain2}
\Mg(7/10-\e) \simeq \Mhs
\end{equation}
such that the induced morphism
$$\Psi^+:\Mg(7/10-\e) \ra \Mg(7/10)$$
is the flip of $\Psi$.
\end{thm}
We thus obtain a modular/GIT interpretation of the flip:
\[
\xymatrix{ \Mg(\frac7{10}+\e) \simeq \Mps \ar[dr]^-{\Psi} & &
\Mg(\frac7{10}-\e) \simeq \Mg^{hs}
\ar[dl]_-{\Psi^+}\\
& \Mg(\frac7{10})\simeq \Mg^{cs} &\\
}
\]

\subsection{Stability results on bicanonical curves}
For c-semistable curves, $\omega_C^{\otimes 2}$ is very ample and has
no higher cohomology (Proposition~\ref{P:embedcstable}). The image
in $\bP^{3g-4}$ is said to be {\em bicanonically embedded}.

\begin{thm}\label{T:Chow} The semistable locus $\Cg2^{ss}$ (resp. stable locus $\Cg2^s$)
corresponds to bicanonically embedded c-semistable (resp. c-stable) curves.
\end{thm}
Unlike in $\Mg$ and $\Mg^{ps}$, nonisomorphic curves may be identified
in the quotient $Chow_{g,2}/\!\!/\SL_{3g-3}$.
For example, if a c-semistable curve $C = D \cup_{p,q} E$ consists of
a genus $g-2$ curve $D$ meeting in two nodes $p, q$ with an
elliptic curve $E$, then it is identified with any tacnodal curve obtained by
replacing $E$ with a tacnode.
In Section \ref{S:hardwork}, we shall give a complete
classification of strictly semistable curves and the curves in their
orbit closure.

\begin{thm}\label{T:Hilbert}
The semistable locus $\Hg2^{ss}$ (resp. stable locus $\Hg2^{s}$)
with respect to the asymptotic
linearization corresponds to bicanonically embedded h-semistable
(resp. h-stable) curves.
\end{thm}
One difference from the case of Chow points is that tacnodal curves may
well be Hilbert stable. For instance, when $g \ge 4$ {\em irreducible}
bicanonical h-semistable curves are necessarily h-stable.
When $g = 3$, a bicanonical h-semistable curve is
Hilbert strictly semistable if and only if it has a tacnode
\cite{HL2}. When $g = 4$, every h-semistable curve is h-stable
and the moduli functor is thus separated.

Since c-stable curves are h-stable (see Proposition~\ref{prop:HilbtoChow})
and pseudostable (see Remark~\ref{R:whatiscstable}), we have
\begin{coro} \label{C:excep-Psi}
$\Psi$ and $\Psi^+$ are isomorphisms over the locus of c-stable curves.
Thus $\Psi$ is a small contraction with exceptional
locus $T_{\bullet}$ and  $\Psi^+$ is a small
contraction with exceptional locus $\mathrm{Tac}$, the h-semistable curves
with tacnodes.
\end{coro}
Thus the geometry of the flip is as indicated in Figure~\ref{fig:GF}:
$\Psi^+(C')=\Psi(C)$ precisely when $C$ is the `pseudostable reduction'
of $C'$.

\subsection{Detailed roadmap for the GIT analysis}
\label{SS:detailed}

The proof of Theorems \ref{T:Chow} and \ref{T:Hilbert} is rather
intricate, so we give a bird's eye view for the reader's
convenience.

\

\noindent (1) The following implications are straightforward:
\begin{itemize}
\item From the definitions, it is clear that:
\[
\mbox{h-semistable $\Rightarrow$ c-semistable}
\]
\item  General results on linearizations of Chow and Hilbert schemes
(Proposition~\ref{prop:HilbtoChow}) imply
\[
\mbox{Hilbert semistable $\Rightarrow$ Chow semistable}
\]
and
\[
\mbox{Chow stable $\Rightarrow$ Hilbert stable.}
\]
\end{itemize}

\noindent (2) We next prove that non c-semistable (resp. non
h-semistable) curves are Chow unstable (resp. Hilbert unstable). The
main tool is the stability algorithm
(Proposition~\ref{prop:stabcrit2}).

\begin{itemize}
\item Non c-semistable curves can be easily destabilized by
one-parameter subgroups (\S \ref{S:UBC}). We obtain
\[
\mbox{Chow semistable $\Rightarrow$ c-semistable}.
\]
\item We show that if a curve $C$ admits an {\em open rosary}
of even length (see
Definition~\ref{D:rosary}),
then there is a 1-PS $\rho$ coming from the automorphism group of the rosary
such that the {\it  $m$th Hilbert point} $[C]_m$
(Definition~\ref{D:mth-hilb-pt}) is unstable with respect to $\rho$ for all $m\ge 2$ (Proposition~\ref{P:hs-o-ros} and  Proposition~\ref{P:hs-cr-1br}).

\item If $C$ admits an elliptic chain, then it is contained in the
{\em basin of attraction}
$A_\rho([C_0]_m)$ (see Definition~\ref{D:basin}) of a curve $C_0$ admitting an
open rosary of even length such that $\mu([C_0]_m, \rho) < 0$.  Hence
such curves are Hilbert unstable
(Propositions~\ref{P:hs-echain} and \ref{P:basin-cr-1br})
and we obtain
\[
\mbox{Hilbert semistable $\Rightarrow$ h-semistable}.
\]
\end{itemize}

\noindent (3) We prove ``c-semistable $\Rightarrow$ Chow
semistability", and use it to establish ``h-semistability
$\Rightarrow$ Hilbert semistability".
\begin{itemize}
\item The
only possible Chow semistable replacement of a c-stable curve
is itself (see Theorem~\ref{T:partial}).
Thus c-stable curves are Chow stable and hence Hilbert stable.

\item We show that any strictly c-semistable curve $C$ is contained
in a basin of attraction
of a distinguished c-semistable
curve $C^{\star}$ with one-parameter isomorphism
such that $\mu(Ch(C^{\star}), \rho) =
0$  (see Proposition~\ref{P:degenerate}).
Indeed, we choose $C^{\star}$ so that it has closed
orbit in the locus of c-semistable points
(cf. Proposition~\ref{P:c-minorbit}).

\item If $C$ is strictly c-semistable, its pseudo-stabilization
$D$ has elliptic
bridges. For any such $D$,
there is a distinguished strictly c-semistable curve $C^{\star}$ such
that its basins of attraction contain every c-semistable
replacement for $D$.
Futhermore, {\em every}
possible Chow-semistable replacement for $D$ is contained in some
basin of attraction $A_{\rho'}(Ch(C^{\star}))$ with
$\mu(Ch(C^{\star}), \rho') = 0$.  Since one of these must be
Chow semistable, every one is Chow semistable (see Lemma~\ref{L:flat limit}).

\item  The Hilbert semistable curves form a subset of the set of Chow semistable 
curves.
We first identify the Chow semistable curves admitting one-parameter subgroups
that are Hilbert-destabilizing.  Then we show that any curve that is Hilbert unstable
but Chow semistable arises in the basin of attraction of such a curve.
These basins of attraction consist of the curves that are c-semistable but not h-semistable.
Thus the h-semistable curves are Hilbert semistable
(\S\ref{SS:HS}).

\end{itemize}

\section{GIT of Chow varieties and Hilbert
schemes}\label{S:GITCH} Let $\bP^N=\bP(V)$ for some
$(N+1)$-dimensional vector space $V$. Throughout this section, let
$\rho : \bG_m \to \GL(V)$ be a one-parameter subgroup and $x_0, \dots,x_N$
be homogeneous coordinates that diagonalize the
$\rho$-action so that
\[
\rho(t).x_i = t^{r_i}x_i, \, i = 0, \dots, N, \, r_0 \geq \dots
\geq r_N=0.
\]
We have the associated one-parameter subgroup
$\rho^{\circ}:\bG_m \ra \SL(V)$
\[
\rho^{\circ}(t). x_i = t^{r_i-(r_0+\ldots+r_N)/(N+1)}x_i.
\]
Given $x\in \bP(V)$,  the {\em Hilbert-Mumford index} is given by
(cf. \cite[2.1]{GIT}):
$$\mu(x,\rho)=\max \{-r_i+(r_0+\ldots+r_N)/(N+1):x_i \neq 0 \}.$$
We say that  $x$ is (semi)stable with respect to $\rho^{\circ}$ if
 $\mu(x,\rho) >(\ge) \,  0$.  A fundamental theorem
of GIT is that $x$ is GIT (semi)stable if and only if it is
(semi)stable with respect to every 1-PS of $\SL(V)$.
 We will sometimes abuse terminology and
say that $x$ is (semi)stable with respect to $\rho$ when it is (semi)stable
with respect to $\rho^{\circ}$.

\subsection{GIT of Chow points}

We briefly recall Mumford's interpretation of Hilbert-Mumford
criterion for Chow stability of projective varieties \cite{M}.
Let $X \subset  \bP(V)$ 
be a projective variety and $\rho$ be a one-parameter subgroup
of $\GL(V)$ with weights $r_0 \ge r_1 \ge \cdots \ge r_N = 0$.
Let ${\mcl I}_{\rho}$ be the ideal sheaf of $\SS_X[t]$ such that
\[
{\mcl I}_{\rho}\cdot\SS_X(1)[t] = \left(\begin{array}{c}\mbox{the
$\SS_X$-submodule of
$\SS_X(1)[t]$ generated by} \\
\mbox{$t^{r_i}x_i$, \, $i = 0, \dots, N$.}
\end{array}\right)
\]

\begin{defn} The {\it Hilbert-Samuel multiplicity} $e_{\rho}(X)$ is
the normalized leading coefficient of $P(n): = \chi(\mcl L^n/\mcl
I_{\rho}^n\mcl L^n)$ where $\mcl L$ is the invertible $\SS_X[t]$-module
$\SS_X(1)[t]$.
\end{defn}

Then the Hilbert-Mumford criterion can be
translated in terms of $e_{\rho}(X)$ as follows:

\begin{thm}\label{T:Mumford}\cite[Theorem 2.9]{M} 
The  Chow point of $X$ is stable (resp. semistable) if and only
if
\[
e_{\rho}(X) < (resp. \leq) \frac {\dim(X)+1}{N+1} \deg(X) \sum r_i
\]
for any one-parameter subgroup $\rho : \bG_m \to \GL(V)$ with
weights $r_0 \geq r_1 \geq \dots \geq r_N = 0$.
\end{thm}

We shall make frequent use of the following lemma which describes
how the Hilbert-Samuel
multiplicity is affected by the singular points:

\begin{lemma}\label{L:balance}\cite[Lemma 1.4]{Sch} Let $X$ be a reduced curve in $\bP(V)$ and
$\nu: \til{X} \to X$ be its normalization.

\begin{enumerate}
\item $e_{\rho}(X) = \sum_{P \in \tilde{X}} e_{\rho}(\til{X})_P$, where
$e_{\rho}(\til{X})_P$ denotes the normalized leading coefficient of
$\dim_k \SS_{\tilde{X}\ex \A^1}/\mcl I_{P\ex \{0\}}^m$.

\item Suppose that $v(\nu^*x_i) + r_i \geq a$ for all $i$ where $v$
is the natural valuation of $\SS_{P, \tilde{X}}$. Then $
e_{\rho}(\til{X}) \geq a^2.$
\end{enumerate}
\end{lemma}

We shall also use the following
 lemma that addresses the case in which $X$ is degenerate and
$\rho$ acts on $X$ trivially.

\begin{lemma}\label{L:degree}\cite[Lemma 1.2]{Sch}
 Let $X$ be an $r$-dimensional variety in $\bP^N$. Let
$\rho$ be a 1-PS of $\GL_{N+1}(\k)$ such that $\rho(t)\cdot x_i =
t^{r_i}x_i$, $r_0 \geq \dots \geq r_N = 0$. Suppose that $x_j,
x_{j+1}, \dots, x_N$ vanish on $X$ and $r_0 = r_1 = \dots =
r_{j-1}$. Then
\[
e_{\rho}(X) = (r+1) r_0 \deg(X).
\]
\end{lemma}

\subsection{GIT of Hilbert points}
\label{subsect:GITHP}

Let $X \subset \bP^N = \bP(V)$ be a projective variety with
Hilbert polynomial $P(t)$.    Choose an integer $m$ sufficiently
large so that
\begin{itemize}
\item{$\cO_X(m)$ has no higher cohomology;}
\item{the natural map
\[
\Sym^mV^* \rightarrow \Gamma(\cO_X(m))
\]
is surjective.}
\end{itemize}

\begin{defn}\label{D:mth-hilb-pt}
 The {\it $m$th Hilbert point $[X]_m$ of $X$} is defined
\[
[X]_m := \left[ \Sym^mV^* \to \Gm(\cO_X(m)) \right] \in
\Gr(P(m),\Sym^mV) \inj \bP(\bigwedge^{P(m)} \Sym^mV).
\]
\end{defn}
Note that $X$ is determined by $[X]_m$ provided $X$ is
cut out by forms of degree $m$.

\begin{defn} $X$ is said to be {\it $m$-Hilbert stable (resp.
semistable)} if $[X]_m$ is GIT stable (semistable) with respect to
the natural $\SL(V)$ action on $\bP(\bigwedge^{P(m)} \Sym^mV)$.
\end{defn}

We refer the reader to \cite{HHL} for detailed discussion of an
algorithm (and a Macaulay 2 implementation) using Gr\"obner basis to
determine whether a variety is $m$-Hilbert (semi)stable with respect
to a given one-parameter subgroup.  We sketch the main results here.

For any given $v \in \bR^{N+1}$, $\prec_{v}$ denotes the
monomial order defined by declaring $x^a \prec_{v} x^b$ if
\begin{enumerate}
\item $\deg x^a < \deg x^b$;
\item $\deg x^a = \deg x^b$ and $v. a < v. b$;
\item $\deg x^a = \deg x^b$, $v. a = v. b$ and $x^a
\prec_{{\rm Lex}} x^b$ in the lexicographic order.
\end{enumerate}
In particular, given a one-parameter subgroup $\rho$ with the
weight vector $w = (r_0,\dots,r_N)$, the monomial order
$\prec_\rho$ means the graded lexicographic order associated to
the weight $w$.  Given a
monomial $x^a = x_0^{a_0}\cdots x_N^{a_N}$, the $\rho$-weight is
defined
$$wt_{\rho}(x^a) := w.a = r_0a_0+\ldots+r_Na_N.$$
For each polynomial $f$, let $in_{\prec_\rho}(f)$ denote the
largest term of $f$ with respect to $\prec_\rho$.  For an ideal
$I\subset \Sym V^*$, we let $in_{\prec_\rho}(I) := \langle
in_{\prec_\rho}(f) \, | \, f \in I \rangle$. Let $I\subset \Sym
V^*$ be a homogeneous ideal with graded pieces $I_m=I\cap \Sym^m
V^*$. The monomials $\{x^{a(1)}, \dots, x^{a(P(m))}\}$ of degree
$m$ not contained in $in_{\prec_\rho}(I)$ form a basis for $\Sym^m
V^*/I_m$.

We reformulate Gieseker's stability
criterion for Hilbert points \cite[pp. 8]{Gies} in these terms:
\begin{prop}
\label{prop:stabcrit2} The Hilbert-Mumford index of $[X]_m$ with
respect to a one-parameter subgroup $\rho : \bG_m \to \GL(V)$
with weights $r_0, r_1, \dots, r_N$ is given by
\begin{equation}\label{eqn:stabcrit}
\mu([X]_m, \rho) = \frac{m P(m)}{N+1} \sum r_i -\sum_{j=1}^{P(m)} wt_{\rho}(x^{a(j)})
\end{equation}
 where
$a(1),\ldots,a(P(m))$ index the monomials of degree $m$ not
contained in $in_{\prec_\rho}(I)$.
 In particular, $[X]_m\in \bP(\bigwedge^{P(m)} \Sym^mV)$ is
stable (resp. semistable) under the natural $\SL(V)$-action if and
only if for any one-parameter subgroup $\rho$ we have
$$\sum_{j=1}^{P(m)} wt_{\rho}(x^{a(j)})
<(\text{resp.} \leq) \frac{m P(m)}{N+1} \sum r_i.$$
\end{prop}

\subsection{Polarizations on Hilbert schemes}
\label{subsect:PHS} Let $Hilb$ be the connected component of the
Hilbert scheme containing $X$, $\cX \subset \bP(V) \times
Hilb$ the universal family, $\pi: \cX \to Hilb$ the natural projection, and $\cO_{\cX}(1)$ the polarization.

A coherent sheaf $\mathcal F$ on $\bP(V)$ is said to be {\em
$M$-regular} in the sense of Castelnuovo and Mumford \cite[ch.
14]{MumCAS} if $H^i({\mathcal F}(M-i))=0$ for each $i>0$. Suppose
that the ideal sheaf ${\mathcal I}_X$ is $M$-regular. It follows
that for each $m\ge M$
\begin{itemize}
\item{
$\Gamma({\mathcal I}_X(m)) \otimes V^*
\ra \Gamma({\mathcal I}_X(m+1))$
is surjective;}
\item{$H^i({\mathcal I}_X(m-i))=0$ for each $i>0$;}
\end{itemize}
and also
\begin{itemize}
\item{$\Sym^{m-1} V^* \ra \Gamma(\cO_X(m-1))$ is surjective;}
\item{$H^i(\cO_X(m-1-i))=0$ for $i>0$.}
\end{itemize}
In particular, $\cO_X$ is $(M-1)$-regular.
Conversely, if $\cO_X$ is $(M-1)$-regular and $M\ge 0$
then ${\mathcal I}_X$ is $M$-regular.  \cite[pp. 68]{EisGoS}

There exists an $M\gg 0$ such that {\em every } $[X] \in Hilb$ has
$M$-regular ideal sheaf \cite[ch.14]{MumCAS}. Then for each $m\ge
M$ we get a closed embedding \cite[ch. 15]{MumCAS}
$$
\begin{array}{rcl}
Hilb & \subset & \Gr(P(m),\Sym^m V)  \subset \bP(\bigwedge^{P(m)} \Sym^mV) \\
 \ [X] & \mapsto & [X]_m
\end{array}
$$
The universal quotient bundle $Q\ra \Gr(P(m),\Sym^m V)$ satisfies
$$Q|Hilb=\pi_* \cO_{\cX}(m),$$
and on taking determinants we find
$$\Lambda_m:=\cO_{\Gr}(1)|Hilb=\det(\pi_* \cO_{\cX}(m)).$$

\subsection{Tautological classes on the Hilbert scheme}
\label{subsect:TCHS} Recall the {\em tautological divisor classes}
developed in \cite{Fog} and \cite[Theorem 4]{Kn}: There exist
Cartier divisors $L_0,\ldots,L_{r+1}$ on $Hilb$ such that
\begin{equation} \label{eqn:taut}
\det(\bR^\bullet\pi_*\cO_{\cX}(m))= \sum_{i=0}^{r+1}
\binom{m}{i} L_i,
\end{equation}
where $r$ is the dimension of subschemes parametrized by $Hilb$. That is, the determinant of cohomology of $\cO_{\cX}(m)$ can be
expressed as a polynomial in the tautological class. This is a
relative version of the Hilbert polynomial of $\cX$ over $Hilb$.
It follows that the polarizations introduced above satisfy:
\begin{equation} \label{eqn:expansion}
\Lambda_m=\det(\pi_*\cO_{\cX}(m))=\sum_{i=0}^{r+1} \binom{m}{i} L_i.
\end{equation}

Using these formulas, we extend our definition:
\begin{defn}
For each $m\in \bZ$, write
$$\Lambda_m=\det(\bR^{\bullet}\pi_*\cO_{\cX}(m))=
\sum_{i=0}^{r+1} \binom{m}{i} L_i.$$
\end{defn}

In many situations the tautological divisors satisfy a
dependence relation:
\begin{prop} \label{prop:depend}
Let $Hilb$ denote a connected component of the Hilbert scheme
parametrizing subschemes in $\bP(V)$ of dimension $r$ and
$L_0,\ldots,L_{r+1}$ the tautological divisors on $Hilb$. Let
$Hilb^{\bullet,1} \subset Hilb$ denote an open subset
corresponding to subschemes $X$ where the following hold:
\begin{itemize}
\item{$\cO_X(1)$ has no higher cohomology;}
\item{the restriction map
$V^* \ra \Gamma(\cO_X(1))$
is an isomorphism.}
\end{itemize}
Over $Hilb^{\bullet,1}$ we have the relation $L_0+L_1=0$.
\end{prop}
In particular, if $r=1$ then Equation~\ref{eqn:expansion} takes the form
\begin{equation} \label{eqn:expansion2}
\Lambda_m=L_0+mL_1+\frac{m(m-1)}{2}L_2=(m-1)(L_1+\frac{m}{2}L_2).
\end{equation}
\begin{proof}
Let $\pi:\cX \ra Hilb^{\bullet,1}$ be the universal family embedded
in $\bP(V) \times Hilb^{\bullet,1}$.  Our first assumption implies
$\pi_*\cO_{\cX}(1)$ is locally free and
$$\Lambda_1=\det(\pi_* \cO_{\cX}(1))=L_0+L_1.$$
The second assumption implies we have a trivialization
(cf. \cite[pp.44]{Vbook})
$$\Gamma(\cO_{\bP(V)}(1))\otimes \cO_{Hilb^{\bullet,1}}\simeq \pi_*\cO_{\cX}(1).$$
In particular, it follows that $L_0+L_1=0$.
\end{proof}

\subsection{Hilbert points and Hilbert schemes}
\label{subsect:HPHS} We have seen that $Hilb$ admits an embedding
into $\Gr(P(m),\Sym^m V)$ for $m\gg 0$.  In practice, we are
usually interested in subsets of $Hilb$, that exclude degenerate
subschemes with very high Castelnuovo-Mumford regularity:
\begin{prop} \label{prop:hphs}
Let $Hilb^{\circ,m} \subset Hilb$ denote the open subset
parametrizing $[X] \in Hilb$ satisfying:
\begin{itemize}
\item{$\cO_X(m)$ has no higher cohomology;}
\item{$\Sym^m V^* \ra \Gamma(\cO_X(m))$ is surjective.}
\end{itemize}
Let $\pi:\cX\ra Hilb^{\circ,m}$ denote the universal family
restricted to this subset. Then we have
$$\det(\pi_*\cO_{\cX}(m))=\Lambda_m|Hilb^{\circ,m}=
\sum_{i=0}^{r+1} \binom{m}{i} L_i|Hilb^{\circ,m}$$ and there
exists a morphism
$$
\begin{array}{rcl}
\phi_m:Hilb^{\circ,m} & \ra &
    \Gr(P(m),\Sym^m V)  \subset \bP(\bigwedge^{P(m)} \Sym^mV) \\
\ [X]   & \mapsto & [X]_m
\end{array}
$$
such that $\phi_m^*\cO(1)=\Lambda_m$.
\end{prop}
\begin{rem}
$Hilb^{\circ,m}$ contains the open subset parametrizing subschemes
$X$ with $m$-regular ideal sheaf $\mathcal I_X$ and
$(m-1)$-regular structure sheaf $\cO_X$.  In particular,
$Hilb^{\circ,m}=Hilb$ for $m\gg 0$.
\end{rem}
\begin{proof}
The first assertion is just an application of
Equation~\ref{eqn:taut} in the special situation when
$\bR^i\pi_*\cO_{\cX}(m)=0$ for $i > 0$.  The morphism $\phi_m$ is just the
classifying map for the surjection of locally free sheaves
$$\Sym^m V^*\otimes \cO_{Hilb^{\circ,m}}
\twoheadrightarrow \pi_*\cO_{\cX}(m).$$
Again, if $Q \ra \Gr(P(m),\Sym^m V)$ is the universal quotient bundle then
$\phi_m^*Q=\pi_*\cO_{\cX}(m)$, and taking determinants gives the
equality of line bundles.
\end{proof}
Applying the functorial properties of the Hilbert-Mumford index
\cite[2.1]{GIT} we obtain:
\begin{coro} \label{coro:compareindex}
Retain the notation of Proposition~\ref{prop:hphs}. Suppose that
$[X] \in Hilb^{\circ,m}$ and $\rho:\bG_m \ra \GL(V)$ is a
one-parameter subgroup as before such that
$$\lim_{t\ra 0} \rho(t).[X] \in Hilb^{\circ,m}.$$
Then we have
\begin{equation}
\mu([X]_m,\rho)=\mu^{\Lambda_m}([X],\rho). \label{eqn:equalindices}
\end{equation}
\end{coro}

\subsection{Chow stability and Hilbert stability}
\label{subsect:CPHS}
We compare the geometric invariant
theory of the Hilbert points $[X]_m,m\gg 0$ with that of the Chow point
$Ch(X)$.

Let $\chow \subset \bP(\otimes^{r+1} \Sym^d V)$ denote the
corresponding {\em Chow variety}, i.e., the image of the Hilbert
scheme under the morphism \cite[\S 5.4]{GIT}
$$
\begin{array}{rcl}
\cycle: Hilb & \ra & \chow \subset \bP(\otimes^{r+1} \Sym^d V) \\
         \left[X\right] & \mapsto & Ch(X)
\end{array}.
$$
This is equivariant under the natural actions of $\SL(V)$.
By \cite[Theorem 4]{Kn}, we obtain the proportionality
\begin{equation} \label{eqn:Chowpolarization} \cycle^*\cO_{\chow}(1)\sim L_{r+1}
\end{equation}
and thus
\begin{equation} \label{eqn:hilbtochow}
\lim_{m\ra \infty} \frac{\Lambda_m}{\binom{m}{r+1}}\sim
\cycle^*\cO_{\chow}(1).
\end{equation}
In other words, the sequence $\{[\Lambda_m]\}$ converges to the
pull back of the Chow polarization in the
projectivized N\'eron-Severi group of the Hilbert scheme.

Let
$$\chow^s\subset \chow^{ss} \subset \chow$$
denote the locus of points
stable and semistable under the $\SL(V)$-action.
For each $m\gg 0$, let
$$Hilb^{s,m}\subset Hilb^{ss,m} \subset Hilb$$
the locus of points stable and semistable under the $\SL(V)$-action
linearized by $\Lambda_m$. The ample cone of $Hilb$ admits a
finite decomposition into locally-closed cells, such that the
semistable locus is constant for linearizations taken from a given
cell \cite[Theorem 0.2.3(i)]{DolHu}.  In particular, $Hilb^{s,m}$
and $Hilb^{ss,m}$ are constant for $m\gg 0$;  these are loci of the points
{\em stable and semistable with respect to the asymptotic
linearization}. While the linearization is not well-defined, its
locus of semistable points is!

Moreover, applying functoriality
of stability \cite[Theorem 2.1]{Rei}, we find
\begin{prop}  \label{prop:HilbtoChow}
Let $X\subset \bP^N$ be a variety.  If $X$ is Chow stable then $X$
is $m$-Hilbert stable for $m\gg 0$.  If $X$ is $m$-Hilbert
semistable for $m\gg 0$ then $X$ is Chow semistable.
\end{prop}
\begin{coro} \label{coro:HilbtoChow}
Assume the cycle map induces an $\SL(V)$-equivariant map
$$\cycle:Hilb^{ss,m} \ra  \chow^{ss},$$
which is the case for $m\gg 0$.  Then we obtain a natural
morphism of GIT quotients
$$\cycle: Hilb^{ss,m}/\!\!/\SL(V) \ra \chow^{ss}/\!\!/\SL(V).$$
\end{coro}

\subsection{Filtered Hilbert polynomials}
\begin{defn}
Given a graded ideal $I \subset \Sym V^*$, the {\em
filtered Hilbert function} is defined
$$H_{\Sym V^*/I,\rho}(m)=\sum wt_{\rho}(x^a)$$
where the sum is taken over the monomials of degree $m$ not
contained in $in_{\prec_\rho}I$. For a closed subscheme $X\subset
\bP(V)$, we define
$$H_{X,\rho}=H_{\Sym V^*/I_X,\rho}.$$
\end{defn}
\begin{prop} \label{prop:fhf}
The filtered Hilbert function $H_{X,\rho}(m)$
is a polynomial $P_{X,\rho}(m)$ for $m\ge M$,
the Castelnuovo-Mumford regularity of $\cO_X$.  This polynomial is called the
{\em filtered Hilbert polynomial}.
\end{prop}
\begin{proof}
Since $m\ge M$
we have an embedding
$$Hilb \subset \Gr(P(m),\Sym^mV)  \inj \bP(\bigwedge^{P(m)} \Sym^mV)$$
and Proposition~\ref{prop:hphs} implies
$$
\mu([X]_m,\rho)=\mu^{\Lambda_m}([X],\rho).$$
Equation~\ref{eqn:expansion} gives
$$\Lambda_m=\sum_{i=0}^{r+1} \binom{m}{i} L_i$$
Fixing the point and the one-parameter subgroup,
$\mu$ is a homomorphism
in the line bundle variable
\cite[2.2]{GIT} and we have
\begin{equation} \label{eqn:expand}
\mu^{\Lambda_m}([X],\rho)=\sum_{i=0}^r \binom{m}{i}
\mu^{L_i}([X],\rho).
\end{equation}
The result follows from Equation~\ref{eqn:stabcrit}.
\end{proof}

\subsection{Hilbert schemes of curves}
In this section, we assume $Hilb$ parametrizes schemes of pure
dimension one.  Here Equation~\ref{eqn:taut} takes the form
$$\Lambda_m=
L_0+mL_1+\binom{m}{2}L_2.
$$
\begin{prop} \label{prop:curveresult}
Let $Hilb^{\bullet} \subset Hilb$ denote the open subset
parametrizing $[X] \in Hilb$ satisfying:
\begin{itemize}
\item{$X$ is connected of pure dimension one;}
\item{$V^* \ra \Gamma(\cO_X(1))$ is an isomorphism;}
\item{$\cO_X$ is $2$-regular.}
\end{itemize}
Then for each $m \ge 2$ we have
\begin{equation}
\mu([X]_m,\rho)=\mu^{\Lambda_m}([X],\rho)=
(m-1)\left[(3-m)\mu^{\Lambda_2}([X],\rho)+(m/2-1)\mu^{\Lambda_3}([X],\rho)\right].
\label{eqn:curveindices}
\end{equation}
\end{prop}
\begin{proof}
Proposition~\ref{prop:depend} gives the relation $L_0+L_1=0$ over
$Hilb^{\bullet}$, and Equation~\ref{eqn:expansion2} gives
$$\Lambda_m|{Hilb^{\bullet}}=
(m-1)(L_1+\frac{m}{2}L_2).$$
Under our regularity hypothesis, Proposition~\ref{prop:hphs}
applies for each $m\ge 2$ and
\begin{equation}\label{E:HisMu}
\mu^{\Lambda_m}([X],\rho) = (m-1) \left[
\mu^{L_1}([X],\rho) + \frac m2 \mu^{L_2}([X],\rho)\right]
\end{equation}
is a polynomial  for $m \ge 2$ (see Proposition~\ref{prop:fhf}). One
can obtain  (\ref{eqn:curveindices}) by expressing
$\mu^{L_1}([X],\rho)$ and $\mu^{L_2}([X],\rho)$ in terms of
$\mu^{\Lambda_2}([X],\rho)$ and $\mu^{\Lambda_3}([X],\rho)$ and
plugging them in (\ref{E:HisMu}).
\end{proof}

\section{Basin  of attraction and equivalences}
\label{S:basin}
\begin{defn}\label{D:basin}
Let $X$ be a variety with $\bG_m$ acting via $\rho:\bG_m \ra Aut(X)$
with fixed points $X^{\rho}$.
For each $x^{\star} \in X^{\rho}$, the {\em basin of attraction}
is defined
\[
A_\rho(x^\star) := \left\{ x\in X \, | \, \lim_{t \to 0} \rho(t).x = x^\star
\right\}.
\]
\end{defn}

When $X$ is smooth and projective this can be interpreted via
the Bia{\l}ynicki-Birula decomposition \cite[Theorem 4.3]{Bia}:
Consider the decomposition $X^{\rho}_i,i\in I$ of the fixed points into connected components.
Then there is a unique locally closed $\rho$-invariant decomposition $X=\cup_{i\in I}X_i$ and
morphisms $\gamma_i:X_i \ra X^{\rho}_i$ such that
\begin{itemize}
\item{$(X_i)^{\rho}=X^{\rho}_i$ for each $i\in I$;}
\item{$\gamma_i$ is an affine bundle;}
\item{for each $x \in X_i^{\rho}$, the tangent space $T_xX_i \subset T_xX$ is the subspace
over which $\rho$ acts with nonnegative weights.}
\end{itemize}
For $x^{\star}\in X^{\rho}_i$ we have $A_{\rho}(x^{\star})=\gamma_i^{-1}(x^{\star})$.

The importance of this decomposition for the analysis of semistable points is clear
from the following proposition which is well known to experts. 
Given a point on a projective variety
$x\in X \subset \bP^N$,
let $x^*\in \bA^{N+1}$ denote an affine lift, i.e., a point in the affine cone over $X$ lying over $x$.
\begin{prop} \label{P:whybasin}
Suppose that $G$ is a reductive linear algebraic group acting on a projective variety $X$ and $L$
is a $G$-linearized ample line bundle.
Suppose $x_1,x_2\in X$ be semistable points mapping to the
same point in the GIT quotient $X/\!\!/G$.  Then there exists a semistable point $x_0 \in X$
with the following properties:
\begin{itemize}
\item{the orbit $Gx_0^*$ is closed, or equivalently, the stabilizer $G_{x_0^*} \subset G$ is reductive;}
\item{there exists $g\in G$, one-parameter subgroups $\rho_1,\rho_2$ of $G_{x_0^*}$,
and lifts $x_1^*$ and $x_2^*$ of $x_1$ and $x_2$
such that
$$x_1^* \in A_{\rho_1}(x^*_0) \quad g\cdot x^*_2 \in A_{\rho_2}(x^*_0).$$}
\end{itemize}
\end{prop}
\begin{proof}
Since $x_1$ and $x_2$ are identified in the GIT quotient, any homogeneous invariant
vanishing on $x_1$ automatically vanishes on $x_2$, and vice versa.
Consider the orbit closures $\bar{Gx_1}$ and $\bar{Gx_2}$ in $X$.
Their orbit closures meet \cite[Proposition 7, pp. 254]{Ses}:
$$ \bar{Gx_1} \cap \bar{Gx_2} \neq \emptyset,$$
and moreover there exist $x_1^*$ and $x_2^*$ lying over $x_1$ and $x_2$
in the affine cone over $X$ such that
$$ \bar{Gx^*_1} \cap \bar{Gx^*_2} \neq \emptyset.$$
(This is essentially the fact that invariants separate
orbit closures in affine space, e.g., \cite[Corollary 1.2, pp. 29]{GIT}.)
Pick $y_0^* \in \bar{Gx^*_1}\cap \bar{Gx^*_2}$ generating a {\em closed}
orbit of the intersection.

Recall Matsushima's Criterion \cite{Mat,Bia63}:  Suppose
$G$ is a reductive algebraic
group and $H\subset G$ a closed subgroup;  the homogeneous space
$G/H$ is affine if and only if $H$ is reductive.
This gives the equivalence of the two conditions on $x^*_0$.

We apply \cite[Theorem 1.4]{Kempf} to the
closed $G$-invariant set $S=\bar{Gy_0^*}=Gy_0^*$:   There
exist one-parameter subgroups $\rho_1$ and $\rho'_2$ such that
$$x_0^*:=\lim_{t\to 0}\rho_1(t)\cdot x^*_1 \in Gy_0^* \quad
\lim_{t\to 0}\rho'_2(t)\cdot x^*_2 \in Gy_0^*.$$
Clearly there exists $g\in G$ such that
$$g\cdot \lim_{t\to 0}\rho'_2(t)\cdot x^*_2=x_0^*.$$
Setting $\rho_2=g\rho'_2 g^{-1}$, we obtain the desired result.
\end{proof}

Also,  as far as stability is concerned, the points in a basin of
attraction are all equivalent if the attracting point is strictly semistable
with respect to the 1-PS:

\begin{lemma}\label{L:flat limit}
Let $G$, $X$, $L$ be as in Proposition~\ref{P:whybasin}.
Let $x \in X$ and suppose there exists an $x_0 \in X$ and a
one-parameter subgroup $\rho$ of $G$ such that $x\in A_{\rho}(x_0)$.
If $\mu^L(x, \rho) = 0$
then $x_0$ is semistable with respect
to $L$ if and only if $x$ is semistable with respect to $L$.
\end{lemma}

\begin{proof}
Assume that $X$ is embedded in $\bP^N$ by sections
of $L$ and $x^*\in \bA^{N+1} \setminus \{0\}$ be an affine lift of
$x$. Since $\mu^L(x, \rho) = 0$, $\rho(t). x^*$ has a
specialization, say $x_0^* \ne 0$, which corresponds to
$x_0\in \bP^N$.  Let $s$ be a $G$-invariant section of $L$. Then $s(x_0^*) =
s(\rho(t).x^*) = s(x^*)$ and it follows that $s(x_0) \ne 0$ if and
only if $s(x) \ne 0$.
\end{proof}

Let $x_0$ be a point
in $X$ and  $\rho : \bG_m \to G$ be a
one-parameter subgroup fixing $x_0$.
It follows directly from the definition of the Hilbert-Mumford
index that if $\mu^L(x_0, \rho) < 0$ then every point in the
basin of attraction $A_\rho(x_0)$ is unstable.  This observation
can be used to classify unstable points in certain situations;
our approach is similar to \cite[\S 4]{thaddeus}:

\begin{prop} \label{P:destabilize}
Let $G$, $X$, $L$ be as in Proposition~\ref{P:whybasin}
and $M$ a second $G$-linearized ample line bundle.
Let $x_1 \in X$ be semistable with respect to $L$ but unstable
with respect to $L\otimes M^{\epsilon}$
for each rational $\epsilon>0$.
Then there exists a point $x_0 \in X$
having the following properties:
\begin{enumerate}
\item{$x_0$ is strictly semistable with respect to $L$;}
\item{there exists a one-parameter subgroup $\rho:\bG_m \ra G_{x_0}$
such that
$$x_1 \in A_{\rho}(x_0);$$}
\item{
$\mu^{L \otimes M^{\epsilon}}(x_0,\rho)<0$.}
\end{enumerate}
\end{prop}
That is, every strictly semistable point that becomes unstable after
perturbing $L$ can be destabilized by a one-parameter
subgroup acting via automorphisms of a point strictly semistable
with respect to $L$. 
\begin{proof}
Let $\rho$ be a one-parameter subgroup with
$\mu^{L \otimes M^{\epsilon}}(x_1,\rho)<0$,
which exists by the Hilbert-Mumford criterion.  Let
$x_0=\lim_{t\ra 0}\rho(t)\cdot x_1$ denote the corresponding
limit point in $X$.  Clearly, $x_1 \in A_{\rho}(x_0)$ and
$\mu^{L \otimes M^{\epsilon}}(x_0,\rho)<0$.
\end{proof}

For the convenience of the reader, we recall the standard
Semistable Replacement Theorem:
\begin{thm} \label{T:ss-replacement}
Retain the assumptions of Lemma~\ref{L:flat limit}.
Assume that $G$ is reductive so the GIT quotient
scheme $X^{ss}/\!\!/G$ exists. Let $B$ be a smooth curve, $0\in B$
a closed point, and $f : B\setminus\{0\} \to X^{ss}$ be a regular
morphism. Then there exists a covering $\a : B' \to B$ branched only over
$0$ and $\gamma : B'\setminus\{0'\} \to G$, $0'= \a^{-1}(0)$,
such that
\begin{itemize}
\item there is a regular morphism $f' : B' \to X^{ss}$;
\item $f(\a(b')) = \gamma(b').f'(b')$ for all $b' \ne 0'$.
\end{itemize}
\end{thm}

\begin{defn} Two c-semistable curves $C_1$ and $C_2$ are
said to be {\it c-equivalent},
denoted $C_1 \sim_c C_2$,
if there exists a curve $C{^\star}$ (which we may assume has reductive
automorphism group) and one-parameter subgroups
$\rho_1,\rho_2$ of $Aut(C^{\star})$
with $\mu(Ch(C^\star),\rho_i) = 0$ such that
the basins of attraction $A_{\rho_1}(Ch(C^\star))$
and $A_{\rho_2}(Ch(C^\star))$
contain Chow-points of curves
isomorphic to $C_1$ and $C_2$ respectively.
\end{defn}
We define {\it h-equivalence}, denoted $\sim_h$,
in an analogous way.   Lemma~\ref{L:flat limit} shows that
these equivalence relations respect the semistable
and unstable loci.  Proposition~\ref{P:whybasin} shows 
that for GIT-semistable curves
$C_1 \sim_c C_2$ if and only if $Ch(C_1)$ and $Ch(C_2)$
yield the same point of $\Mcs$;  the analogous statement
holds for h-equivalence.

\section{Computations over the moduli space of stable curves}
\label{S:CMS}

Let $\pi:\FCg \ra \FMg$ denote the universal curve over the
moduli stack of curves of genus $g$.  For each $n\ge 1$ we have
the vector bundle
$$E_n=\pi_*\omega_{\pi}^n,$$
of rank
$$r(n)=\begin{cases} g  & \text{ if } n=1, \\
        (2n-1)(g-1) & \text{ if } n>1.
    \end{cases}$$
Write
$$\lambda_n=c_1(E_n)$$
and use $\lambda$ to designate $\lambda_1$.

Consider the multiplication maps
\begin{equation} \label{eqn:multiplication}
\Sym^m E_n \ra E_{mn}
\end{equation}
for each $m>1$.
We have the Chern-class identities
\begin{eqnarray*}
c_1(\Hom(\Sym^m(E_n),E_{mn})) &=& \text{rk}(\Sym^m(E_n)) c_1(E_{mn})-
                \text{rk}(E_{mn})c_1(\Sym^m(E_n)) \\
                  &=& \binom{m+r(n)-1}{m}\lambda_{mn} -
                r(mn)\binom{m+r(n)-1}{m-1}\lambda_n \\
                  &=&\frac{(m+r(n)-1)!}{m! \ r(n)!}
                \left(r(n)\lambda_{mn}-r(mn)m\lambda_n\right) \\
                &\sim& r(n)\lambda_{mn}-r(mn)m\lambda_n,
\end{eqnarray*}
where $\sim$ designates proportionality.  These
divisor classes were introduced by
Viehweg \cite[\S 1.4]{V} and Cornalba-Harris \cite[\S 2]{CH}
and their significance is explained by the following fact:

\begin{prop} \label{prop:descendHilb}  Assume that $n\ge 2$ when $g=2$.
Consider the Hilbert scheme $Hilb$ of degree $2n(g-1)$ curves of
genus $g$ in $\bP^{r(n)-1}$.  Let $\Hgn \subset Hilb$ denote the
closure of the $n$-pluricanonically embedded smooth curves of genus $g$.
Suppose that $\Lambda_m$ (introduced in \S\ref{subsect:GITHP}) is
well-defined and ample on $\Hgn$.

Consider the open subsets
$$V_{g,n}^{s,m} \subset \Hgn^{s,m} \subset \Hgn, \quad
V_{g,n}^{ss,m} \subset \Hgn^{ss,m} \subset \Hgn$$
corresponding to $n$-canonically embedded
Deligne-Mumford stable curves that are
GIT stable and GIT semistable with respect to $\Lambda_m$.  Let
$$\cU^{s,m}_{g,n} \subset \cU^{ss,m}_{g,n} \subset \FMg$$
denote their images in moduli.  Then $\Lambda_m$ descends to
a multiple of
$r(n)\lambda_{mn}-r(mn)m\lambda_n$ along $\cU^{ss,m}_{g,n}$.  This
restricts to an ample divisor on the coarse moduli space
$U^{s,m}_{g,n}$.
\end{prop}
\begin{proof} (cf. \cite[\S 1.6]{Vbook})
We illustrate how $\Lambda_m$ descends to $\cU^{ss,m}_{g,n}$.
Let $\varpi:\cX \ra \Hgn$ denote the universal family.
The multiplication map (\ref{eqn:multiplication}) on moduli is obtained
by descent from the multiplication map over $\Hgn$
$$\Sym^m(\varpi_*\cO_{\cX}(1)) \ra  \varpi_*\cO_{\cX}(m).$$
As in Proposition~\ref{prop:depend}, we have a trivialization
$$\Gamma( \cO_{\bP^{r(n)-1}}(1))\otimes \cO_{V_{g,n}^{ss,m}}
\simeq \varpi_*\cO_{\cX}(1)|V_{g,n}^{ss,m}.$$
Thus the divisor class
$$c_1\left(\Hom(\Sym^m(\varpi_*\cO_{\cX}(1)),\varpi_*\cO_{\cX}(m))|V_{g,n}^{ss,m}\right)$$
is proportional to
$$c_1(\varpi_*\cO_{\cX}(m)|V_{g,n}^{ss,m})=\Lambda_m|V_{g,n}^{ss,m}.$$

As for the ampleness, the coarse moduli space $U^{s,m}_{g,n}$ of $\cU^{s,m}_{g,n}$
can be identified with an open subset of the GIT quotient
$$\Hgn^{ss,m}/\!\!/\SL_{r(n)}.$$
$\Lambda_m$ descends to a polarization of this quotient.
\end{proof}

Mumford \cite[Theorem 5.10]{M} showed that
the Grothendieck-Riemann-Roch formula gives
\begin{equation} \label{eqn:grr}
\lambda_n=(6n^2-6n+1)\lambda-\binom{n}{2}\delta,\quad  n>1.
\end{equation}
We therefore find \Small
\begin{equation}\label{E:Lm}
r(n)\lambda_{mn}-r(mn)m\lambda_n =
    \begin{cases} \lambda+(m-1)((4g+2)m-g+1)\lambda -\frac{gm}{2}\delta)
    \, \text{ if } n=1, \\
            (m-1)(g-1)((6mn^2-2mn-2n+1)\lambda -
                \frac{mn^2}{2}\delta)
                \, \text{ if } n>1.
    \end{cases}
\end{equation}\normalsize Asymptotically as $m\ra \infty$, we obtain the
proportionality
\[ \lim_{m\ra \infty}
\left(r(n)\lambda_{mn}-r(mn)m\lambda_n\right) \sim
    \begin{cases} (4g+2)\lambda -\frac{g}{2}\delta
            \text{ if } n=1 \\
            (6n-2)\lambda -
                \frac{n}{2}\delta
            \text{ if } n>1.
    \end{cases}
\]\normalsize Combining Proposition~\ref{prop:descendHilb} and
Equation~\ref{eqn:hilbtochow}, we obtain:
\begin{prop} \label{prop:descendChow}  Assume that $n\ge 2$ when $g=2$.
Consider the Chow variety $\chow$ of degree $2n(g-1)$ curves of
genus $g$ in $\bP^{r(n)-1}$.  Let $\Cgn \subset \chow$ denote the
closure of the $n$-pluricanonically embedded curves of genus $g$.

Consider the open subsets
$$V^{s,\infty}_{g,n} \subset \Cgn^s  \subset \Cgn, \quad
V^{ss,\infty}_{g,n} \subset \Cgn^{ss} \subset \Cgn$$
corresponding to $n$-canonically embedded
Deligne-Mumford stable curves that are
Chow stable and Chow semistable respectively.  Let
$$\cU^{s,\infty}_{g,n} \subset \cU^{ss,\infty}_{g,n} \subset \FMg$$
denote their images in moduli.  Then the polarization descends to
a multiple of
$$
(4g+2)\lambda -\frac{g}{2}\delta
                        \text{ if } n=1
$$
or
$$
                        (6n-2)\lambda -
                                \frac{n}{2}\delta
                        \text{ if } n>1.
$$
The restriction to the coarse moduli space $U^{s,\infty}_{g,n}$ is ample.
\end{prop}

\begin{rem}[Application to polarizations on $\Mg$]
Mumford has proven that $\cU^{s,\infty}_{g,n}=\FMg$ for each $n\ge 5$
\cite[Theorem 5.1]{M}.  Proposition~\ref{prop:HilbtoChow}
then guarantees that $\cU^{s,m}_{g,n}=\FMg$ for all $m\gg 0$.
Proposition~\ref{prop:descendChow} then implies
that $a\lambda-\delta$ is ample for $a>11.2$ \cite[Corollary 5.18]{M}.
Cornalba and Harris \cite{CH} established the sharp result:
$a\lambda-\delta$ is ample if and only if $a>11$.
\end{rem}

We are primarily interested in situations where not all Deligne-Mumford
stable curves have stable Hilbert/Chow points.  Here
GIT yields alternate birational models of the moduli space.

\begin{thm} \label{thm:getcontraction}
Retain the notation of Propositions~\ref{prop:descendHilb} and
\ref{prop:descendChow} with the convention that
$m=\infty$ in the Chow case.  Suppose that
\begin{itemize}
\item{
the complement to the Deligne-Mumford stable curves in
the GIT-semistable locus $\Hgn^{ss,m}$
(resp. $\Cgn^{ss}$)
has codimension $\ge 2$;}
\item{there exist Deligne-Mumford stable curves in
the GIT-stable locus $\Hgn^{s,m}$
(resp $\Cgn^s$).}
\end{itemize}
Then there exists a birational contraction
$$F:\Mg \dashrightarrow  \Hgn^{ss,m}/\!\!/\SL_{r(n)} \ \text{(resp. }
\Cgn^{ss}/\!\!/\SL_{r(n)})
$$
regular along the Deligne-Mumford stable curves with GIT-semistable
Hilbert (resp. Chow) points.

If ${\mcl L}_m$ is the polarization on the GIT quotient induced by
$\Lambda_m$ then the moving divisor
\begin{equation} \label{eqn:PicCongruence}
F^*{\mcl L}_m\sim \,
r(n)\lambda_{mn}-r(mn)m\lambda_n \!\!\pmod{\mathrm{Exc}(F)},
\end{equation}
where $\mathrm{Exc}(F)\subset \mathrm{Pic}(\Mg)$ is the subgroup
generated by $F$-exceptional divisors.
\end{thm}
A rational map of proper normal varieties is said to be a
{\em birational contraction} if it is birational and its
inverse has no exceptional divisors.
Note that Propositions~\ref{prop:descendHilb} and \ref{prop:descendChow}
cover the case where $F$ is an isomorphism.
\begin{proof}
Our assumptions can be written in the notation of
Propositions~\ref{prop:descendHilb} and
\ref{prop:descendChow}:
\begin{itemize}
\item{
$V^{ss,m}_{g,n} \subset \Hgn^{ss,m}$
(resp. $V^{ss,\infty}_{g,n} \subset \Cgn^{ss}$)
has codimension $\ge 2$;}
\item{$V^{s,m}_{g,n}\neq \emptyset$.}
\end{itemize}

The GIT quotient morphism
$$V^{s,m}_{g,n} \ra \cU^{s,m}_{g,n}$$
identifies the stack-theoretic quotient
$[V_{g,n}^{s,m}/\SL_{r(n)}]$ with $\cU_{g,n}^{s,m}$.
This gives a birational map
$$\Hgn^{ss,m}/\!\!/\SL_{r(n)}\ \text{(resp. }
\Cgn^{ss}/\!\!/\SL_{r(n)}) \dashrightarrow \Mg;$$
we define $F$ as its inverse.

We establish that $F$ is regular along $U_{g,n}^{ss,m}$:
We have an $\SL_{r(n)}$-equivariant morphism
$$\Hgn^{ss,m} \ra \Hgn^{ss,m}/\!\!/\SL_{r(n)},$$
which descends to
$$\cU^{ss,m}_{g,n} \ra \Hgn^{ss,m}/\!\!/\SL_{r(n)}.$$
Recall the universal property of the coarse moduli space:
Any morphism from a stack to a scheme factors through its
coarse moduli space.  In our context, this gives
$$U^{ss,m}_{g,n} \ra \Hgn^{ss,m}/\!\!/\SL_{r(n)}
\text{(resp. } \Cgn^{ss}/\!\!/\SL_{r(n)}).$$
Furthermore, the total transform of $\Mg \setminus U_{g,n}^{ss,m}$
is contained in the complement  $\Hgn^{ss,m} \setminus V_{g,n}^{ss,m}$
(\text{resp.} $\Cgn^{ss} \setminus V_{g,n}^{ss,\infty}$),
which has codimension $\ge 2$.
Thus any divisorial components of
$\Mg \setminus U_{g,n}^{ss,m}$ are $F$-exceptional divisors.
Similary, $F^{-1}$ has no exceptional divisors:
These would give rise to divisors in the complement
to $V_{g,n}^{ss,m}$ in the semistable locus.

We now analyze $F^*\mcl L_m$ in the rational Picard group of $\Mg$.
(Since $\Mg$ has quotient singularities, its Weil divisors
are all $\bQ$-Cartier.)
If $\mcl L_m^a$ is very ample on the GIT quotient then
$F^*\mcl L_m^a$ induces $F$, i.e., $F^*\mcl L_m^a$ has no fixed components
and is generated by global sections over $U^{ss,m}_{g,n}$.
Now $F^*\mcl L_m$ is proportional to
$r(n)\lambda_{mn}-r(mn)m\lambda_n$
over $U^{ss,m}_{g,n}$ and Formula (\ref{eqn:PicCongruence}) follows.
\end{proof}

\section{Properties of c-semistable and h-semistable curves}
\label{S:propchstable}

\subsection{Embedding c-semistable curves}
\begin{prop} \label{P:embedcstable}
 If $g \ge 3$ and $C$ is a c-semistable curve of genus $g$ over $k$,
then $H^1(C, \o^{\ten n}_C) = 0 $ and $\o_C^{\ten n}$ is
very ample for $n \ge 2$.
\end{prop}
\begin{rem} For the rest of this paper, when we refer to the Chow
or Hilbert point of a c-semistable curve $C$ it is with respect to
its bicanonical embedding in $\bP(\Gamma(C,\omega_C^{\ten 2})^*)$.
\end{rem}

\begin{proof}  Our argument follows \cite[Theorem 1.2]{DM}.

By Serre Duality, $H^1(C, \o^{\ten n}_C)$ vanishes if
$H^0(C, \o^{\ten 1-n}_C)$ vanishes.  The restriction of
$\o^{\ten 1-n}_C$ to each irreducible component $D\subset C$ has
negative degree because $\omega_C$ is ample.  It follows that
$\Gamma(D,\o^{\ten 1-n}_C|D)=0$, hence
$\Gamma(C,\o^{\ten 1-n}_C)=0$.

To show that $\o_C^{\ten n}$ is very ample for $n\ge 2$, it
suffices to prove for all $x, y \in C$ that
\begin{equation}\label{E:va}
\Hom(\mfk m_x\mfk m_y, \o_C^{\ten -n}) =  0, \quad n \ge 1.
\end{equation}
Let $\pi:C'\to C$ denote the partial normalization of any
singularities at $x$ and $y$.
When $x$ is singular, a local computation gives
$$
\Hom(\mfk m_x, \mcl L) \simeq \Gamma(C', \pi^*\mcl L).
$$
If $x$ is a cusp and $x' \in C'$ its preimage then
$$
\Hom(\mfk m_x^2, \mcl L) \simeq \Gamma(C', \pi^*\mcl L(2 x')).
$$
If $x$ is a node or tacnode and $x_1,x_2 \in C'$ the
preimage points then
$$
\Hom(\mfk m_x^2, \mcl L) \simeq \Gamma(C', \pi^*\mcl L(x_1+x_2)).
$$
Thus in each case we can express
$$\Hom(\mfk m_x \mfk m_y, \o_C^{\ten -n})=\Gamma(C', \mcl M)$$
for a suitable invertible sheaf $\mcl M$ on $C'$.  Moreover,
we have an inclusion
$$\pi^*\omega_C^{-n} \hookrightarrow \mcl M$$
with cokernel $Q$ supported in $\pi^{-1}\{x,y\}$ of length $\ell(Q)\le 2$.
For instance, if both $x$ and $y$ are smooth then
$$\mcl M=\omega_C^{-n}(x+y);$$
if both $x$ and $y$ are singular and $x\neq y$ then
$$\mcl M=\pi^*\omega_C^{-n}.$$

Suppose that for each irreducible component $D'\subset C'$,
the degree $\deg \mcl M|D'<0$.  Then $\Gamma(C', \mcl M)=0$ and
the desired vanishing follows.  We therefore classify
situations where
$$\deg \mcl M|D'=-n \deg \pi^*\omega_C|D'+\ell(Q|D') \ge 0,$$
which divide into the following cases:
\begin{enumerate}
\item[(a)]{$\deg \pi^*\omega_C|D'=1$, $n=1$, $\ell(Q|D')=1$;}
\item[(b)]{$\deg \pi^*\omega_C|D'=1$, $n=1,2$, $\ell(Q|D')=2$;}
\item[(c)]{$\deg \pi^*\omega_C|D'=2$, $n=1$, $\ell(Q|D')=2$.}
\end{enumerate}
We write $D=\pi(D')\subset C$.

We enumerate the various possibilities.  We use the assumption
that $C$ is c-semistable and thus has no elliptic tails. In cases (a)
and (b), $D$ is necessarily isomorphic to $\bP^1$ and meets the
rest of $C$ in either three nodes or in one node and one tacnode.
After reordering $x$ and $y$, we have the following subcases:
\begin{enumerate}
\item[(a1)]{$x=y \in D$ a node or tacnode of $C$;}
\item[(a2)]{$x \in D$ a node or tacnode of $X$ and $y\in D$
        a smooth point of $C$;}
\item[(a3)]{$x\in D$ a smooth point of $C$ and $y \not \in D$.}
\item[(b1)]{$x,y \in D$ smooth points of $C$.}
\end{enumerate}
In case (c), $D$ may have arithmetic genus zero or one:
\begin{enumerate}
\item[(c1)]{$D\simeq \bP^1$ with $x,y \in D$ smooth points of $C$;}
\item[(c2)]{$D$ of arithmetic genus one with $x,y \in D$ smooth
points of $C$;}
\item[(c3)]{$D$ of arithmetic genus one, $x=y$ a node or cusp of $D$,
and $D'\simeq \bP^1$.}
\end{enumerate}
In subcase (c1), $D$ meets the rest
of $C$ in either four nodes, or in two nodes and one tacnode, or in
two tacnodes.  In subcases (c2) and (c3), $D$ meets the rest of $C$
in two nodes.  Except in case (c3),
$\pi:D' \ra D$ is an isomorphism.

For subcases (b1), (c1), and (c2), $\pi$ is an isomorphism.
Moreover, $Q$ is supported along $D$
so $\mcl M$ has negative degree along any other irreducible components
of $C$.  There are other components because the genus of $C$ is
at least three.  Thus elements of $\Gamma(C,\mcl M)$ restrict to
elements of $\Gamma(D,\mcl M|D)$ that vanish at the points where
$D$ meets the other components, i.e., in at least
two points.  Since $\deg \mcl M|D=0$ or $1$, we conclude
$\Gamma(C,\mcl M)=0$.

For subcase (c3), $\pi$ is not an isomorphism but $Q$ is still
supported along $D'$.  As before, $\mcl M$ has negative degree
along other irreducible components of $C'$, and elements of
$\Gamma(C',\mcl M)$ restrict to elements of $\Gamma(D',\mcl M|D')$
vanishing where $D'$ meets the other components.
There are at least two such points but $\deg \mcl M|D'=0,1$,
so we conclude that $\Gamma(C',\mcl M)=0$.

In case (a), we have $\deg \mcl M|D =0$.
Subcases (a1) and (a2) are similar
to (b1) and (c1):
$Q$ is supported along $D'$ so elements in $\Gamma(C',\mcl M)$
restrict to elements of $\Gamma(D',\mcl M|D')$ vanishing at the points
where $D'$ meets the other components.  There is at least
one such point, e.g., the singularity not lying over $x$,
hence $\Gamma(C',\mcl M)=0$.

Subcase (a3) is more delicate.  If $D'$ is the unique component
such that $\deg(\mcl M|D')\ge 0$ then the arguments of the previous
cases still apply.  However, the support of $Q$ might not be confined
to a single component.  We suppose there are two
components $D'_1$ and $D'_2$
as described in (a3),
such that $\deg(\mcl M|D'_i)\ge 0$.
Since the genus of $C$ is $>2$, $C$ cannot just be the union
of $D'_1$ and $D'_2$;  there is at least one additional component meeting
each $D'_i$ at some point $z_i$, and the restriction of $\mcl M$
to this component has negative degree.  Thus elements of $\Gamma(C',\mcl M)$
restrict to elements of $\Gamma(D'_i,\mcl M|D'_i)$ vanishing at $z_i$,
which are necessarily zero.
\end{proof}

\begin{coro}\label{C:versal}
Let $C\subset \bP^{3g-4}$ be a c-semistable bicanonical curve.
\begin{itemize}
\item{$\cO_C$ is 2-regular.}
\item{The Hilbert scheme is smooth at $[C]$.}
\item{Let $p_1,\ldots,p_n$ denote the singularities of $C$ and
$\mathrm{Def}(C,p_i),i=1,\ldots,n$ their versal deformation spaces.  Then there
exists a neighborhood $U$ of $[C]$ in the Hilbert scheme such that
$$U \ra \prod_i \mathrm{Def}(C,p_i)$$
is smooth.}
\end{itemize}
\end{coro}
\begin{proof}
Proposition~\ref{P:embedcstable} yields
$$H^1(C,\cO_C(1))=H^1(C,\omega_C^{\ten 2})=0$$
which gives the regularity assertion.  This vanishing also implies
\cite[I.6.10.1]{Kol}
$$H^1(C,\mathrm{Hom}(I_C/I_C^2,\cO_C))=0;$$
since the singularities of $C$ are local complete intersections we have
$$\mathrm{Ext}^1(I_C/I_C^2,\cO_C)=H^1(C,\mathrm{Hom}(I_C/I_C^2,\cO_C))=0$$
thus the Hilbert scheme is unobstructed at
$[C]$ (see \cite[I.2.14.2]{Kol}).   The assertion about the
map onto the versal deformation spaces is \cite[I.6.10.4]{Kol}.
\end{proof}

\begin{coro}\label{C:all-m} Let $C \subset \bP^{3g-4}$ be a bicanonical
c-semistable curve and $C^\star$ denote the curve to which
$\rho(t).C$ specializes.  If $C^\star$ is a bicanonical
c-semistable curve then
$$
\mu([C]_m,\rho)=
(m-1)\left[(3-m)\mu([C]_2,\rho)+(m/2-1)\mu([C]_3,\rho)\right].
$$
Thus $[C]_m$ is stable (resp. strictly semistable, resp. unstable)
with respect to $\rho$ for each $m\ge 2$ if and only if
$\mu([C]_3,\rho) \ge 2 \mu([C]_2,\rho) > 0 $ (resp.
$\mu([C]_3,\rho) = \mu([C]_2,\rho) = 0$, resp.
$\mu([C]_3,\rho) \le
2 \mu([C]_2,\rho) < 0 $.)
The Chow point $Ch(C)$ is stable (resp. strictly semistable, resp.
unstable) with respect to $\rho$ if and only if
$\mu([C]_3,\rho)-2\mu([C]_2,\rho) >0$ (resp. $=0$, resp. $<0$.)
\end{coro}

\begin{proof} By Proposition~\ref{P:embedcstable}, a bicanonical
c-semistable is 2-regular and the assertion on the Hilbert
points  follows immediately from
Equation (\ref{eqn:curveindices}).
Equations (\ref{eqn:Chowpolarization}) and (\ref{eqn:expand}) allow us
to interpret the Hilbert-Mumford index of the Chow point
in terms of the leading coefficient of $\mu([C]_m,\rho)$ as a polynomial
in $m$.
\end{proof}

\subsection{Basic properties of tacnodal curves}
Let $C$ be a curve with a tacnode $r$, i.e., a singularity
with two smooth branches intersecting with simple tangency.
Let $\nu:D\ra C$ be the partial normalization of $C$ at $r$
and $\nu^{-1}(r)=\{p,q\} \subset D$ the conductor.  The
descent data from $(D,p,q)$ to $(C,r)$ consists of a choice
of isomorphism
$$\iota:T_pD \stackrel{\sim}{\ra} T_qD$$
identifying the tangent spaces to the branches.
Functions on $C$ pull back to functions $f$ on $D$ satisfying
$f(p)=f(q)$ and $\iota(f'(p))=f'(q)$.

Varying the descent data gives a one-parameter family
of tacnodal curves:
\begin{prop}\label{P:const-tac} Let $D$ be a reduced curve and
$p,q \in D$ distinct smooth points with local parameters
$\sigma_p$ and $\sigma_q$.  Each invertible linear transformation
$T_pD\ra T_qD$ can be expressed
$$\iota(t):\frac{\partial}{\partial \sigma_p} \mapsto t
\frac{\partial}{\partial \sigma_q}$$
for some $t\neq 0$;
let
$\bG_m   \simeq   \mathrm{Isom}(T_pD,T_qD)$
denote the corresponding identification.
Then there exists a family $\cC \ra \bG_m$,
a section $r: \bG_m \ra \cC$, and
a morphism
\[
\xymatrix{
D\times \bG_m  \ar[rr]^-{\nu}\ar[dr] & & \cC\ar[dl] \, \\
                      &  \bG_m & &
}
\]
such that
\begin{enumerate}
\item{$\nu$ restricts
to an isomorphism
$$D \setminus \{p,q\} \times  \bG_m
 \stackrel{\sim}{\ra} \cC \setminus r;$$}
\item{for each $t\in \bG_m$,
$r_t \in C_t$ is a tacnode and $\nu_t$
its partial normalization;}
\item{the descent data from $(D,p,q)$ to $(C_t,r_t)$ is given by $\iota(t)$.}
\end{enumerate}
Every tacnodal curve normalized by $(D,p,q)$
occurs as a fiber of $\cC \ra \bG_m.$
\end{prop}
If $D$ is projective of genus $g-2$ then each $C_t$ has genus $g$.

We sketch the construction of $\cC$:
$\iota(t)$ tautologically yields an identification over $\bG_m$
\begin{equation} \label{eqn:iotamap}
\iota:T_{p \times \bG_m} D \times \bG_m/\bG_m  \stackrel{\sim}{\longrightarrow}
T_{q \times \bG_m} D \times \bG_m/\bG_m
\end{equation}
which is the descent data from $D\times \bG_m$ to $\cC$.
Fiber-by-fiber, we get the universal family of tacnodal
curves normalized by $(D,p,q)$.

We will extend $\cC \ra \bG_m$ to a family of tacnodal curves
$\cC' \ra \bP^1$.  First, observe that the graph construction
gives an open embedding
$$\bG_m\simeq \mathrm{Isom}(T_pD,T_qD)
    \subset \bP(T_pD\oplus T_qD)\simeq \bP^1,$$
where $t=0$ corresponds to $[1,0]$ and $t=\infty$ corresponds to
$[0,1]$.  However, the identification (\ref{eqn:iotamap}) fails to
extend over all of $\bP^1$;  indeed, it is not even defined at
$p\times [0,1]$ and its inverse is not defined at $q\times [1,0]$.
We therefore blow up
$$\cD'=\mathrm{Bl}_{p\times [0,1],q\times [1,0]}D \times \bP^1$$
and consider the sections
$$\mfk p,\mfk q: \bP^1 \ra \cD'$$
extending $p \times \bG_m$ and $q\times \bG_m$.  Now
(\ref{eqn:iotamap}) extends to an identification
$$
\iota':T_{\mfk p} \cD'/\bP^1  \stackrel{\sim}{\longrightarrow}
T_{\mfk q} \cD'/\bP^1.
$$

\begin{prop}\label{P:limit-tac}  Retain the notation of
Proposition~\ref{P:const-tac}.
There exists an extension
$$
\begin{array}{rcccl}
                & \cC & \subset & \cC' &  \\
                & \downarrow & & \downarrow & \\
\mathrm{Isom}(T_pD,T_qD) \simeq &\bG_m & \subset  & \bP^1 & \simeq
                            \bP(T_pD \oplus T_qD)
\end{array}
$$
 where $\cC'\ra \bP^1$ denotes the family of
curves obtained from $\cD'$ and $\iota'$ by descent, $r':\bP^1 \ra
\cC'$ the tacnodal section, and $\nu':\cD' \ra \cC'$ the resulting
morphism. The new fiber $(C'_{0},r'(0))$ (resp.
$(C'_{\infty},r'(\infty))$) is normalized by $(D'_0=D\cup_q
\bP^1,p,\mfk q(0))$ (resp. $(D'_{\infty}=D\cup_p
\bP^1,\mfk{p}(\infty),q)$.)
\end{prop}
We say that the tacnodes in the family $\{C'_t, r_t\}_{t \in
\bP^1}$ are {\it compatible}, and that two curves are {\it compatible}
if one can be obtained from the other by replacing some tacnodes
by compatible tacnodes.

\section{Unstable bicanonical curves}
\label{S:UBC}
In this section, we show
that if a curve is not c-semistable then it has unstable Chow
point:
\begin{prop}\label{P:chow-imply-c}
If $Ch(C) \in \Cg2$ is Chow semistable then
$C\subset \bP^{3g-4}$ is c-semistable.
\end{prop}
We prove this by finding one-parameter subgroups destabilizing
curves that are not c-semistable.
Many statements in this section are fairly direct
generalizations of results in \cite{M} and \cite{Sch}.

\subsection{Badly singular curves are Chow unstable}
A Chow semistable bicanonical curve $C$ cannot have a triple
point, since $\frac d{N+1} = \frac{4g-4}{3g-3} < \frac 32$ and
this implies that $C$ is Chow unstable by Proposition 3.1 of
\cite{M}. We need to show that among the double points, only
nodes, ordinary cusps and tacnodes are allowed.

\begin{lemma} If $C$ has a non-ordinary cusp, then it is Chow
unstable.
\end{lemma}
\begin{proof} Suppose that $C$ has a non-ordinary cusp at $p$. Let
$\nu: \, \til{C} \to C$ be the normalization, $p' = \nu\inv(p)$ and
assume $p = [1, 0, \dots, 0]$.
 Recall that the singularity at $p$ is determined by the vanishing
 sequence $\left(a_i(\nu^*|\o_C^{\ten 2}|, p')\right)_{i=1}^{N+1}$
 which is the strictly
 increasing sequence determined by the condition
\[
\{ a_i(\nu^*|\o_C^{\ten 2}|, P) \, | \, i = 1, 2, \dots, N+1 \} = \{
\ord_{p'}(\sm) \, \, | \,\, \sm \ne 0 \in \nu^*|\o_C^{\ten 2}| \}.
\]
$C$ has a cusp at $p$ if and only if the vanishing sequence $(
a_i(\nu^*|\o_C^{\ten 2}|, p') )$ is of the form
$
(0, 2, \geq 3)$,
and it has an ordinary cusp if it is of the form $(0, 2, 3, \geq
4)$.

Hence if $C$ has a non-ordinary cusp at $p$, then we can choose
coordinates $x_0, \dots, x_N$ such that $\ord_{p'} x_0 = 0$,
$\ord_{p'}x_1 = 2$, $\ord_{p'}x_2 = 4$, and $ \ord_{p'} x_i \geq 5,
i = 3, 4, \dots, N.$ Let $\rho : \bG_m \to \GL_{N+1}(\k)$ be the
one-parameter subgroup such that $\rho(t).x_i = t^{r_i} x_i$, where
the weights are:
\[
(r_0, r_1, \dots, r_N) = (5, 3, 1, 0, \dots, 0).
\]
Then $\ord_{p'}x_i + r_i \geq 5$ for all $i$, and it follows from
Lemma~\ref{L:balance} that
\[
e_{\rho}(C) = e_{\rho}(\til{C}) \geq e_{\rho}(\til{C})_{p'} \geq 5^2 = 25,
\]
while $\frac {2d}{N+1}\sum r_i = \frac {2\cdot
4(g-1)}{3(g-1)}\cdot 9 = 24.$ The assertion now follows from
Theorem~\ref{T:Mumford}.
\end{proof}

\begin{lemma} Suppose $C$ has a singularity at $p$ such that
\[
\widehat{\cO}_{C,p} \simeq k[x,y]/(y^2 - x^{2s}), \, s \geq 3.
\]
Then $C$ is unstable.
\end{lemma}

\begin{proof} Let $\nu : \til{C} \to C$ be the normalization,
$\nu\inv(p) = \{ p_1, p_2 \}.$ Since the two branches of $C$ agree
to order $s$ at $p$, we may choose coordinates $x_0, \dots, x_N$
such that
\[
(\ord_{p_i} x_0, \dots, \ord_{p_i}x_N) = (0, 1, 2, \geq 3), \quad i
= 1, 2.
\]

Let $\rho$ be the one-parameter subgroup of $\GL_{N+1}(\k)$ with
weights $(r_0, \dots, r_N) = (3, 2, 1, 0, \dots, 0)$. Then we have
\[
\ord_{p_i} x_j + r_j \geq 3, \quad i = 1, 2 \, \,  \mbox{and} \, \,
j = 0, 1, \dots, N,
\]
and by Lemma~\ref{L:balance},
\[
e_{\rho}(C) = e_{\rho}(\til{C}) \geq e_{\rho}(\til{C})_{p_1} +
e_{\rho}(\til{C})_{p_2} \geq 2 \cdot 3^2 = 18,
\]
which is strictly greater than $\frac {2d}{N+1}\sum r_i = \frac
{2\cdot 4(g-1)}{3(g-1)}\cdot 6 = 16.$

\end{proof}

\begin{lemma} If $C$ has a multiple component, $C$ is
Chow unstable.
\end{lemma}

\begin{proof} Let $C_1$ be a component of $C$ with multiplicity
$n \geq 2$. Choose a smooth non-flex point $p \in C_1^{red}$ such that $p$
does not lie in any other component. Since $p$ is smooth on
$C_1^{red}$, we may choose coordinates $x_0, \dots, x_N$ such that
\[
(\ord_{p} x_0, \dots, \ord_{p}x_N) = (0, 1, 2, \geq 3).
\]
Let $\rho$ be the one-parameter subgroup of $\GL_{N+1}(\k)$ with
weights $(r_0, \dots, r_N) = (3, 2, 1, \dots, 0)$. Then we have
\[
\ord_{p} x_i + r_i \geq 3.
\]
This yields the inequality
\[
e_{\rho}(C) \geq n\cdot e_{\rho}(C_1) \geq 2 \cdot 3^2 = 18
\]
whereas $\displaystyle{\frac {2d}{N+1}\sum r_i = \frac 83 \cdot 6 =
16.}$
\end{proof}

\subsection{Polarizations on semistable limits of bicanonical curves}
We first prove that the semistable limit of a one parameter family
of smooth bicanonical curves is bicanonical:
\begin{prop}\label{P:pol-ss-limit} Let $\mcl C \to
\spec k[[t]]$ be a family of Chow semistable curves of genus $g$
such that the generic fibre $\mcl C_\eta$
 is smooth. If $\Phi : \mcl C \to \bP^{3g-4}_{k[[t]]}$
is an embedding such that $\Phi_\eta^*(\cO(1)) = \o_{\mcl
C_\eta/k[[t]]}^{\ten 2}$ then $\cO_{\mcl C}(1) = \o_{\mcl
C/k[[t]]}^{\ten 2}$.
\end{prop}
By \cite[4.15]{M}, nonsingular bicanonical curves are Chow stable.
Hence any Chow semistable curve is a limit of nonsingular bicanonical
curves and Proposition~\ref{P:pol-ss-limit}  implies that
if $C$ is not bicanonical, then $Ch(C) \not\in
Chow_{g,2}^{ss}$. In particular, a Chow semistable curve does not
have a smooth rational component meeting the rest of the curve in
$< 3$ points.
 Mumford proved the
statement for the $n$-canonical curves for $n \ge 5$, and his
argument can be easily modified to suit our purpose. It is an easy
consequence of (ii) of the following proposition, which, in
Mumford's words, says that the degrees of the components of $C$
are roughly in proportion to their {\it natural} degrees.

\begin{prop}[Proposition~5.5, \cite{M}] Let $C \subset \mbb P^{3g-4}$
be a connected curve of genus $g$ and degree $4g-4$. Then
\begin{itemize}
\item[(i)] $C$ is embedded by a non-special complete linear system.
\item[(ii)] Let $C = C_1 \cup C_2$
be a decomposition of $C$ into two sets of components such that $W =
C_1\cap C_2$ and $w = \# W$ (counted with multiplicity). Then
\[
\left| \, \deg C_1 - 2 \, \deg_{C_1}\o_C \, \right| \le \frac w2.
\]
\end{itemize}
\end{prop}
Mumford's argument goes through in the bicanonical case except for
the proof of $H^1(C_1, \cO_{C_1}(1)) = 0$. If $H^1(C_1,
\cO_{C_1}(1)) \ne 0$, then by Clifford's theorem we have
\[
h^0(C_1, \cO_{C_1}(1)) \le \frac{\deg(C_1)}2 + 1
\]
and the Chow semistability of $C$ forces
\[
w + 2 \, \deg C_1 \le \frac{2 \,\deg C}{3g-3} \, h^0(C_1,
\cO_{C_1}(1)).
\]
Combining the two, we obtain
\[
\deg C_1 \le 4 - \frac32 w.
\]
If $w \ne 0$, then $\deg C_1 \le 2$, hence $C_1$ is rational and
$H^1(C_1, \cO_{C_1}(1)) = 0$. If $w = 0$, then $\deg C_1 \le 4$
which is absurd since $C_1 = C$ and $\deg C_1 = 4g-4$.

We need to justify our use of Clifford's theorem here, as Chow
semistable bicanonical curves have cusps and tacnodes. We shall
sketch the proof of Gieseker and Morrison \cite{Gies} and highlight
the places where modifications are required to accommodate the
worse singularities (\cite{Gies} assumes that $C$ has only nodes).

\begin{thm}[Clifford's Theorem] Let $C \subset \bP^N$
be a reduced curve with nodes, cusps and tacnodes. Let $L$ be a line
bundle generated by sections. If $H^1(C, L) \ne 0$, then there is a
subcurve $C_1 \subset C$ such that
\[
h^0(C, L) \le \frac{\deg_{C_1}L}2 + 1.
\]
\end{thm}
\begin{proof}[Sketch proof] Suppose that $H^1(C,
L)\ne 0$ and $\varphi\ne 0 \in Hom(L, \o_C)$. Let $C_1$ be the
union of components where $\varphi$ does not vanish entirely and
$p_1, \dots, p_w$ be the intersection points of $C_1$ and
$\overline{C - C_1}$. Assume that $p_i$'s are ordered so that
$p_1, \dots, p_\ell$ are 
tacnodes. Then we have $$\o_C|_{C_1}(-2 \sum_{i=1}^\ell p_i - \sum_{i =
\ell+1}^w p_i) = \o_{C_1}.$$ We claim that $\varphi$ restricts to
give a homomorphism from $L_{C_1}$ to $\o_{C_1}$: Let $p_i$ be a
tacnode and let $D \not\subset C_1$ be the irreducible component
containing $p_i$. Since $\varphi$ vanishes entirely on $D$,
$\varphi$ must vanish to order $\ge 2$ at $p_i$ on $C_1$.
Likewise, $\varphi$ must vanish at each node. It follows that
 $\varphi|_{C_1}$ factors through $\o_C|_{C_1}(-2 \sum_{i=1}^\ell
p_i - \sum_{i = \ell+1}^w p_i)$. Let $s_1, \dots, s_r$ be a basis of
$Hom(L_{C_1}, \o_{C_1})$ such that $s_1 = \varphi$, and let $t_1,
\dots, t_p$ be a basis for $H^0(C, L)$ such that $t_1$ does not
vanish at the support of $s_1$ and at any singular points. It is
shown in \cite{Gies} that
\[
\begin{array}{cccccccccc}
[s_1, t_1], & [s_1, t_2], &  [s_1, t_3], & \dots, & [s_1, t_p], \\
            & [s_2, t_1], &  [s_3, t_1], & \dots, & [s_r, t_1]\\
\end{array}
\]
are linearly independent sections of $H^0(C_1, \o_{C_1})$, which
implies that $p + r - 1 \le p_a(C_1) + 1$. Combining it with the
Riemann-Roch gives the desired inequality.
\end{proof}

\subsection{Elliptic subcurves meeting the rest of the curve in one point}

Let $C$ be a Deligne-Mumford stable curve with an elliptic
tail $E\subset C$.  Then $\omega_C^{\ten 2}$ fails to be very ample along $E$
and thus $C$ does not admit a bicanonical embedding.
In particular, $C$ does not arise in GIT quotients of the
Chow variety/Hilbert scheme of bicanonical curves.

Here, we focus on curves with an elliptic subcurve meeting the
rest of the curve in a {\em tacnode}.

\begin{prop}\label{P:g1-1tac}
Let $C = E \cup_p R \cup_q D$ be a bicanonical curve consisting of a
rational curve $E$ with one cusp, a rational curve $R$ and  a genus
$g-2$ curve $D$ such that $p$ is a tacnode and $q$ is a node. Then
$C$ is Chow unstable with respect to a one-parameter subgroup coming
from its automorphism group.
\end{prop}

\begin{proof}
Restricting $\o_C^{\ten
2}$ we get
\[
\o_C^{\ten 2}|E \simeq \cO_E(4p), \quad \o_C^{\ten 2}|R \simeq
\cO_R(2), \quad \o_C^{\ten 2}|D \simeq \o_D^{\ten 2}(2q).
\]
Since $h^0(\o_C^{\ten 2}|D) = 3(g-2)-3+2 = 3g-7$, we can choose
coordinates so that
\[
E\cup_p R \subset \{x_6 = x_7 = \cdots = x_{3g-4} = 0\}
\]
and $ D \subset \{x_0 = x_1 = x_2 =  x_3 = 0\}$. $E$ and $R$ can be
parametrized by
\[
[s, t] \mapsto [s^4, s^2t^2, st^3, t^4, 0, \dots, 0]
\] and
\[
[u,v] \mapsto [0, 0, uv, u^2, v^2, 0, \dots, 0].
\]
The cusp is at $[1,0, \dots, 0]$, $p = [0,0,0,1,0,\dots,0]$, and
$q = [0,0,0,0,1,0,\dots, 0]$. 
Let $\rho$ be the 1-PS with weight $(0, 2, 3, 4, 2, \dots, 2)$. We
have
\[
e_\rho(C) \ge e_\rho(E)_p + e_\rho(R)_p + e_\rho(R)_q + e_\rho(D).
\]
On $E$ (and $R$), we have $v_p(x_i) + r_i \ge 4$ for all $i$ where
$v_p$ is the valuation of $\cO_{E,p}$ (and $\cO_{R, p}$ respectively)
and
$r_i$ are weights of $\rho$. By Lemma~\ref{L:balance}, $e_\rho(E)_p \ge 4^2$
and $e_\rho(R)_p \ge 4^2$. On $R$, $v_q(x_i) + r_i \ge 2$ and
$e_\rho(R)_q \ge 2^2$. Since $\rho$ acts trivially on $D$ with
weight $2$, we use Lemma~\ref{L:degree} and obtain
\[
e_\rho(D) = 2\cdot 2\cdot \deg D = 4(4g-10).
\]
Combining them all, we obtain
\[
e_\rho(C) \ge 36 +  16g - 40 > 2\cdot \frac43\sum_{i=1}^{3g-4}r_i =
16g - \frac{40}3.
\]

\end{proof}

\begin{coro} Let $C' = E'\cup_p D$ be a bicanonical curve consisting
of a genus one curve $E'$ and a genus $g-2$ curve $D$ meeting in one
tacnode $p$. Then $C'$ is Chow unstable.
\end{coro}

\begin{proof} In view of Proposition~\ref{P:g1-1tac}, it suffices
to show that $C'$ is in the basin of attraction
 of $E\cup_p R\cup_q D$ with respect to $\rho$. Consider the
induced action on the local versal deformation space of the cusp
$
[1, 0, \dots, 0]
$
which is given by
\[
y^2 = x^3 + a x + b
\]
where $y = x_2/x_0$ and $x = x_1/x_0$. The $\bG_m$ action is given
by
\[
t.(a, b) = (t^4 a, t^6 b)
\]
and the basin of attraction contains arbitrary smoothing of the cusp.
On the other hand, the local versal deformation space of the tacnode $p = [0, 0, 0, 1, 0, \dots, 0]$ is given by
\[
y^2 = x^4 + a x^2 + b x + c
\] where $x = x_2/x_3$ so that $\bG_m$ acts on $(a, b, c)$
with weight $(-2, -3, -4)$ and the basin of attraction
does not contain any smoothings of the tacnode. At the node
$q = [0,0,0,0,1,0,\dots, 0]$, the local versal deformation space is
\[
x y = c_0
\]
where $x$ may be taken to be $x_2/x_4$
and $\bG_m$ acts with weight $+1$ on the
branch of $R$ and trivially on $D$. Thus the induced action
on the deformation space has weight $+1$, and the basin
of attraction contains arbitrary smoothing of the node.
\end{proof}

\subsection{Hilbert unstable curves}\label{S:Hilb-unstable}
Let $C$ be a bicanonical curve.
By Proposition~\ref{prop:HilbtoChow}, $C$ is Chow semistable if
it is Hilbert semistable.   Note that by definition,
 if $C$ does not admit an elliptic
chain, then $C$ is c-semistable if and only if it is h-semistable.
Combining this with Proposition~\ref{P:chow-imply-c}, we obtain:

\begin{prop}\label{P:hilb-imply-h}
If a bicanonical curve is Hilbert semistable and does not
admit an elliptic chain, then it is h-semistable.
\end{prop}

We shall have completed the implication
\[
\mbox{Hilbert semistable} \Rightarrow \mbox{h-semistable}
\]
once we prove that a Hilbert semistable curve
does not admit an elliptic chain. We accomplish this in 
Proposition~\ref{P:echain-unstable} and Corollary~\ref{C:basin-cr-1br}.

\section{Classification of curves with automorphisms}
\label{S:automorphisms}

In this section, we classify
c-semistable curves with infinite automorphisms.

\subsection{Rosaries}
\begin{defn}\label{D:rosary}
An {\it open rosary}\footnote{This name was suggested
to us by Jamie Song.} $R_r$ of length $r$ is a
two-pointed connected curve $(R_r, p, q)$
such that
 \begin{itemize}
 \item $R_r = L_1\cup_{a_1}L_2\cup_{a_2}\cdots\cup_{a_{r-1}}
L_r$ where $L_i$ is a smooth rational curve,
 $i = 1, \dots, r$;
 \item   $L_i$ and $L_{i+1}$ meet each other in a single
tacnode $a_i$, for $i = 1, \dots, r-1$;
\item $L_i\cap L_j = \emptyset$ if $|i-j| > 1$;
 \item $p \in L_1$ and $q \in L_r$ are smooth points.
\end{itemize}
\end{defn}
\begin{rem}
An open rosary of length $r$ has arithmetic genus $r-1$.
Note that an open rosary
of length $r=2r'$ is naturally an open elliptic chain of length $r'$.
\end{rem}

\begin{defn} We say that a curve $C$ {\it admits an open rosary}
or length $r$ if there is a $2$-pointed open rosary $(R_r, p, q)$ and
a morphism $\iota : R_r \to C$ such that
\begin{itemize}
\item $\iota$ is an isomorphism onto its image over $R_r\setminus\{p,q\}$;
\item $\iota(p), \iota(q)$ are nodes of $C$;  we allow the case $\iota(p) = \iota(q)$.
\end{itemize}
A {\em closed rosary} $C$ is a curve admitting $\iota : C' \to C$ as above with the second condition replaced by
\begin{itemize}
\item
$\iota(p)=\iota(q)$ at a tacnode of $C$.
\end{itemize}
A {\em closed rosary with broken beads} is a curve expressible as
a union of open rosaries.
\end{defn}
\begin{rem} \label{rem:rosary}
If $C$ admits an open rosary of length $r\ge 2$ then $C$ admits a
weak elliptic chain.  If $r$ is even then $C$ admits an elliptic chain.
Thus a closed rosary of even length is also a closed weak elliptic chain.
\end{rem}

\begin{figure}[!htb]
  \begin{center}
    \includegraphics[width=3.5in]{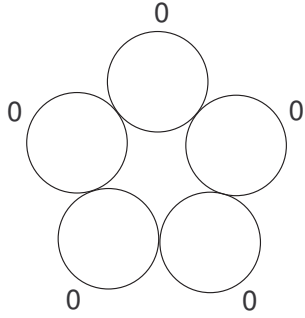}
  \end{center}

  \caption{Closed rosary of genus six}
  \label{F:rosary}
\end{figure}

\begin{prop}\label{P:rosary-aut}
Consider the closed rosaries of genus $r+1$.
If the genus is even then there is a unique closed rosary $C$
(of the given genus) 
and the automorphism group $Aut(C)$ is finite.  If the genus
is odd then the closed rosaries depend
on one modulus and the connected component
of the identity $Aut(C)^{\circ}$ is isomorphic
to $\bG_m$.

There is a unique open rosary
$(R,p,q)$ of length $r$.
If $Aut(R,p,q)$ denotes the automorphisms fixing $p$
and $q$  then
$$Aut(R,p,q)^{\circ} \simeq \bG_m.$$
It acts on tangent spaces of the endpoints with weights
satisfying
$$wt_{\bG_m}(T_pR)=(-1)^r wt_{\bG_m}(T_qR).$$
\end{prop}

\begin{proof}
Let $C$ be a closed $r$-rosary obtained by gluing $r$ smooth
rational curves $\{[s_i, t_i]\}$  so that
\[
\frac{\partial}{\partial (s_r/t_r)} = \a_r
\frac{\partial}{\partial (t_{1}/s_{1})}; \quad
\frac{\partial}{\partial (s_i/t_i)} = \a_i
\frac{\partial}{\partial (t_{i+1}/s_{i+1})}, \quad i = 1, 2,
\dots, r-1. \]
 Let $C'$ be another such rosary with the gluing data
\[
\frac{\partial}{\partial (s_r'/t_r')} = \a_r'
\frac{\partial}{\partial (t_{1}'/s_{1}')}; \quad
\frac{\partial}{\partial (s_i'/t_i')} = \a_i'
\frac{\partial}{\partial (t_{i+1}'/s_{i+1}')}, \quad i = 1, 2,
\dots, r-1. \]

 Consider the morphism $f : \til{C} \to \til{C}'$ between the
 normalizations of $C$ and $C'$
given by
$
[s_i, t_i] \mapsto [\b_i s_i', t_i'].$
For $f$ to descend to an isomorphism from $C$ to $C'$, the
following is necessary and sufficient:

\SMALL
\[
df\left(\frac{\partial}{\partial (s_i/t_i)}\right)
= \frac{\partial}{ \b_i\partial (s_i'/t_i')}
 =  \frac{\a_i'}{\b_i} \frac{\partial}{\partial
(t_{i+1}'/s_{i+1}')}  =  \a_i\b_{i+1}\frac{\partial}{\partial
(t'_{i+1}/s'_{i+1})} = df\left(\a_i\frac{\partial}{\partial (t_{i+1}/s_{i+1})}\right)
\]
\normalsize
 This gives rise to $
\b_i\b_{i+1} = \a'_i/\a_i$ and $ \b_r\b_1 = \a_r'/\a_r$.
Solving for $\b_i$, we get
\[
\b_i = \begin{cases} \frac{\a_i'\a_{i+1}\a_{i+2}' \cdots
\a_r}{\a_i\a'_{i+1}\a_{i+2} \cdots \a'_r}\b_1, \quad \mbox{if
$r-i$ is odd} \\
\ \\
\frac{\a_i'\a_{i+1}\a_{i+2}' \cdots \a_r'}{\a_i\a'_{i+1}\a_{i+2}
\cdots \a_r}\b_1\inv, \quad \mbox{if $r-i$ is even} \end{cases}
\]

When $r$ is odd, there is no constraint and all $r$-rosaries are
isomorphic.  When $r = 2k$,
\[
(\b_1\b_2)(\b_3\b_4)\cdots(\b_{2k-1}\b_{2k}) =
(\b_2\b_3)(\b_4\b_5)\cdots(\b_{2k}\b_1)
\]
forces the condition
\begin{equation}\label{E:para}
\frac{\a_1'\a_3'\cdots\a_{2k-1}'}{\a_1\a_3\cdots\a_{2k-1}} =
\frac{\a_2'\a_4'\cdots\a_{2k}'}{\a_2\a_4\cdots\a_{2k}}.
\end{equation}
This means that the $2k$-rosaries are parametrized by
\[
\frac{\a_1\a_3\cdots\a_{2k-1}}{\a_2\a_4\cdots\a_{2k}} \in \bG_m.
\]

To describe the automorphisms we take $C'=C$.  When $r$ is odd we get
$\beta_i=\beta^{-1}_i$ for each $i$ which implies that $Aut(C)^{\circ}$
is trivial.  When $r=2k$ we get a unique solution
$$\b_1=\b_2^{-1}=\b_3=\ldots=\b_{2k}^{-1}$$
and thus $Aut(C)^{\circ} \simeq \bG_m$.

The open rosary case entails exactly the same analysis, except that
we omit the gluing datum
\[
\frac{\partial}{\partial (s_r/t_r)} = \a_r
\frac{\partial}{\partial (t_{1}/s_{1})}
\]
associated with the end points.  Thus we get a $\bG_m$-action regardless of
the parity of $r$.  Our assertion on the weights at the distinguished points $p$ and $q$ 
follows from the computation above of the action on tangent spaces.
\end{proof}

\begin{defn} By \emph{breaking the $i$th bead} of a rosary (open or closed),
we mean replacing $L_i$ with a union $L_i'\cup L_i''$ of smooth
rational curves meeting in a node such that $L_i'$ meets $L_{i-1}$
in a tacnode $a_{i-1}$ and $L_i''$ meets $L_{i+1}$ in a tacnode
$a_{i+1}$ (Figure~\ref{F:breaking-bead}).

\begin{figure}[!htb]
  \begin{center}
    \includegraphics[width=5in]{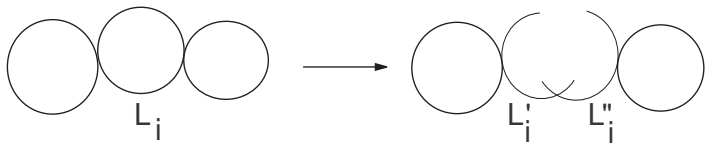}
  \end{center}

  \caption{Breaking a bead of a rosary}
  \label{F:breaking-bead}
\end{figure}
\end{defn}

\subsection{Classification of automorphisms}

\begin{prop}\label{P:hsinf} A c-semistable curve $C$ of genus $ \ge 4$
has infinite automorphisms if and only if
\begin{enumerate}
\item $C$ admits an open rosary of length $\ge 2$, or
\item $C$ is a closed rosary of odd genus (possibly with broken
beads).
\end{enumerate}
\end{prop}

\begin{proof} We have already seen in Proposition~\ref{P:rosary-aut}
that closed rosaries of odd
genus have infinite automorphisms.

Let $C$ be a c-semistable curve of genus $g
\ge 4$ that is not a closed rosary. For $C$ to have infinitely
many automorphisms, it must have a smooth rational component, say
$C_1$. To satisfy the stability condition and still give rise to
infinite automorphisms,
$C_1$ has to meet the
rest of the curve in one node and a tacnode, or in two tacnodes.
We examine each case below:

\smallskip

\ndt (1) \, $C_1$ meets the rest in one node $a_0$ and in a tacnode
$a_1$: For the automorphisms of $C_1$ to extend to automorphisms of
$C$, the irreducible component $C_2 (\ne C_1)$ containing $a_1$ must
be a smooth rational component - this follows easily from  that an
automorphism of $C$ lifts to an automorphism of its normalization.
Also, $C_2$ has to meet the rest of the curve in one point $a_2$
other than $a_1$ since otherwise $C_1\cup C_2$ would be an elliptic
tail (or $a_1 = a_0$ and $C$ is of genus two).

\smallskip

\ndt (2) \,  $C_1$ meets the rest in two tacnodes $a_0$ and $a_1$: For
the automorphisms to extend to $C$, the components $C_0 \ne C_1$
containing $a_0$ and  $C_2\ne C_1$ containing $a_1$ must be smooth
rational curves. Hence $C$ contains $C_0\cup C_1\cup C_2$ which is
a rosary of length three. Moreover, $C_0$ and $C_2$ do not intersect: If
they do meet, say at $a_2$, then either $C = C_0\cup C_1\cup C_2$
and the genus of $C$ is of genus three (if $a_2$ is a node) or $C$
is a closed rosary if $a_2$ is a tacnode.

Iterating, we eventually produce an open rosary $\iota:R_r \ra C$
of length $r\ge 2$ containing $C_1$ as a bead.
\end{proof}

\begin{coro}\label{C:hsinf1}
An h-semistable curve $C$ of genus $\ge 4$ has infinite
automorphisms if and only if
\begin{enumerate}
\item $C$ admits an open rosary of odd length $\ge 3$, or
\item $C$ is a closed rosary of odd genus (possibly with broken beads).
\end{enumerate}
\end{coro}

Let $C$ be a c-semistable curve and suppose
$D$ is a Deligne-Mumford stabilization of $C$.  In other words,
there exists a smoothing of $C$
$$\varpi:\cC \ra T$$
such that $D=\lim_{t\ra t_0}\cC_t$ in the moduli space of stable curves.
Here, a smoothing is a flat proper
morphism to a smooth curve with distinguished point $(T,t_0)$
such that $\varpi^{-1}(t_0)=C$ and the generic fiber is smooth.

Our classification result (Proposition~\ref{P:hsinf}) has the
following immediate consequence:
\begin{coro} \label{C:hsinf}
Suppose $C$ is a c-semistable curve with infinite automorphism group.
Then the Deligne-Mumford stabilization and the pseudo-stabilization of $C$
 admit an elliptic bridge.
\end{coro}
Indeed, $C$ necessarily admits a tacnode which means that its
stabilization contains a connected subcurve of genus one
meeting the rest of the curve in two points.

\section{Interpreting the flip via GIT}
\label{S:FGITQ}
We will eventually give a complete description of the semistable
and stable points of $\Cg2$ and $\Hg2$.  For our immediate purpose,
the following partial result will suffice:
\begin{thm} \label{T:partial}
If $C$ is c-stable, i.e., a pseudostable curve admitting no
 elliptic bridges,
then $Ch(C) \in \Cg2^s$.  Thus the Hilbert point
$[C]_m \in \Hg2^{s,m}$ for $m\gg 0$.
\end{thm}
\begin{proof}
The GIT-stable loci $\Cg2^s$ and $\Hg2^{s,m},m\gg 0$,
contain the nonsingular curves by \cite[4.15]{M}.
(We discussed the relation between Chow and asymptotic
Hilbert stability in Section~\ref{subsect:CPHS}.)

Recall that Proposition~\ref{P:embedcstable} guarantees
that c-semistable curves admits bicanonical embeddings.
In particular, this applies to pseudostable curves without
elliptic bridges.

Suppose that $C$ is a singular pseudostable curve
without elliptic bridges.  Assume that $Ch(C)$ is
{\em not} in $\Cg2^s$.   If $Ch(C)$ is strictly semistable then
it is c-equivalent to a semistable curve $C'$ with infinite automorphism
group.  It follows that $C$ is a pseudo-stabilization
of $C'$, and we get a contradiction to Corollary~\ref{C:hsinf}.
Suppose $Ch(C)$ is unstable and let $C'$ denote a semistable replacement
(see Theorem~\ref{T:ss-replacement}).  By uniqueness of the
pseudo-stabilization, $C'$ is not pseudostable
but has $C$ as its pseudo-stabilization.
It follows that $C'$ has a tacnode.  However,
the pseudo-stabilization of such a curve necessarily
contains an  elliptic bridge.
\end{proof}

With our current partial understanding of the GIT of bicanonical curves, we are ready
to prove Theorem~\ref{T:main12}.
Our main task is to establish Isomorphisms (\ref{eqn:Tmain1}) and
(\ref{eqn:Tmain2}).  Proposition~\ref{P:makePsi} established the existence
of a birational contraction morphism
$\Psi:\Mg^{ps} \ra  \Mg(7/10)$.
The first step here is to show that $\Mg^{cs}$ and $\Mg^{hs}$
are birational contractions of $\Mg$ and small contractions of $\Mg^{ps}$.
In particular, we may identify the divisor class groups of these GIT
quotients with the divisor class group of $\Mg^{ps}$ (which in turn is a
subgroup of the divisor class group of $\Mg$).  Furthermore, we obtain
$$
\Gamma(\Mps, n(K_{\Mps} + \alpha \dps))
\simeq \Gamma(\Mhs, n(K_{\Mhs} + \alpha \dhs))
\simeq \Gamma(\Mcs, n(K_{\Mcs} + \alpha \dcs))
$$
and Lemma~\ref{L:logterm} gives
\begin{equation} \label{E:smallproj}
\begin{array}{rcl}
\Mg(7/10)&\simeq & Proj \left( \oplus_{n \ge 0} \Gamma(n(K_{\Mcs}+7/10 \dcs)) \right) \\
\Mg(7/10-\epsilon) & \simeq & Proj
\left( \oplus_{n \ge 0} \Gamma(n(K_{\Mhs}+(7/10-\epsilon) \dhs)) \right).
\end{array}
\end{equation}

The second step is to compute
the induced polarizations of $\Mg^{cs}$ and $\Mg^{hs}$
in the divisor class group of $\Mg^{ps}$.  This
will show that $K_{\Mcs}+7/10\delta^{cs}$
(resp. $K_{\Mhs}+(7/10-\epsilon)\delta^{hs}$)
is ample on $\Mg^{cs}$ (resp. $\Mg^{hs}$).
Isomorphisms (\ref{eqn:Tmain1}) and
(\ref{eqn:Tmain2}) then follow from
(\ref{E:smallproj}).

To realize our GIT quotients as contractions of  $\Mg$, we apply Theorem~\ref{thm:getcontraction}
in the bicanonical case.
Consider the complement of
the Deligne-Mumford stable curves $V_{g,2}^{ss,\infty}$ in the GIT-semistable locus $\Cg2^{ss}$;
we must show this has codimension $\ge 2$.
Since $\cycle(\Hg2^{ss,m}) \subset \Cg2^{ss}$
and $\cycle|\Hg2^{s,m}$ is an isomorphism where
$\cycle$ denotes the cycle class map from $Hilb$ to $Chow$, the analogous statement for the Hilbert scheme follows
immediately.

Proposition~\ref{P:chow-imply-c} implies  $\Cg2^{ss} \setminus V_{g,2}^{ss,\infty}$
parametrizes
\begin{itemize}
\item{pseudostable curves that are not Deligne-Mumford stable, i.e., those with cusps; and}
\item{c-semistable curves with tacnodes.}
\end{itemize}
The cuspidal pseudostable curves have codimension two in moduli;  the tacnodal curves
have codimension three.  Indeed, a generic tacnodal curve of genus $g$ is determined by a
two-pointed curve $(C', p, q)$ of genus $g-2$ and an isomorphism
$T_pC' \simeq T_qC'$.
We conclude there exist rational contractions $F^{cs}:\Mg \dashrightarrow \Mg^{cs}$
and $F^{hs}:\Mg \dashrightarrow \Mg^{hs}$.

It remains to show that we have small contractions
$G^{cs}:\Mg^{ps} \dashrightarrow \Mg^{cs}$
and $G^{hs}:\Mg^{ps} \dashrightarrow \Mg^{hs}$.
To achieve this, we must establish that $\Delta_1$
is the unique exceptional divisor of $F^{cs}$ (resp. $F^{hs}$).
The exceptional locus of $F^{cs}$ (resp. $F^{hs}$) lies in the complement
to the GIT-stable curves in the moduli space
$$\Mg \setminus U_{g,2}^{s,\infty} \quad (\text{resp.} \,
\Mg \setminus U_{g,2}^{s,m}).$$
Chow stable points are asymptotically Hilbert stable
(cf. Proposition~\ref{prop:HilbtoChow}), i.e.,
$U_{g,2}^{s,\infty} \subset U_{g,2}^{s,m}$ when $m \gg 0$.
It suffices then to observe that $\delta_1$ is the unique divisorial
component of
$\Mg \setminus U_{g,2}^{s,\infty}$, which is guaranteed by Theorem~\ref{T:partial}.

Theorem~\ref{thm:getcontraction} gives moving divisors on $\Mg$
inducing the contractions $F^{cs}$ and $F^{hs}$.  We apply Equation (\ref{E:Lm}),
which in our situation takes the form
$$
r(2)\lambda_{m2}-r(2m)m\lambda_2 =
            (m-1)(g-1)((20m-3)\lambda - 2m \delta) \sim (10-\frac{3}{2m})\lambda - \delta;$$
this approaches $10\lambda - \delta$ as $m\to \infty$.
Thus we have
$$(F^{hs})^*{\mcl L}_m \sim (10-\frac{3}{2m})\lambda - \delta\pmod{\delta_1}, \quad  m \gg 0$$
and
$$(F^{cs})^*{\mcl L}_{\infty} \sim 10\lambda - \delta \pmod{\delta_1}.$$
Using the identity
$$K_{\FMg} = 13 \lambda - 2 \delta$$
we obtain (for $m\gg 0$)
$$(F^{hs})^*{\mcl L}_m \sim K_{\FMg} + (7/10 - \e(m)))\delta \pmod{\delta_1},\quad  \e(m)=39/(200m-30)$$
and
$$(F^{cs})^*{\mcl L}_{\infty} \sim K_{\FMg} + 7/10\delta \pmod{\delta_1}.$$
It follows then that
$$(G^{hs})^*{\mcl L}_m \sim K_{\Mps}+ (7/10 - \e(m)))\delta^{ps}$$
and
$$(G^{cs})^*{\mcl L}_{\infty} \sim K_{\Mps}+ 7/10\delta^{ps}.$$

The proof of Theorem~\ref{T:main12} will be complete if we can show
that $\Psi^+$ is the flip of $\Psi$.  More precisely, for small
positive $\e \in \bQ$, $\Psi^+$ is a
small modification of $\Mps$ with
$K_{\Mhs}+(7/10-\e)\dhs$ ample.
Since $\Mcs$ and $\Mhs$ ar
both small contractions of $\Mps$, $\Psi^+$ is small as well.
And the polarization we exhibited on $\Mhs$ gives the desired
positivity, which completes the proof of Theorem~\ref{T:main12}.


\section{Stability under one-parameter subgroups}
\label{S:hardwork}

In this section, we analyze whether c-semistable curves
are GIT-semistable with respect to the one-parameter
subgroups of their
automorphism group.
We shall also use deformation theory to
classify the curves that belong to basins of attraction of such
curves.

Our analysis will focus primarily on the Hilbert points.
Indeed, Corollary~\ref{C:all-m} shows that we can recover
the sign of the Hilbert-Mumford index of the Chow point from
the indices of the Hilbert points.  And
in view of the cycle map $\cycle: Hilb_{g,2} \to Chow_{g,2}$,
if $[C]_m \in A_{\rho}([C^{\star}]_m)$ for $m \gg 0$ then
$Ch(C) \in A_\rho(Ch(C^\star))$.

\subsection{Stability analysis: Open rosaries}\label{S:open-rosaries}
\begin{prop}\label{P:hs-o-ros}
Let $C = D\cup_{a_0,a_{r+1}}R$ be a c-semistable curve of genus $g$
consisting of a genus $g-r-1$ curve $D$ meeting the genus $r$
curve $R$ in two nodes $a_0$ and $a_{r+1}$ where
$$ R := L_1\cup_{a_1} L_2 \cup_{a_2}  \cdots \cup_{a_{r}} L_{r+1}$$
is a rosary of length $r+1$, and $D\cap L_1 = \{a_0\}$ and $D \cap
L_{r+1} = \{a_{r+1}\}$. There is a one-parameter subgroup $\rho$
coming from the automorphisms of $C$ of the rosary $R$ such that for all $m \ge 2$,
\begin{enumerate}
\item  $\mu([C]_m,\rho) = 0$ if $r$ is even;
\item $\mu([C]_m,\rho) = -m + 1$ if $r$ is odd. \end{enumerate}
In particular, $C$ is Hilbert unstable if $R$ is of even
length and strictly semistable otherwise.
\end{prop}
\noindent An application of Corollary~\ref{C:all-m} then yields:
\begin{coro} \label{C:cs-o-ros} Let $C$ and $\rho$ be as in
Proposition~\ref{P:hs-o-ros}. Then $C$ is Chow strictly semistable
with respect to $\rho$ and $\rho^{-1}$.
\end{coro}

\begin{proof} Upon restricting $\o_C$ to $D$ and each component
of $L$, we get
\begin{itemize}
\item $\o_C|_{D} \simeq \o_D(a_0 + a_{r+1})$;
\item $\o_C|_{L_1} \simeq \o_{L_1}(a_0 + 2 a_1)$;
\item $\o_C|_{L_{r+1}} \simeq \o_{L_{r+1}}(a_{r+1} + 2 a_{r})$;
\item $\o_C|_{L_i} \simeq \o_{L_i}(2 a_{i-1} + 2 a_i)$, \quad $2
\le i \le r$.
\end{itemize}
Hence we may choose coordinates $x_0, \dots, x_N$, $N = 3g-4$,
such that
\begin{enumerate}
\item $L_1$ is parametrized by $$[s_1, t_1] \mapsto
[s_1^2, s_1t_1, t_1^2, 0, \dots, 0];$$
\item $L_{r+1}$ is parametrized by
$$[s_{r+1}, t_{r+1}]\mapsto [\underbrace{0, \dots, 0}_{3r-2},
s_{r+1}t_{r+1}, s_{r+1}^2, t_{r+1}^2, 0, \dots, 0].$$
\item For $2\le j \le r$, $L_j$ is parametrized by
$$[s_j,t_j] \mapsto [ \underbrace{0, \dots, 0}_{3j-5}, s_j^3t_j, s_j^4, s_j^2t_j^2,
s_jt_j^3, t_j^4, 0, \dots, 0].$$
\item $D$ is contained in the linear subspace
$$ x_1 = x_2 = \cdots = x_{3r-1} = 0$$
and $a_0 = [1, 0, \dots, 0]$ and $a_{r+1} = [\underbrace{0, \dots,
0}_{3r}, 1, 0, \dots, 0]$.
\end{enumerate}

From the parametrization, we obtain a set of generators for the
ideal of $L$:
\begin{equation}\label{E:gen1}
\begin{array}{l}
x_1^2 - x_0x_2 - x_2x_3, x_0x_3, x_0x_4, \dots, x_0x_{3r}, \\
x_ix_{i+5}, x_ix_{i+6}, \dots, x_ix_{3r}, \quad i = 1, 2, \dots, 3r-5, \\
\end{array}
\end{equation}
and for $j = 1, 2, \dots, r-1$,
\begin{equation}\label{E:gen2}
\begin{array}{c}
x_{3j-1}x_{3j+3}, \, x_{3j}x_{3j+3}, \, x_{3j}x_{3j+4}, \,
x_{3j+1}^2 -
x_{3j}x_{3j+2} - x_{3j+2}x_{3j+3}, \, x_{3j}^2 - x_{3j-1}x_{3j+2}, \, \\
x_{3j-2}x_{3j+1} - x_{3j-1}x_{3j+2}, \,  x_{3j-2}x_{3j} -
x_{3j-1}x_{3j+1}, \, x_{3j-2}x_{3j+2} - x_{3j}x_{3j+1}.
\end{array}
\end{equation}
In Proposition~\ref{P:rosary-aut}
we showed that $\bG_m$ acts
on the open rosary via automorphisms.
With respect to our coordinates, this is
the one-parameter subgroup $\rho$ with weights:
\[
\begin{cases}
(2,1,0,2,3,4,2,1,0, \dots, 2,3,4, \underbrace{2, 2, \dots, 2}_{N-3r}), \quad \mbox{if $r$ is even} \\
(2,1,0,2,3,4,2,1,0, \dots, 2,1,0, \overbrace{2, 2, \dots, 2}),
\quad \mbox{if $r$ is odd.}
\end{cases}
\]
By considering the parametrization, it is easy to see that $C$ is stable
under the action of $\rho$.

Now we shall enumerate the degree two monomials in the initial
ideal of $C$. From (\ref{E:gen1}) and (\ref{E:gen2}), we get the
following monomials in $x_0, \dots, x_{3r}$:
\begin{equation}\label{E:inL-deg2}
\begin{array}{l}
x_0x_2, x_0x_3, x_0x_4, \dots, x_0x_{3r},  \\
x_ix_{i+5}, x_ix_{i+6}, \dots, x_ix_{3r}, \quad i = 1, 2, \dots,
3r-5, \\
x_{3j-1}x_{3j+3}, x_{3j}x_{3j+3}, x_{3j}x_{3j+4}, x_{3j}x_{3j+2},
\\ x_{3j-1}x_{3j+2}, x_{3j-2}x_{3j+1}, x_{3j-2}x_{3j},
x_{3j-2}x_{3j+2}, \quad  j = 1, 2,
\dots, r-1.\\
\end{array}
\end{equation}
The weights of these $(9r^2-5r)/2$ monomials sum up to give
\[
\begin{cases}
 18r^2 - 10r, \quad \mbox{if $r$ is even;} \\
18r^2 - 19 r + 7, \quad \mbox{if $r$ is odd}.
\end{cases}
\]
The total weight $\displaystyle{\sum_{i\le j, \, 0\le i,j \le 3r}
wt_\rho(x_ix_j) } $ of all degree two monomials in $x_0, \dots,
x_{3r}$ is
\[
\begin{cases}
(3r+2)(6r+2), \quad \mbox{if $r$ is even;} \\
(3r+2)(6r-1), \quad \mbox{if $r$ is odd}.
\end{cases}
\]
Hence the degree two monomials in $x_0, \dots, x_{3r}$ that are
{\it not} in the initial ideal contributes, to the total weight,
\begin{equation}\label{E:L-deg2}
\begin{cases}
 (3r+2)(6r+2) - (18r^2 - 10r) = 28r + 4,  \quad \mbox{if $r$ is even;} \\
 (3r+2)(6r-1) - (18r^2 - 19 r + 7) = 28r - 9,\quad \mbox{if $r$ is odd}.
\end{cases}
\end{equation}
The rest of the contribution comes from the monomials supported on
the component $D$: These are the degree two monomials in $x_0,
x_{3r}, x_{3r+1}, \dots, x_N$ that vanish at $a_0$ and $a_{r+2}$.
The number of such monomials is,
 by Riemann-Roch,
\begin{equation}\label{E:D-deg2} h^0(D, \cO_{D}(2)(-a_0 - a_{r+2})) =
7(g-r-1)-1.
\end{equation}
 Since $wt_\rho(x_i) = 2$ for all $i = 0, 3r, 3r+1, \dots, N$, these monomials contribute
 $28g - 28r-32$ to the sum. Combining (\ref{E:L-deg2}) and (\ref{E:D-deg2}), we find
 the sum of the weights of the degree two monomials not in $in(C)$ to be
\[
\begin{cases}
28g - 28, \quad \mbox{if $r$ is even;} \\
28g - 41,  \quad \mbox{if $r$ is odd}.
\end{cases}
\]
On the other hand, the average weight is
\[
\frac{2 \cdot P(2)}{N+1} \sum_{i=0}^N
wt_\rho(x_i) = \begin{cases}
28g - 28, \quad \mbox{if $r$ is even;} \\
28g - 42,  \quad \mbox{if $r$ is odd}.
\end{cases}
\]
Hence by (\ref{eqn:stabcrit}), we find that $\mu([C]_2, \rho) = 0$ if $r$ is even and $\mu([C]_2,\rho) = -1$ if $r$ is odd.

We enumerate the degree three monomials in the same way: The
degree three monomials in $x_0, \dots, x_{3r}$ that are in the
initial ideal are the multiples of (\ref{E:inL-deg2}) together
with
\begin{equation}
x_{3j+1}^2x_{3j+3}, \quad j = 0, 1, \dots, r-1,
\end{equation}
that come from the linear relation
\[
x_{3j+2}(x_{3j}x_{3j+3}) - x_{3j+3}(x_{3j+1}^2 - x_{3j}x_{3j+2} -
x_{3j+2}x_{3j+3}) = 0.
\]
From this, we find that the degree three monomials in $x_0, \dots,
x_{3r}$ that are not in the initial ideal contribute
\[
\begin{cases}
66r+6, \quad \mbox{if $r$ is even;} \\
66r - 25,  \quad \mbox{if $r$ is odd}.
\end{cases}
\]
The contribution from $D$ is
\[
6 h^0(D, \cO_D(3)(-a_0-a_{r+2})) = 6 (11g - 11r - 12)
= 66g - 66r - 72.
\]
Hence the grand total is
\[
\sum_{j=1}^{P(3)} wt_\rho(x^{a(j)}) =
\begin{cases}
66g - 66, \quad \mbox{if $r$ is even;} \\
66g - 97,  \quad \mbox{if $r$ is odd}
\end{cases}
\]
On the other hand, the average weight is
\[
\frac{3P(3)}{N+1} \sum_{i=0}^N
wt_\rho(x_i) =
\begin{cases}
66g-66, \quad \mbox{if $r$ is even;} \\
66g-99,  \quad \mbox{if $r$ is odd}.
\end{cases}
\]
Using (\ref{eqn:stabcrit}), we compute
\[
\mu([C]_3, \rho) = \begin{cases}
0, \quad \mbox{ if $r$ is even;}\\
-2, \quad \mbox{if $r$ is odd.}
\end{cases}
\]
Corollary~\ref{C:all-m} implies
\[
\mu([C]_m,\rho) = \begin{cases}0, \quad \mbox{ if $r$ is even;}\\
-m+1, \quad \mbox{if $r$ is odd.}
\end{cases}
\]
for each $m\ge 2$.
\end{proof}


\subsection{Basin of attraction: Open rosaries}\label{S:ba-o-ros}
Let $C$ and $R$ be as in the previous section.  Let $x_i, y_i$ be homogeneous coordinates on $L_i$. We may assume  that
\[
a_0 = [0, 1]; \, a_{r+1} = [1, 0]; \,
a_i = \begin{cases} [1,0] = \infty \, \, {\mbox on} \, \, L_i\\
[0,1] = 0 \, \, {\mbox on} \, \, L_{i+1} \\
\end{cases} \\
\]
Consider the $\bG_m$ action associated to $\rho$.  The action on
$R$ is given by $(t, [x_i,y_i]) \mapsto [x_i,t^{(-1)^{i-1}}
y_i]$ on each $L_i$. Hence it induces an action on the tangent space
$T_{a_i}L_i$ given by
\[
\displaystyle{\left(t, \frac{\partial}{\partial (y_i/x_i)}\right)
\mapsto
\frac{\partial}{\partial\left(t^{(-1)^{i-1}}y_i/x_i\right)} =
t^{(-1)^i} \frac{\partial}{\partial(y_i/x_i)}}.
\]

There is an induced
$\bG_m$ action on the Hilbert scheme and $\Hg2$.  Corollary~\ref{C:versal}
asserts that a neighborhood of $[C]$
in the Hilbert scheme dominates the product of the versal deformation
spaces.  These inherit a $\bG_m$ action as well, which we shall
compute explicitly.

\ndt (A) \, \textit{$\bG_m$ action on the versal deformation spaces
of nodes $a_0$ and $a_r$}: Let $z$ be a local parameter at $a_0$ on
$D$. We have $x_1/y_1$ as a local parameter at $a_0$ on $L_0$ and
the local equation at $a_0$ on  $C$ is $z\cdot (x_1/y_1) = 0$. Hence
the action on the node $a_0$ is given by $(z, x_1/y_1) \mapsto (z,
t^{-1} x_1/y_1)$ and the action on the versal deformation space is
\[
c_0 \mapsto t^{-1} c_0.
\]
Likewise, at $a_{r+1}$, the action on the node is
\[
(y_{r+1}/x_{r+1}, z') \mapsto (t^{(-1)^{r+1}} y_{r+1}/x_{r+1}, z')
\]
where $z'$ is a local parameter at $a_{r+1}$ on $D$,
and the action on the versal deformation space is
\[
c_0 \mapsto t^{(-1)^{r}} c_0.
\]

\

\ndt (B) \, \textit{$\bG_m$ action on the versal deformation space
of a tacnode $a_i$}: At $a_i$, the local analytic equation is of
the form $y^2 = x^4$ where $x := (y_i/x_i, x_{i+1}/y_{i+1})$ and
$y := ((y_i/x_i)^2, -(x_{i+1}/y_{i+1})^2)$ in $ k[[y_i/x_i]]\oplus
k[[x_{i+1}/y_{i+1}]]$
 and the $\bG_m$ action at the tacnode is
given by
\[
\begin{array}{c}
t. x = (t^{(-1)^{i-1}} y_i/x_i, x_{i+1}/(t^{(-1)^i}y_{i+1})) =
t^{(-1)^{i-1}} x\\
t. y = t^{2 (-1)^{i-1}} y.
\end{array}
\]
Therefore the action on the versal deformation space is
\[
(c_0,c_1,c_2) \mapsto (t^{4(-1)^{i-1}} c_0, t^{3(-1)^{i-1}} c_1,
t^{2(-1)^{i-1}} c_2).
\]

\

\begin{figure}[!t]
  \begin{center}
    \includegraphics[width=6.5in]{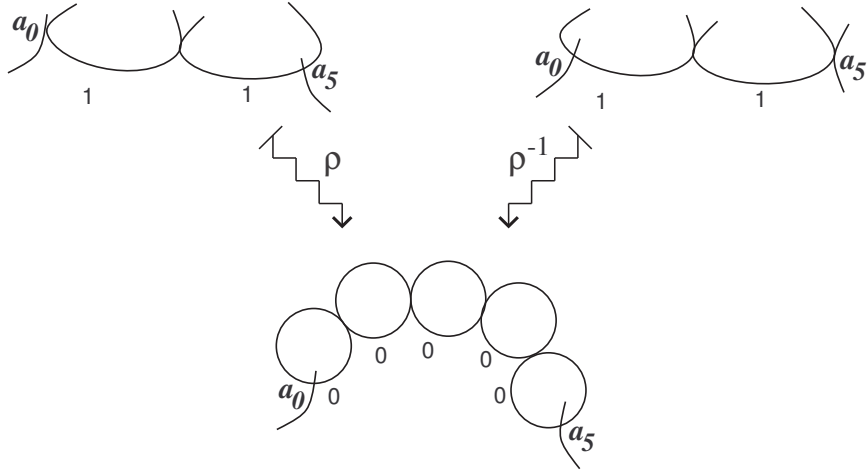}
  \end{center}

  \caption{Basin of attraction of an open rosary of length five}
  \label{F:basin-o-ros-even1}
\end{figure}

\begin{figure}[!t]
  \begin{center}
    \includegraphics[width=6in]{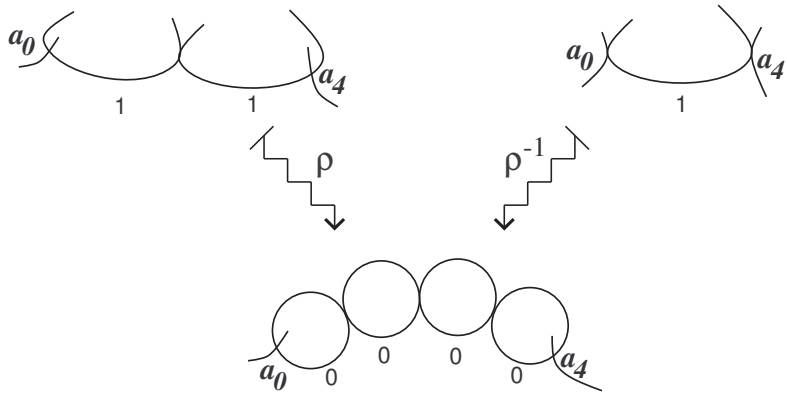}
  \end{center}

  \caption{Basin of attraction of an open rosary of length four}
  \label{F:basin-o-ros-odd}
\end{figure}

\ndt From these observations, we conclude that the basin of
attraction of $C$ with respect to $\rho$ contains
arbitrary smoothings of $a_{2k+1}$ but no smoothing of $a_{2k}$ for
all $0 \le k < \lceil {(r+1)}/2 \rceil$. We have established:

\begin{prop} \label{P:hs-echain}
Retain the notation of Proposition~\ref{P:hs-o-ros}
and assume that $m\gg 0$.
\begin{enumerate}
\item
If $r$ is even (i.e., the length of the rosary is odd)  then $A_\rho([C]_m)$ (resp.
$A_{\rho^{-1}}([C]_m)$) parametrizes the curves
consisting of $D$ and a weak elliptic chain $C'$ of length $r/2$
meeting $D$ in a node at $a_0$ and in a tacnode at $a_{r+1}$ (resp.
in a tacnode at $a_0$ and in a node at $a_{r+1}$)
(Figure~\ref{F:basin-o-ros-even1});

\item If $r$ is odd (i.e., the length of the rosary is even) then $A_\rho([C]_m)$
(resp. $A_{\rho^{-1}}([C]_m)$)
parametrizes the curves consisting of $D$ and an
elliptic chain $C'$ of length $(r+1)/2$ (resp. length $(r-1)/2$)
meeting $D$ in a node (resp. tacnode) at $a_0$ and $a_{r+1}$
(Figure~\ref{F:basin-o-ros-odd}).  When $r=1$, $A_{\rho}([C]_m)$
consists of tacnodal curves normalized by $D$.
\end{enumerate}
\end{prop}

It follows from Proposition~\ref{P:hs-echain} and Proposition~\ref{P:hs-o-ros} that

\begin{prop}\label{P:echain-unstable}
If a bicanonical curve admits an open elliptic chain
then it is Hilbert unstable.
In particular, a bicanonical curve with an elliptic bridge is
Hilbert unstable.
\end{prop}
\noindent  The closed case can be found in Proposition~\ref{P:basin-cr-1br} and Corollary~\ref{C:basin-cr-1br}.

\medskip

\subsection{Stability analysis: Closed rosaries}
\begin{prop}\label{P:hss-c-ros}
Let $C$ be a bicanonical closed rosary of even length $r$.
Then $C$ is Hilbert strictly semistable with respect to the one-parameter
subgroup $\rho : \bG_m \to \SL(3r)$ arising from $Aut(C)$.
\end{prop}
The relevant one-parameter subgroup was introduced in
Proposition~\ref{P:rosary-aut}.

\begin{proof}
Restricting $\o_C^{\ten 2}$ to each component $L_i$, we find
that each $L_i$ is a smooth conic in $\bP^{3g-4}$. We can choose
coordinates $x_0, \dots, x_N$ such that $L_i$ is parametrized by
\begin{itemize}
\item  $[s_i, t_i] \mapsto [\underbrace{0, \dots, 0}_{3(i-1)}, s_i^3t_i, s_i^4,
s_i^2t_i^2, s_it_i^3, t_i^4, 0, \dots, 0]$, \quad $i= 1, \dots,
r-1$;
\item $[s_r, t_r] \mapsto
[s_rt_r^3, t_r^4,  0, \dots, 0, s_r^3t_r, s_r^4, s_r^2t_r^2]$
\end{itemize}
The normalization of $C$ admits the automorphisms given by
\[
[s_i, t_i] \mapsto [\a^{sgn(i)}s_i, \, \a^{1-sgn(i)}t_i], \quad
sgn(i) := i - 2\lfloor i/2 \rfloor,
\]
for $i = 1, \dots, r-1$ and $ [s_r, t_r] \mapsto [s_r, \a t_r] $.
The one-parameter subgroup $\rho$ associated to this automorphism
has weights
\[
(3,4,2,1,0,2,\cdots,3,4,2,1,0,2)
\]
The sum of the weights $\sum_{i=1}^N wt_\rho(x_i)$ is $6r$ if $r$
is even and $6r+3$ if $r$ is odd.

From the parametrization, we obtain a set of generators for the
ideal of $C$:
\begin{equation}\label{E:CRgen}
\begin{array}{l}
x_0x_5, x_0x_6, \dots, x_0x_{3r-4}, \, x_1x_5, x_1x_6, \dots, x_1x_{3r-4}, \\
x_ix_{i+5}, \, x_ix_{i+6}, \, \dots \, , \, x_ix_{3r-1}, \quad i = 2, \dots, 3r-6, \\
x_{3j-2}x_{3j+2}, x_{3j-1}x_{3j+2}, x_{3j-1}x_{3j+3},
x_{3j}x_{3j+2},  \\
x_0^2-x_1x_2-x_1x_{3r-1}, \, x_{3j}^2 - x_{3j-1}x_{3j+1}-x_{3j+1}x_{3j+2}, \\
x_{3r-1}^2 - x_1x_{3r-2}, \, x_{3j-1}^2 - x_{3j-2}x_{3j+1},  \\
x_0x_{3r-3}-x_{1}x_{3r-2}, \, x_{3j-3}x_{3j} - x_{3j-2}x_{3j+1}, \\
x_0x_{3r-1} - x_1x_{3r-3}, x_0x_{3r-2} - x_{3r-3}x_{3r-1}, \,
x_{3j-3}x_{3j-1} - x_{3j-2}x_{3j}, \quad j = 1, 2, \dots, r-1, \\
x_{3j}x_{3j+3}, x_{3j}x_{3j+4},  \quad j = 1, 2, \dots, r-2.
\end{array}
\end{equation}
We have the following $(9r^2-11r)/2$ degree two monomials  that
are in the initial ideal:
\begin{equation}\label{E:inCR}
\begin{array}{l}
x_0^2, x_0x_5, x_0x_6, \dots, x_0x_{3r-1},   \\
x_1x_5, x_1x_6,  \dots, x_1x_{3r-4}, x_1x_{3r-2}, \\
x_ix_{i+5}, \, x_ix_{i+6}, \, \dots \, , \, x_ix_{3r-1}, \quad i = 2, \dots, 3r-6, \\
x_{3j-3}x_{3j-1}, x_{3j-3}x_{3j}, x_{3j-2}x_{3j+1},
x_{3j-2}x_{3j+2},
x_{3j-1}x_{3j+1}, x_{3j-1}x_{3j+2}, \quad j = 1, 2, \dots, r-1, \\
 x_{3j-1}x_{3j+3},
x_{3j}x_{3j+4}, \quad j = 1, 2, \dots, r-2. \\
\end{array}
\end{equation}
The sum of the weights of these monomials is $ 34r-18r^2$. It
follows that the sum of the weights of the monomials {\it not} in
the initial ideal is
\[
(3r-1)6r - (34r-18r^2) = 28r
\]
which is precisely $\frac{2P(2)}{N+1}\sum_{i=0}^N wt_\rho(x_i).$
Hence $\mu([C]_2, \rho) = 0$.

We shall now enumerate the degree three monomials in the initial
ideal of $C$. Together with the monomials divisible by the
monomials from (\ref{E:inCR}), we have the  initial terms
\begin{equation}\label{E:inL-syz-deg3}
x_{3r-3}^2x_{3r-1}, \, x_1x_{3r-3}^2, \,  x_{3j-2}x_{3j}^2, \quad j
= 1, 2, \dots, r-1 \end{equation} that come from the Gr\"obner basis
members
\[
x_1x_{3r-3}^2 - x_{3r-1}^3, \, x_{3r-3}^2x_{3r-1} -
x_{3r-2}x_{3r-1}^2, \, x_{3j-2}x_{3j}^2 - x_{3j-1}^3, \quad j = 1,
2, \dots, r-1.
\]
The degree three monomials in (\ref{E:inL-syz-deg3}) and the
degree three monomials divisible by monomials in (\ref{E:inCR})
have total weight $66r$. This agrees with the average weight
$\frac{3P(3)}{N+1}\sum_{i=0}^N wt_\rho(x_i) = \frac{3\cdot
11(g-1)}{3g-3}\frac{3g-3}6(3+4+2+1+0+2) = 66(g-1) = 66r.$
Therefore, $\mu([C]_2, \rho) = 0 = \mu([C]_3, \rho) $ and $C$ is
$m$-Hilbert strictly semistable by Corollary~\ref{C:all-m}.

Since $Aut(C)\simeq \bG_m$, a one-parameter subgroup coming from $Aut(C)$ is of the form $\rho^{a}$ for some $a \in \bZ$, and we have
\[
\mu([C]_m, \rho^a) = a \, \mu([C]_m, \rho) = 0.
\]
\end{proof}

\subsection{Basin of attraction: Closed
rosaries}\label{S:basin-closed}

\begin{prop}\label{P:gen-ss-c} Retain the notation of
Proposition~\ref{P:hss-c-ros}.
Then the basin of attraction $A_\rho([C]_m)$ parametrizes
the closed weak elliptic chains of
length $r/2$ (Figure~\ref{F:basin-c-ros}).
\end{prop}
\begin{proof} We use the parametrization from the proof of
Proposition~\ref{P:hss-c-ros}.  $C$ has tacnodes $a_i =
[\underbrace{0, \dots, 0}_{3i+1}, 1, 0, \dots, 0]$, $i = 1, \dots,
k$. From the parametrization, we find that the local parameters
$x_{3i}/x_{3i+1}$ at $a_i$ to the two branches is acted upon by
$\rho$ with weight $(-1)^{i-1}$.
 It follows that $\rho$ acts on the versal
deformation space $(c_0,c_1,c_2)$ of the tacnode $a_i$ with weights
$(4(-1)^{i-1}, 3(-1)^{i-1}, 2(-1)^{i-1})$. Hence the basin of
attraction $A_\rho([C])$ has arbitrary smoothings of $a_i$ for odd
$i$ but no nontrivial deformations of $a_i$ for even $i$.

\begin{figure}[!t]
  \begin{center}
    \includegraphics[width=6in]{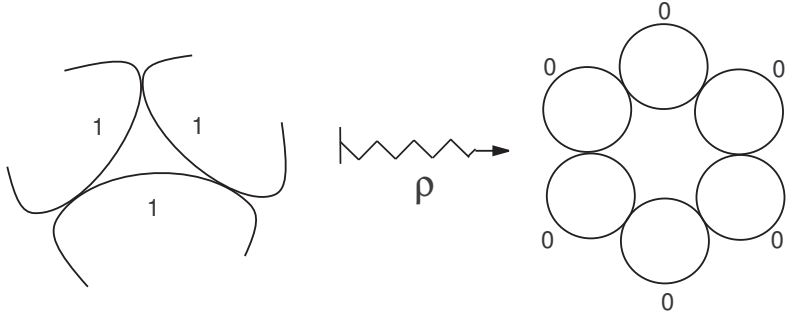}
  \end{center}
  \caption{Basin of attraction of a closed rosary of length six}
  \label{F:basin-c-ros}
\end{figure}
\end{proof}

\subsection{Stability analysis:  Closed rosaries with a broken bead}\label{S:rosary-broken}
Closed rosaries with broken beads of {\em even} genus
are unstable:
\begin{prop}
\label{P:hs-cr-1br} Let $r \ge 3$ be an odd number and
$C_r$ the curve obtained from a closed rosary of length $r$
by breaking a bead.  Then there exists a one-parameter subgroup
$\rho$ of $Aut(C_r)$ with
$\mu([C_r]_m,\rho) = 1-m $ for each $m\ge 2$, so
$[C_r]_m$ is unstable.  Furthermore, $Ch(C_r)$
is strictly semistable with respect to $\rho$.
\end{prop}

\begin{figure}[!t]
  \begin{center}
    \includegraphics[width=3.5in]{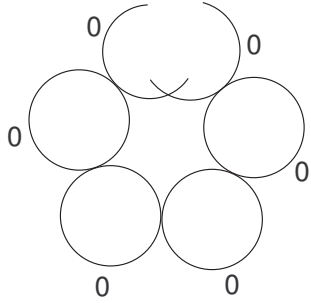}
  \end{center}

  \caption{Closed rosary of genus six with a broken bead}
  \label{F:broken-rosary}
\end{figure}

\begin{proof}
Note that $C_r$ is unique up to isomorphism
and it can parametrized by \small
\begin{equation}\label{E:CR1-para}
\begin{array}{lcl}
\bullet \, (s_0,t_0) & \mapsto  &(s_0t_0, s_0^2, t_0^2, 0, \dots, 0);\\
\bullet \, (s_1,t_1)&  \mapsto  & (0, 0, s_1^2, s_1t_1, t_1^2, 0, \dots, 0);\\
\bullet \, (s_i,t_i) & \mapsto & (\underbrace{0, \dots,
0}_{3(i-1)}, s_i^3t_i, s_i^4, s_i^2t_i^2, s_it_i^3,
t_i^4, 0, \dots, 0), \quad i = 2, \dots, r-1; \\
\bullet \, (s_r,t_r) & \mapsto  & (s_rt_r^3, t_r^4, 0, \dots, 0,
s_r^3t_r, s_r^4, s_r^2t_r^2).
\end{array}
\end{equation}\normalsize
We give the set of monomials the  graded $\rho$-weighted
lexicographic order, where $\rho$ is the one-parameter subgroup
with the weight vector
\begin{equation}\label{E:weight-even-g}
(1,0,2,1,0,2,3,4,2,1,0,2,3,4,2, \dots, 1,0,2,3,4,2).
\end{equation}
A Gr\"obner basis for $C_r$ is:
\[
\begin{array}{c}
x_0x_3, x_0x_4, \dots, x_0x_{3r-4}; \, x_1x_3, x_1x_4, \dots,
x_1x_{3r-4}; \,
x_2x_5, x_2x_6, \dots, x_2x_{3r-1};\\
x_0x_{3r-2} - x_{3r-3}x_{3r-1}, \,
x_{3r-1}^2 - x_1x_{3r-2}, x_{3r-1}^2 - x_0x_{3r-3}, \, \\
x_0^2 - x_1x_2 - x_1x_{3r-1}, \,
x_0x_{3r-1} - x_1x_{3r-3};\\
x_{6j+2}x_{6j+4} - x_{6j+3}^2 + x_{6j+4}x_{6j+5}, \quad j = 0, 1, \dots, \frac r2 - 1;\\
\end{array}
\] \normalsize and for $j = 1, 2, \dots, r-2$,
\[
\begin{array}{c}
x_{3j}x_{3j+2} - x_{3j+1}x_{3j+3}; \, x_{3j+2}^2 -
x_{3j+1}x_{3j+4}; \,
x_{3j+2}^2 - x_{3j}x_{3j+3}; \\
x_{3j}x_{3j+4} - x_{3j+2}x_{3j+3}; \, x_{3j+2}x_{3j+4} -
x_{3j+3}^2 + x_{3j+4}x_{3j+5}; \,
\\
x_{3j}x_{3j+5}, x_{3j}x_{3j+6}, \dots, x_{3j}x_{3r-1}; \,
x_{3j+1}x_{3j+5}, x_{3j}x_{3j+6}, \dots, x_{3j+1}x_{3j-1}; \\
x_{3j+2}x_{3j+5}, x_{3j+2}x_{3j+6}, \dots, x_{3j+2}x_{3j-1},
\end{array}
\]
together with the following degree three polynomials
\begin{equation}\label{E:deg3-gb}
\begin{array}{c}
x_1x_{3r-3}^2 - x_{3r-1}^3, x_{3r-3}^2x_{3r-1} - x_{3r-2}x_{3r-1}^2, \\
x_{3j+1}x_{3j+3}^2 - x_{3j+2}^3;  \quad j = 1, 2, \dots, r-2.
\end{array}
\end{equation}
The degree two initial monomials are:
\begin{equation}\label{E:1bead-deg2-1}
\begin{array}{c}
x_0^2, x_0x_3, x_0x_4, \dots, x_0x_{3r-1}; \,
x_1x_3, x_1x_4, \dots, x_1x_{3r-4}, x_1x_{3r-2}; \\
x_2x_4, x_2x_5, x_2x_6, \dots, x_2x_{3r-1};
\end{array}
\end{equation}\normalsize
and for $j = 1, 2, \dots, r-2$,
\begin{equation}\label{E:1bead-deg2-2}
\begin{array}{c}
x_{3j}x_{3j+2}, x_{3j}x_{3j+3}, \dots, x_{3j}x_{3r-1};
x_{3j+1}x_{3j+4}, x_{3j}x_{3j+5}, \dots, x_{3j+1}x_{3r-1}; \\
x_{3j+2}x_{3j+4}, x_{3j+2}x_{3j+5}, \dots, x_{3j+2}x_{3r-1}.
\end{array}
\end{equation}\normalsize
The sum of the weights of the monomials in (\ref{E:1bead-deg2-1})
is $27 r - 33$, whereas the monomials in (\ref{E:1bead-deg2-2})
contribute $18r^2 - 58r + 43$ to the total weight of the monomials
in the initial ideal.

The total weight of all degree two monomials is $18r^2 - 3r - 3$.
Hence the weights of all degree two monomials not in the initial
ideal  sum up to
\[
18r^2 - 3r - 3 - (18r^2-58r+43) - (27r-33) = 28r - 13
\]
On the other hand, the average weight $\frac{2 P(2) \sum r_i}{N+1}$
is $28 r - 14$. It follows from Proposition~\ref{prop:stabcrit2}
that $\mu([C_r]_2, \rho) = (-(28r-13) + 28r-14) =
-1$.

The degree three monomials  divisible
 by the ones
in the lists (\ref{E:1bead-deg2-1}), (\ref{E:1bead-deg2-2})
contribute $27r^3 + \frac{27}2r^2 - \frac{159}2 r^2 +26$ to the
total weight of the monomials in the initial ideal. On the other
hand, the monomials
\begin{equation}\label{E:1bead-deg3}
\begin{array}{l}
x_1x_{3r-3}^2, \, x_{3r-3}^2x_{3r-1}; \, x_{3j+2}^3, \quad j = 1, 2, \dots, r-2\\
\end{array}
\end{equation}
coming from the degree three Gr\"obner basis members
(\ref{E:deg3-gb}) contribute $6r + 2$. The sum of the weights of
all degree three monomials is $27r^3 + \frac{27}2 r^2 - \frac{15}2
r - 3$. Hence the total weight of the degree three monomials not
in the initial ideal is
\[
27r^3 + \frac{27}2 r^2 - \frac{15}2 r - 3 - \left( 27r^3 +
\frac{27}2r^2 - \frac{159}2 r^2 +26\right) - (6r + 2) = 66r - 31.
\]
On the other hand, the average weight is
\[
\frac{3 P(3)}{N+1} (6r-3) = 66r-33.
\]
By Proposition~\ref{prop:stabcrit2}, the Hilbert-Mumford index is
$\mu([C_r]_3, \rho) =  - (66r - 31) + 66r - 33)
= - 2.$ Since $\mu([C_r]_2, \rho) = 2 \mu([C_r]_3, \rho) < 0$, it
follows from Corollary~\ref{C:all-m} that $C_r$ is $m$-Hilbert
unstable for all $m \ge 2$.  Indeed, we find that
$$\mu([C_r]_m,\rho)=1-m$$
for each $m\ge 2$ and $\mu(Ch(C_r),\rho)=0$.
\end{proof}

\subsection{Basin of attraction: Closed rosary with a broken bead}

\begin{prop}\label{P:basin-cr-1br}
Let $C_r$ and $\rho$ be as in Proposition~\ref{P:hs-cr-1br}. Then
the basin of attraction $A_\rho([C]_m)$ parametrizes
closed elliptic chains
$(C', p, q)$ of length $(r+1)/2$ such
that $\iota(p) = \iota(q)$ (Figure~\ref{F:basin-c-ros-1broken}).

\begin{figure}[!t]
  \begin{center}
    \includegraphics[width=6in]{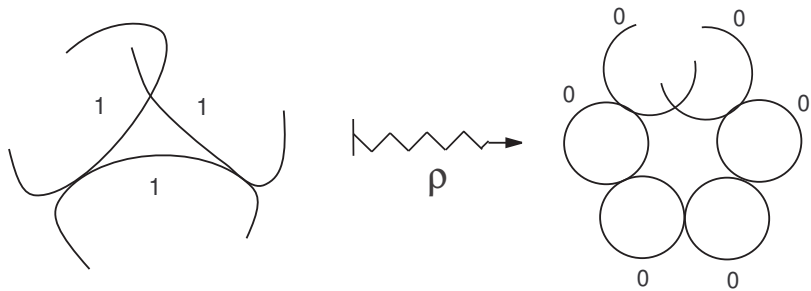}
  \end{center}
  \caption{Basin of attraction of a closed rosary with a broken bead}
  \label{F:basin-c-ros-1broken}
\end{figure}
\end{prop}
\begin{coro}\label{C:basin-cr-1br} A closed elliptic chain
is Hilbert unstable.
\end{coro}

\begin{proof} At the node, the local analytic equation is given by
\[
\frac{x_0}{x_2}\cdot \frac{x_3}{x_2} = 0
\]
and $\bG_m$ acts on the local parameters $x: = x_0/x_2$ and $y :=
x_3/x_2$ with weight $-1$. Hence $\bG_m$ acts on the local versal
deformation space (defined by $xy = c_0$) with weight $-2$. At the
adjacent tacnode,
$\bG_m$ acts on the tangent space to the two branches
with positive weights:  The tangent lines are traced by $x_3/x_4$ and
$x_6/x_4$. In fact, $\bG_m$ acts on the local versal deformation of
the tacnode (defined by $y^2 = x^4 + c_2 x^2 + c_1 x + c_0$) with a
positive weight vector $(2, 3, 4)$. Similar analysis reveals that
$\bG_m$ acts on the subsequent tacnode with a negative weight vector
$(-2, -3, -4)$.  Using the symmetry of the rosary, we can conclude
that $\bG_m$ acts on the local versal deformation space of the
tacnodes with weight vector alternating between $(2,3,4)$ and
$(-2,-3,-4)$.  The assertion now follows.
\end{proof}

\section{Proofs of semistability and applications}
\label{S:PSA}

Our main goal is
a complete description of all c-equivalence and h-equivalence
classes.  Throughout, each c-semistable curve $C$ is embedded
bicanonically (cf. Proposition~\ref{P:embedcstable})
$$C\hookrightarrow \bP^{3g-4},$$
and we consider the corresponding Chow points
$Ch(C) \in \Cg2$ and Hilbert points $[C]_m \in \Hg2, m \gg 0$.
To summarize:
\begin{itemize}
\item If $C$ is c-stable (resp.
h-stable), then the equivalence class of $C$ is trivial:
It coincides with the $\SL_{3g-3}$ orbit of $Ch(C)$ (resp. $[C]_m$ for $m\gg 0$).
\item  If $C$ is strictly c- or h-semistable its equivalence class is
nontrivial. We shall identify the  unique closed orbit curve and describe
all equivalent curves.
\item Since closed orbit curves are separated in a good quotient \cite[1.5]{Ses1}, we
have a complete classification of curves identified in the quotient spaces $\Hg2/\!\!/\SL_{3g-3}$ and $\Cg2/\!\!/\SL_{3g-3}$. 
\end{itemize}

\subsection{Elliptic bridges and their replacements}
\begin{figure}
  \begin{center}
    \includegraphics[width=3in]{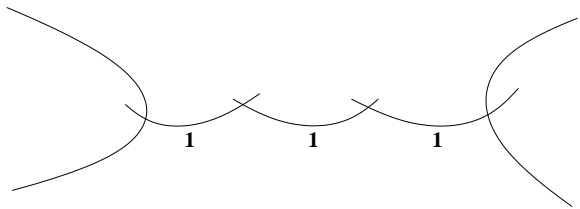}
  \end{center}

  \caption{A generic elliptic bridge of length three}
  \label{F:eb-length3}
\end{figure}

\begin{defn} An {\it elliptic bridge of length $k$} is a two-pointed curve $(C', p, q)$ (Figure~\ref{F:eb-length3}) such that
\begin{itemize}
\item $C' = E_1\cup_{a_1} \cdots \cup_{a_{k-1}} E_k$ consists of
connected genus-one
curves $E_1, \dots, E_k$ such that $E_i$ meets $E_{i+1}$ in a
node $a_i$, $i = 1, 2, \dots, k-1$;
\item $E_i \cap E_j \ne \emptyset$ if $|i-j|\ne 1$;
\item $p \in E_1$ and $q \in E_k$ are smooth points.
\end{itemize}
An ordinary elliptic bridge is an elliptic bridge of length one.
\end{defn}
\begin{figure}
  \begin{center}
    \includegraphics[width=2in]{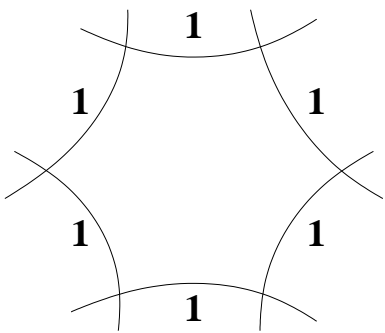}
  \end{center}

  \caption{A generic closed elliptic bridge of length six and genus seven}
  \label{F:closed-bridge}
\end{figure}

Let $C$ be a strictly c-semistable curve that is pseudostable,
i.e., $C$ has no tacnodes.  Let
$E_1,\ldots,E_N$ be the genus-one subcurves of $C$
arising as components of elliptic bridges.

\begin{figure}[!htb]
  \begin{center}
    \includegraphics{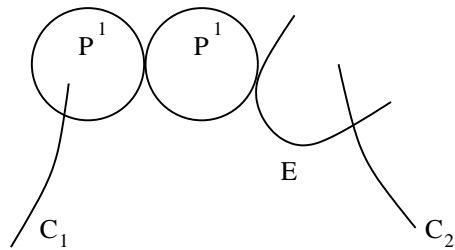}
  \end{center}
  \caption{A configuration corresponding to $C_1-\bP^1=\bP^1=E-C_2$}
  \label{F:graph}
\end{figure}

\begin{lemma} \label{L:possibilities}
Every c-semistable curve $C'$ admitting $C$ as a
pseudostable reduction can
be obtained from the following procedure:
\begin{enumerate}
\item{Fix a subset
$$\{E_i \}_{i \in I} \subset \{E_1,\ldots,E_N\}$$
of the genus-one subcurves arising in elliptic bridges.}
\item{Choose a subset of the nodes of $C$ lying on $\cup_{i\in I}E_i$
consisting of points of the following types:
\begin{itemize}
\item{If $E_i \cap E_{i'} \neq \emptyset$ for some distinct $i,i'\in I$ then
the node where they intersect must be included.}
\item{Nodes where the $E_i, i\in I$ meet other components
 may be included.}
\end{itemize}}
\item{
Replace each of these nodes by a smooth $\bP^1$ (for any
point of our subset) or by a chain of
two smooth $\bP^1$'s (only for points of the first type).  Precisely,
let $Z$ denote the curve obtained by normalizing our set of nodes
and then joining each pairs of glued points with a
$\bP^1$ or a chain of two $\bP^1$'s with one component meeting
each glued point.}
\item{
Let $E'_i$ denote the proper transform of
$E_i, i\in I$, which are pairwise disjoint in $Z$.
Replace each $E'_i$ with a
tacnode.  Precisely, write
$$D=Z \setminus \cup_{i\in I}E'_i$$
and consider a morphism
$$\nu:Z \ra C'$$
such that
\begin{itemize}
\item{
$\nu|D$ is an isomorphism and
$\nu|\overline{D} \ra C'$
is the normalization;}
\item{for $i\in I$, $\nu$ contracts
$E'_i$ to a tacnode of $C'$.}
\end{itemize}}
\end{enumerate}

The generic curve $C'$ produced by this procedure does not admit components
isomorphic to $\bP^1$ containing a node of $C'$.
We introduce $\bP^1$'s in Step (3) only to separate
two adjacent contracted elliptic components.
\end{lemma}
When enumerating the c-semistable curves,
it is convenient to use a graph that is similar to
the dual graph: We use a single line to denote a node and a double line to denote a
tacnode (e.g. Figure~\ref{F:graph}).

\begin{eg}
Let $C$ be the elliptic bridge of length one
$$C_1 - E - C_2.$$
The possible $Z$ are
$$C_1 - E - C_2, \quad C_1 - \bP^1 - E - C_2, \quad C_1 - E - \bP^1 - C_2 , \quad
C_1 - \bP^1 -  E - \bP^1 - C_2$$
and the possible $C'$ are
$$C_1-E-C_2, \quad C_1 = C_2, \quad C_1 - \bP^1 = C_2, \quad C_1 = \bP^1 - C_2 , \quad
C_1 - \bP^1 = \bP^1 - C_2.$$
The first two are the generic c-semistable configurations.

If $C$ is an elliptic bridge of length two
$$C_1 - E_1 - E_2 - C_2$$
then the possible $Z$ are
$$\begin{array}{cc}
C_1-E_1 - E_2 - C_2, \quad  C_1- \bP^1 - E_1 - E_2 - C_2,
\quad  C_1 - E_1 - \bP^1 - E_2 - C_2 \\
C_1 - E_1 - E_2 - \bP^1 - C_2,\quad
C_1 - \bP^1 - E_1 - \bP^1 - E_2 - C_2,  \\
C_1 - E_1 - \bP^1 - \bP^1 - E_2 - C_2, \quad
C_1 - E_1 - \bP^1 - E_2 - \bP^1 - C_2 , \\
C_1 - \bP^1 - E_1 - \bP^1 - \bP^1 - E_2 - C_2, \quad
C_1 - \bP^1 - E_1 - \bP^1 - E_2 - \bP^1 - C_2, \\
C_1 - E_1 - \bP^1 - \bP^1 - E_2 - \bP^1 - C_2 , \quad
C_1 - \bP^1 - E_1 - \bP^1 - \bP^1 - E_2 - \bP^1 - C_2.
\end{array}
$$
The generic c-semistable configurations are
$$
C_1-E_1-E_2-C_2, \quad C_1=E_2-C_2, \quad C_1 - E_1= C_2 ,
\quad C_1=\bP^1=C_2.
$$
\end{eg}
\begin{proof} (of Lemma~\ref{L:possibilities})
Our hypotheses give a flat
family
\[
\mcl C' \to B := \spec k[[t]] \tag{$\dag$}
\]
whose generic fibre $\mcl C'_\eta$ is smooth and the special fibre
$\mcl C'_0$ is $C'$.  Furthermore, after a base change
$$
\begin{array}{rcl}
B & \ra & B_1= \spec k[[t_1]] \\
t & \mapsto & t_1^s
\end{array}
$$
there exists a birational modification over $B_1$
$$\psi: \mcl C \dashrightarrow \mcl C' \times_B B_1$$
such that ${\mcl C}_0$ is $C$.
In other words, $\mcl C \to B_1$ is the {\em pseudostable reduction}
of $\mcl C' \ra B$;  we
replace each tacnode by an elliptic
bridge and contract any rational component that meets the rest of
the curve in fewer than three points.

Let $\mcl Z$ be the
normalization of the graph of $\psi$,
with $\pi_1$ and $\pi_2$ the
projections to $\mcl C$ and $\mcl C' \times_B B_1$ respectively:
\[
\xymatrix{ & \mcl Z \ar[dl]_-{\pi_1}\ar[rd]^-{\pi_2} & \\
\mcl C & & \mcl C' \times_B B_1 \\
}
\]
By \cite[4.4]{Sch}, $\mcl Z$ is flat over $B$ and $Z=\mcl Z_0$ is
reduced.  An argument similar to \cite[4.5-4.8]{Sch} yields
\begin{itemize}
\item  The exceptional locus of $\pi_2$ is a disjoint union
of connected genus-one subcurves
$$\sqcup_{i \in I} E_i' \subset Z$$
that arise as proper transforms of components of elliptic bridges
in $C$.  Each component is mapped to a tacnode of $C'$.
\item  The exceptional locus of $\pi_1$ is a union
of chains of rational curves of length one or two
$$\sqcup \bP^1 \sqcup (\bP^1 \cup \bP^1) \subset Z$$
that arise as proper transforms
of rational components of $C'$ meeting the rest of the curve in two points
(either two tacnodes or one node and one tacnode).
Each component is mapped
to a node of $C$ contained in an elliptic bridge.
\end{itemize}
This yields the schematic description for the possible combinatorial
types of $C'$.

We analyze the generic curves arising from our procedure.  Suppose
there is a component isomorphic to $\bP^1$ meeting the rest of the
curve in a node and a tacnode.  Corollary~\ref{C:versal} implies we can
smooth the node to get a c-semistable curve.  The smoothed curve
also arises from our procedure.
\end{proof}

\begin{rem} \label{R:subsets}
Lemma~\ref{L:possibilities} yields a bijection between subsets
$$\{E_i \}_{i\in I} \subset \{E_1,\ldots,E_N \}$$
and generic configurations of the locus of curves arising from
our procedure.  Indeed,
there is a unique generic configuration contracting the curves
$\{E_i \}_{i \in I }$.
\end{rem}

\begin{prop} \label{P:degenerate}
Let $C$ be strictly c-semistable without tacnodes and
$E_1,\ldots,E_N$ the genus-one subcurves of $C$
arising as components of elliptic bridges.  Let $C^{\star}$
be the curve obtained from $C$ by replacing each $E_i$
with an open rosary $(R_i,p_i,q_i)$ of length two.  Then there exists
a one-parameter subgroup
$$\rho:\bG_m \ra Aut(C^{\star})$$
such that $Ch(C) \in A_{\rho}(Ch(C^{\star}))$ and
$\mu(Ch(C^{\star}),\rho)=0$.

If $C'$ is another c-semistable curve with
pseudostable reduction $C$  then there exists
a one-parameter subgroup
$$\varrho':\bG_m \ra Aut(C^{\star})$$
such that $Ch(C') \in A_{\varrho'}(Ch(C^{\star}))$ and
$\mu(Ch(C^{\star}),\varrho')=0$.
\end{prop}
\begin{proof}
The assmption that $C$ is strictly c-semistable without tacnodes
ensures it contains an elliptic chain of length one, i.e., an elliptic
bridge.

The analysis of Proposition~\ref{P:rosary-aut} makes clear that
our description of $C^{\star}$ determines it uniquely up to isomorphism.
Furthermore, we have
$$Aut(C^{\star})^{\circ} \simeq \bG_m^N$$
with basis $\{\rho_1,\ldots,\rho_N\}$;  here
$\rho_i$ denotes the one-parameter subgroup
acting trivially on $R_j,j\neq i$
and with weight $1$ on the tangent spaces $T_{p_i}R_i$ and
$T_{q_i}R_i$.  As explained in \S~\ref{S:ba-o-ros},
it acts with negative weights on the versal deformation
space of the tacnode of $R_i$.

Consider the one-parameter subgroup
$$\rho=\prod_{i=1}^N \rho^{-1}_i,$$
which acts with positive weights on each of the tacnodes.
The basin of attraction analysis of Proposition~\ref{P:hs-echain}
shows that $A_{\rho}(Ch(C))$ parametrizes those curves
obtained by replacing
each open rosary of $C^{\star}$ with an elliptic bridge/chain of length one.
This includes our original curve $C$.

Now for any one-parameter subgroup
$$\varrho'=\prod_{i=1}^N \rho^{-e_i}_i,$$
we can compute
$$\mu(Ch(C^{\star}),\rho')=-\sum_{i=1}^N e_i \, \mu(Ch(C^{\star}),\rho_i)=0$$
using Corollary~\ref{C:cs-o-ros}.  In particular, we have
$$\mu(Ch(C^{\star}),\rho)=0.$$

Section~\ref{S:ba-o-ros} gives
the action of $\rho'$ on the versal deformations of the singularities of
$C^{\star}$.
It acts with weights $(2e_i,3e_i,4e_i)$ on the versal deformation
space of the tacnode on $R_i$.  At a node ($p_i$ or $q_i$)
lying on a single
open rosary $R_i$ of length two, it acts with weight $-e_i$.
For nodes on two open rosaries $R_i$ and $R_j$, it acts with
weight $-(e_i+e_j)$.

Restrict attention to one-parameter subgroups with weights
$e_i \neq 0$ for each $i,j=1,\ldots,N$.  These naturally
divide up into $2^N$ equivalence classes, depending on the signs
of the $e_i$.
Let $I\subset \{1,\ldots,N\}$ denote those indices with $e_i<0$.
Just as in the proof of Proposition~\ref{P:hs-echain},
the basin of attraction $A_{\varrho'}(Ch(C^{\star}))$ does not contain
smoothings of the tacnodes in $R_i, i\in I$ but does contain all smoothings of
the remaining tacnodes.  Choosing the negative $e_i$
suitably large in absolute value, we can assume each $-(e_i+e_j)>0$,
so the nodes where two rosaries meet are smoothed provided
at least one of the adjacent tacnodes is {\em not} smoothed.

Thus $A_{\varrho'}(Ch(C^{\star}))$ consists of the c-semistable
curves obtained by smoothing all the tacnodes {\em not} indexed
by $I$, as well as the nodes on the rosaries containing one
of the remaining tacnodes (indexed by $I$).
The generic member of the basin equals
the generic configuration indexed by $I$, as described in
Remark~\ref{R:subsets}.
It follows that {\em each}
curve $C'$ enumerated in Lemma~\ref{L:possibilities} appears in the
the basin of attraction of $Ch(C^{\star})$ for a suitable
one-parameter subgroup $\varrho'$.
\end{proof}

\subsection{Chow semistability of c-semistable curves}\label{SS:CSBC}
Here we prove that bicanonical c-semistable curves are Chow
semistable. By Theorem~\ref{T:partial}, it suffices to consider
curves that are not c-stable.

Let $C'$ denote a strictly c-semistable curve, with tacnodes and/or
elliptic bridges.
Assume that $C'$ is Chow unstable and let
\[
\mcl C' \to B := \spec k[[t]]
\]
be a smoothing.  Let $C''$ be a
Chow semistable reduction of this family
(see Theorem~\ref{T:ss-replacement}) and $C$ the pseudostable reduction.

Reversing the steps outlined in Lemma~\ref{L:possibilities}, we see
that $C$ is obtained by replacing each tacnode of $C'$ (or $C''$)
with an elliptic bridge and then pseudo-stabilizing.  Let
$C^{\star}$ denote the curve obtained from $C$ in
Proposition~\ref{P:degenerate}, which guarantees that $Ch(C')$ and
$Ch(C'')$ are contained in basins of attraction
$A_{\rho'}(Ch(C^{\star}))$ and $A_{\rho''}(Ch(C^{\star}))$
respectively.  Moreover, since
$$\mu(Ch(C^{\star}),\rho')=
\mu(Ch(C^{\star}),\rho'')=0$$ Lemma~\ref{L:flat limit} implies that
$C'$ (resp. $C''$)  is Chow semistable iff $C$ is Chow semistable.
This contradicts our assumption that $C'$ is Chow unstable.

Next, we give a characterizaztion of the  closed orbit curves
in c-equivalence classes of strictly semistable curves:

\begin{prop}\label{P:c-minorbit} A strictly c-semistable curve has a closed orbit if and only if
\begin{itemize}
\item each tacnode is contained in an open rosary;
\item each open rosary has length two; and
\item there are no elliptic bridges other than length two rosaries.
\end{itemize}
\end{prop}

Since each length-two rosary has one tacnode and
contributes a $\bG_m$-factor to
$Aut(C)$, we have:
\begin{coro} If $C$ is a strictly c-semistable curve with closed orbit then
\[
Aut(C)^{\circ} \simeq \bG_m^{\tau}.
\]
where $\tau$ is the number of tacnodes and the superscript $\circ$ denotes the connected component of the identity.
\end{coro}

\begin{proof}[Proof of Proposition~\ref{P:c-minorbit}]
Assume that $C'$ is a strictly semistable curve with closed orbit.
Let $C$ be a pseudostable reduction and $C^{\star}$ the curve
specified in Proposition~\ref{P:degenerate}, so the Chow point
of $C'$ is in the basin of attraction of the Chow point of $C^{\star}$.
Since $C^{\star}$ is Chow semistable, we conclude that $C'=C^{\star}$.

Conversely, suppose $C'$ is a curve satisfying the three conditions
of Proposition~\ref{P:c-minorbit}.  Again, let $C$ be a pseudostable
reduction of $C'$ and $C^{\star}$ the curve obtained in Proposition~\ref{P:degenerate},
so that $C'$ is in the basin of attraction of $C^{\star}$ for some one-parameter
subgroup $\rho'$.   Note that $C^{\star}$ also satisfies the conditions of
Proposition~\ref{P:c-minorbit}.  The basin of attraction
analysis in Section~\ref{S:ba-o-ros} implies that any nontrivial deformation
of $C^{\star}$ in $A_{\rho'}(Ch(C^{\star}))$ induces a nontrivial deformation
of at least one of the singularities
of $C^{\star}$ sitting in an open rosary.

There are three cases to consider:  First, we could deform the tacnode on one of the
rosaries $R_i$.  However, then the rosary $R_i$ deforms to an elliptic bridge in $C'$
that is not a length two rosary, which yields a contradiction.  Therefore, we may
assume that none of the tacnodes in $C^{\star}$ is deformed in $C'$.
Second, we could smooth a node where length two rosaries meet.  However, this
would yield a rosary in $C'$ of length $>2$.  Finally, we could smooth a node where
a length two rosary $R_i$ meets a component not contained in an rosary.  However,
the tacnode of $R_i$ then deforms to a tacnode of $C'$ not on any length two rosary.
\end{proof}

\begin{figure}[!t]
  \begin{center}
    \includegraphics[width=6.5in]{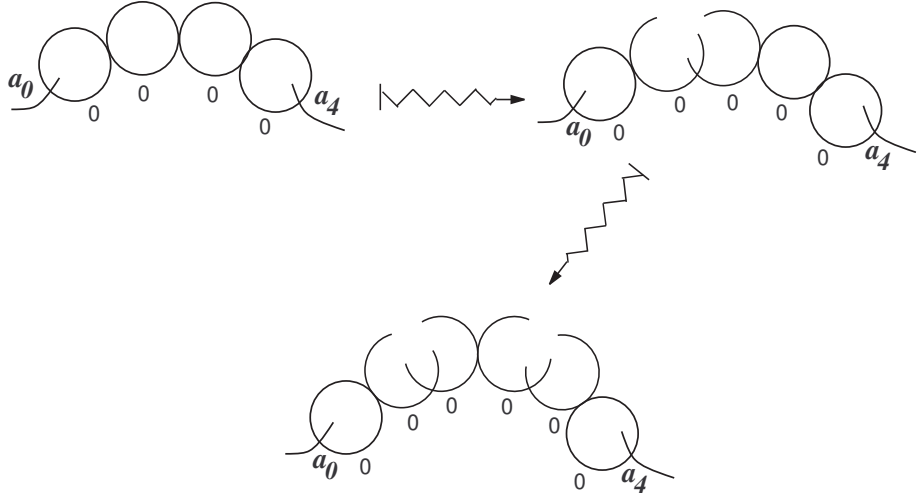}
  \end{center}
  \caption{Degeneration to the c-semistable closed orbit curve}
  \label{F:min-o-ros1}
\end{figure}


\subsection{ Hilbert semistability of h-semistable curves}\label{SS:HS}
Suppose that $C$ is a h-semistable bicanoncal curve.  By definition it is
also c-semistable and thus Chow-semistable by the analysis of
Section~\ref{SS:CSBC}.  Of course, strictly Chow-semistable points
can be Hilbert unstable, and we classify these in two steps.
First, we enumerate the curves $C_0$ with strictly semistable
Chow point such that there exists a one-parameter subgroup
$\rho: \bG_m \hookrightarrow Aut(C_0)$
destabilizing the Hilbert point of $C_0$, i.e., with $\mu([C_0]_m,\rho)<0$
for $m \gg 0$.  Second, we list the curves that are in the basins of attraction
$A_{\rho}([C_0]_m)$, which are also guaranteed to be Hilbert unstable by the
Hilbert-Mumford one-parameter subgroup criterion.
Proposition~\ref{P:destabilize} shows these
are {\em all} the Hilbert unstable curves:
If $C$ is a c-semistable bicanonical curve that is
Hilbert unstable then its Chow point is contained in the basin of attraction
$A_\rho(Ch(C_0))$ of a Chow semistable curve $C_0$ with closed orbit such that
$\mu([C_0]_m, \rho) < 0.$

If the genus is odd and
$C_0$ is a closed rosary (without broken beads) then
$C_0$ is Hilbert semistable
with respect to any 1-PS coming from $Aut(C_0)$
(Proposition~\ref{P:hss-c-ros}).

Suppose that $C_0$ has open rosaries $S_1, \dots, S_\ell$.
Each contributes $\bG_m$ to the automorphism group of $C_0$ and
$Aut(C_0)^{\circ} \simeq
\bG_m^{\times \ell}$.
Let $p_i, q_i$ denote the nodes in the
intersection $S_i\cap \overline{C_0-S_i}$.
The automorphism coming from $S_i$ gives rise to a one-parameter subgroup
\[
\rho_i : \bG_m \stackrel{\simeq}{\longrightarrow} \{1\} \times \cdots \times \underbrace{\bG_m}_{\mbox{$i$th}}\times \cdots \times \{1\} \inj \bG_m^{\times \ell} \simeq G_{Ch(C_0)}^{\circ}
\]
where the second $\bG_m$ means the $i$th copy in the product $\bG_m^{\times \ell}$ and $G_{Ch(C_0)}$ is the stabilizer group.
We assume that $S_1, \dots, S_k$ are open rosaries
of even length and $S_{k+1}, \dots, S_\ell$ are of odd length.
For $i \le k$, the weights of $\rho_i$ on the versal deformation
spaces of $p_i$ and $q_i$ have the same sign
(see \S\ref{S:ba-o-ros}).
We normalize $\rho_i$ so that this weight is negative.

Given our one-parameter subgroups $\rho:\bG_m \ra Aut(C_0)^{\circ}$ with
$\mu([C_0]_m,\rho)<0$, we can expand
\[
\rho = \prod_{i=1}^k \rho_i^{a_i} \times \prod_{i=k+1}^\ell \rho_i^{b_i}, \quad a_i,b_i \in \bZ
\]
so that
\[
\mu([C_0]_m, \rho)  
 = \sum_{i=1}^k  a_i \, \mu([C_0]_m, \rho_i) + \sum_{i=k+1}^\ell b_i \, \mu([C_0]_m, \rho_i)
<0.
\]
We have already computed these terms:  Proposition~\ref{P:hs-o-ros} implies
that $\mu([C_0]_m,\rho_i)=0$ for $i=k+1,\ldots,\ell$ and
$\mu([C_0]_m,\rho_i)=1-m$ for $i=1,\ldots,k$.  Thus in order for the sum to be
negative, we must have $a_i>0$ for some $i=1,\ldots,k$.  In particular, there
is at least one rosary of even length.   Proposition~\ref{P:hs-echain} implies
that the basin of attraction $A_{\rho}([C_0]_m)$ contains curves with
elliptic chains, which are not h-semistable.

We are left with the case of a closed rosary $C_r$ of even genus with one broken bead.
There is a unique one-parameter subgroup $\rho$ of the automorphism group,
and we choose the sign so that it destabilizes $C_r$ (cf.
Proposition~\ref{P:hs-cr-1br}).  The basin of attraction analysis in
Proposition~\ref{P:basin-cr-1br} again shows that the curves
with unstable Hilbert points admit elliptic chains.

Thus curves with unstable Hilbert points are not h-semistable, which completes
our proof that h-semistable curves are Hilbert semistable.

\

We shall now prove that if $C$ is h-stable then it is Hilbert
stable. If $C$ is Hilbert strictly semistable, then it belongs to
a basin of attraction $A_\rho([C_0]_m)$ where $C_0$ is a Hilbert
semistable curve  with infinite automorphisms and $\rho$ is a 1-PS
coming from $Aut(C_0)$. By Corollary~\ref{C:hsinf}, $C_0$ admits
an open rosary of odd length $\ge 3$ or is a closed rosary of even length $\ge 4$. But
we showed in Propositions \ref{P:hs-echain} and
\ref{P:gen-ss-c} that any curve in the basin of such $C_0$ has a
weak elliptic chain and hence is not h-stable.

Finally, we characterize the closed orbits of strictly h-semistable curves.
These do not admit elliptic chains, and in particular,
do not admit open rosaries of even length
(see Remark~\ref{rem:rosary}).

\begin{prop} \label{P:h-minorbit} A strictly h-semistable curve has a closed orbit if and only
if
\begin{itemize}
\item it is a closed rosary of odd genus; or
\item each weak elliptic chain is contained in a chain of open rosaries of length three.
\end{itemize}
\end{prop}
Since each length three open rosary has two tacnodes
and contributes $\bG_m$ many automorphisms
to $Aut(C)$,
\begin{coro} If $C$ is a strictly
h-semistable curve with closed orbit then
\[
Aut(C)^{\circ} \simeq \bG_m^{\tau/2}
\]
where $\tau$ is the number of tacnodes.
\end{coro}
\begin{proof}[Proof of Proposition~\ref{P:h-minorbit}]
Suppose $C'$ is strictly h-semistable.  We shall show that there exists a
curve $C^{\ast}$ satisfying the conditions of Proposition~\ref{P:h-minorbit}
and a one-parameter subgroup $\rho'$ of $Aut(C^{\ast})$
such that $[C']_m \in A_{\rho'}([C^{\ast}]_m)$ for $m \gg 0$ and
$\mu([C^{\ast}]_m,\rho')=0$.

Assume first that $C'$ is a closed weak
elliptic chain with $r$ components, with arithmetic genus
$2r+1$.  Let $C^{\ast}$ denote a closed rosary with beads $L_1,\ldots,L_{2r}$
and tacnodes $a_1,\ldots,a_{2r}$.
Proposition~\ref{P:rosary-aut} implies
$Aut^{\circ}(C^{\ast})\simeq \bG_m$, generated by a one-parameter subgroup $\rho$
acting on the versal deformation spaces of the $a_{2j}$
with positive weights and the $a_{2j-1}$ with negative weights.
Proposition~\ref{P:gen-ss-c} implies $A_{\rho}([C^{\ast}]_m)$ contains the closed weak elliptic
chains of length $r$.

Now assume that $C'$ is not a closed weak elliptic chain but contains maximal
closed weak elliptic chains $C''_1,\ldots,C''_s$ of lengths $\ell_1,\ldots,\ell_s$.
Let $p_j$ (resp. $q_j$) denote the node (resp. tacnode) where $C''_j$ meets
the rest of the curve.
Let $C^{\ast}$ be the curve obtained from $C'$
by replacing each $C_j''$ with a chain of $\ell_j$ open rosaries of length three.
Precisely, write
$$D=C'\setminus (\bigcup_{j=1}^s C''_j \setminus \{p_j,q_j \})$$
and let $S_j,j=1,\ldots,s$ denote a chain of $\ell_j$ open rosaries
of length three joined end-to-end. Then $C^{\ast}$ is obtained by
gluing $S_j$ to $D$ via nodes at $p_j$ and $q_j$.
One special case requires further explanation:  If $C'$ admits an irreducible
component $\simeq \bP^1$ meeting the rest of $C'$ at two points $q_i$ and $q_j$ then
we contract this component in $C^{\ast}$.

\begin{eg} There are examples where the construction of $C^{\ast}$ involves components being contracted.
Let $C_1$ and $C_2$ be smooth and connected of genus $\ge 2$ and let $E_1$ and $E_2$ be elliptic.
Consider the curve $C'$
$$C_1-_{p_1}E_1=_{q_1}\bP^1=_{q_2}E_2-_{p_2}C_2.$$
Replacing the weak elliptic chains with rosaries
of length three yields
$$C_1 -_{p_1} \bP^1 = \bP^1 = \bP^1 -_{q_1} \bP^1 -_{q_2} \bP^1 = \bP^1 = \bP^1 -_{p_2} C_2$$
which is not h-semistable.  Contracting the middle $\bP^1$, we obtain $C^{\ast}$:
$$C_1 -_{p_1} \bP^1 = \bP^1 = \bP^1 - \bP^1 = \bP^1 = \bP^1 -_{p_2} C_2.$$

There are examples where $D$ fails to be pure-dimensional.  Start with the curve $C'$
$$C_1=_{q_1}E_1-_{p_1=p_2}E_2=_{q_2}C_2$$
with the $C_i$ and $E_i$ as above.  Then $C^{\ast}$ is equal to
$$C_1 -_{q_1} \bP^1 = \bP^1 = \bP^1-_{p_1=p_2}\bP^1 = \bP^1 = \bP^1 -_{q_2}C_2.$$
\end{eg}

We return to our proof:  The curve $C^{\ast}$ has
$$Aut(C^{\ast})^{\circ} \simeq  \bG_m^N, \quad N=\sum_{j=1}^s \ell_j.$$
Essentially repeating the argument of Proposition~\ref{P:degenerate}, using the
one-parameter subgroup analysis of Proposition~\ref{P:hs-o-ros} and the
basin-of-attraction analysis of Proposition~\ref{P:hs-echain}
(or Proposition~\ref{P:hs-cr-1br} and \ref{P:basin-cr-1br} in the generate case),
we obtain a one-parameter subgroup  $\rho'$ in the automorphism group
such that $[C']_m \in A_{\rho'}([C^{\ast}]_m)$ for $m \gg 0$ and
$\mu([C^{\ast}]_m,\rho')=0$.

We now show that the curves enumerated in Proposition~\ref{P:h-minorbit} all have
closed orbits.  Due to \cite[Theorem 1.4]{Kempf},
it suffices to show none of these are contained in the
basin of attraction of any other.  Suppose that $C^{\ast}_1$ and $C^{\ast}_2$
are such that
$$[C^{\ast}_2]_m \in A_{\rho}([C^{\ast}_1]_m)$$
for some one-parameter subgroup $\rho$ of $Aut(C^{\ast}_1)^{\circ}$.  A nontrivial
deformation of $C^{\ast}_1$ necessarily deforms one of the singularities of $C^{\ast}_1$.
If the singularity is a tacnode on a length-three open rosary, the resulting deformation
admits a weak elliptic chain that is not contained in a chain of length-three rosaries.
If the singularity is a node where two length-three open rosaries meet then the deformation
admits a weak elliptic chain {\em not} contained in a chain of length-three open rosaries.

However, there is one case that requires special care:  Suppose that $C^{\ast}_1$ is a closed
chain of $r$ rosaries of length three
$$R_1 -_{p_{12}} R_2 -_{p_{23}} \cdots -_{p_{r-2\, r-1}} R_{r-1} -_{p_{r-1r}} R_r -_{p_{r1}} R_1
$$
where $R_1$ and $R_r$ meet at a node $p_{r1}$;
this has arithmetic genus $2r+1$.
Let $C^{\ast}_2$ denote a closed rosary of genus $2r+1$, which is a deformation of
$C^{\ast}_1$.  We need to insure that
\begin{equation} \label{eqn:tricky}
[C^{\ast}_2]_m \not \in A_{\rho}([C^{\ast}_1]_m)
\end{equation}
for any one-parameter subgroup $\rho$ of $Aut([C^{\ast}_1]_m)^{\circ}$.
We can express
$$\rho = \prod_{j=1}^r \rho_j^{e_j}$$
where $\rho_j$ acts trivially except on $R_j$ and has weights $+1$ and $-1$
on $T_{p_{j-1j}}R_j$ and $T_{p_{jj+1}}R_j$.  (Here $\rho_r$ acts with
weights $+1$ and $-1$ on $T_{p_{r-1r}}R_r$ and $T_{p_{r,1}}R_r$.)
However, assuming $\rho$ is nontrivial, one of the following differences
$$e_1-e_2,\ldots, e_r-e_1$$
is necessarily negative;  for simplicity, assume $e_1-e_2<0$.  It follows
that $\rho$ acts with negative weight on the versal deformation of
the node $p_{12}$, thus deformations in $A_{\rho}([C^{\ast}_1]_m)$ cannot
smooth $p_{12}$.  We conclude that deformations in the basin
of attraction of $C^{\ast}_1$ cannot smooth each node,
which yields  (\ref{eqn:tricky})
\end{proof}

\end{document}